\documentclass{amsart}
\usepackage[english]{babel}
\usepackage{graphicx}
\usepackage[hidelinks]{hyperref} 
\usepackage{amssymb}
\usepackage{amsmath}
\usepackage[foot]{amsaddr}
\usepackage{amsfonts}
\usepackage{xcolor}
\usepackage{mathtools}
\usepackage{algpseudocode}
\usepackage{algorithm}
\usepackage[top=1in, bottom=1.25in, left=1.25in, right=1.25in]{geometry}
\usepackage{caption}
\usepackage{subcaption}

\graphicspath{{images/}}
\input{packages.sty}

\begin{document}
\title[Reduced order models for the buckling of hyperelastic beams]{Reduced order models for the buckling of\\ hyperelastic beams}

\author{Federico Pichi$^1$, Gianluigi Rozza$^2$}

\address{$^1$ Chair of Computational Mathematics and Simulation Science, École Polytechnique Fédérale de Lausanne, 1015 Lausanne, Switzerland}
\address{$^2$ mathLab, Mathematics Area, SISSA, via Bonomea 265, I-34136 Trieste, Italy}

\begin{abstract}
In this paper, we discuss reduced order modelling approaches to bifurcating systems arising from continuum mechanics benchmarks.
The investigation of the beam's deflection is a relevant topic of investigation with fundamental implications on their design for structural analysis and health.
When the beams are exposed to external forces, their equilibrium state can undergo to a sudden variation. This happens when a compression, acting along the axial boundaries, exceeds a certain critical value. 
Linear elasticity models are not complex enough to capture the so-called beam's buckling, and nonlinear constitutive relations, as the hyperelastic laws, are required to investigate this behavior, whose mathematical counterpart is represented by bifurcating phenomena.
The numerical analysis of the bifurcating modes and the post-buckling behavior, is usually unaffordable by means of standard high-fidelity techniques such (as the Finite Element method) and the efficiency of Reduced Order Models (ROMs), e.g.\ based on Proper Orthogonal Decomposition (POD), are necessary to obtain consistent speed-up in the reconstruction of the bifurcation diagram.
The aim of this work is to provide insights regarding the application of POD-based ROMs for buckling phenomena occurring for 2-D and 3-D beams governed by different constitutive relations. The benchmarks will involve multi-parametric settings with geometrically parametrized domains, where the buckling's location depends on the material and geometrical properties induced by the parameter. Finally, we exploit the acquired notions from these toy problems, to simulate a real case scenario coming from the Norwegian petroleum industry.
\end{abstract}

\maketitle

\section{Introduction and motivation}
A huge variety of physical phenomena can be modelled by the so-called Partial Differential Equations (PDEs), which exploits the notion of differential operators to mathematically express the relationships between spatial and temporal unknowns of the field one wants to investigate. 

The generalizability of such a description depends intrinsically on its capability at reproducing the real physical behavior of the experiment at hand. This can be achieved by introducing physical or geometrical parameters characterizing the model's properties.
Indeed, it is intuitive that an object made of rubber behaves in a drastically different way with respect to a steel one.

In this work we are interested on a special class of parametrized PDEs aiming at modelling the equilibrium state(s) of continuum mechanics problems related to the bucking of beams. Thus, we focus the analysis on steady-state systems subject to nonlinear constitutive laws.
These are mathematical relations that link the physical and the mechanical properties of the material, determining the way in which it responds to external forces.

The simplest model in the theory of elasticity \cite{ciarlet1988three,Timoshenko,ciarlet1997mathematical} is the famous Hooke's law for elastic object, assuming the stress and strain proportionally related through some elasticity constants.
However, given its linear nature, this model fails to describe accurately many real-life applications, where the material is subject to large deformation, plasticity or viscous behavior.

It is thus crucial to introduce nonlinear terms in the models, which are able to capture complex evolution of the mechanical problem. Unfortunately, this complicates the numerical setting, both from the theoretical and computational point of view. On one side, already at the mathematical level, the investigation of nonlinear system could be an ill-posed problem, admitting several stable states w.r.t.\ the same physical configuration. This is related to the well-known theory of bifurcation \cite{seydel2009practical,kuznetsov2004elements,kielhofer2006bifurcation,Brezzi1980}, i.e.\ the lack of uniqueness of the PDE's solution, whose goal is the identification of the critical points and the reconstruction of the bifurcation diagram.
From the computational point of view, one needs to face many difficulties:
\begin{itemize}
    \item[(i)] the detection of the bifurcation point, as the critical value from which one have multiple coexisting states; 
    \item[(ii)]  the design of a stable and robust nonlinear solver, capable of dealing with the system's ill-posedness in a neighborhood of the bifurcation; 
    \item[(iii)] the discovery of all the possible solution branches originating from the model, and thus the implementation of a continuation method through which obtaining the bifurcation diagram.
\end{itemize}

Even assuming to being able to deal with all these issues, obtaining  accurate enough numerical simulations is usually an unbearable task. Indeed, in order to capture all the effects of the nonlinear terms, and thus the bifurcating behavior, one needs to consider a fine mesh discretization which makes it impossible to fully analyze the evolution of the branches in the whole parameter space. 

Although computational resources are expanding, Reduced Order Models (ROMs) \cite{BennerModelOrderReduction2020,benner2017model} are fundamental for this kind of applications. ROMs constitute a class of methodologies that aims at obtaining efficient and reliable evaluations of the unknowns in a real-time context. The projection-based models we analyze in this work, based on solid mathematical foundations, exploit the construction of a new basis, e.g.\ by means of Proper Orthogonal Decomposition (POD) or Greedy algorithms \cite{hesthavenCertifiedReducedBasis2015,quarteroniReducedBasisMethods2016,patera07:book}, which spans a low-dimensional space capturing most of the information.
A few works have analyzed the ROMs in the context of elasticity models \cite{RozzaReducedBasisApproximation2008,NiroomandiModelOrderReduction2013,NiroomandiModelOrderReduction2010,ZanonReducedBasisMethod2013,veroy03:_phd_thesis,zanon:_phd_thesis,Huynh2018,pichi_phd}. 

In order to be able to obtain efficient evaluation of the PDE's solution, an essential component is the so-called offline-online paradigm. It is indeed crucial to decouple the online phase, in which we want to perform the real-time many-query task, from the high dimensionality of the original system, which burdens the computational effort.
Reduced order models have been utilized in several bifurcating contexts in computational mechanics, from structural \cite{PichiReducedBasisApproaches2019,PichiChapterReducedBasis2022,NOOR198367} to quantum mechanics \cite{PichiReducedOrderModeling2020}, with a special focus on fluid dynamics applications \cite{pitton_quaini_2017,PichiArtificialNeuralNetwork2023,PintoreEfficientComputationBifurcation2021,HessLocalizedReducedorderModeling2019,Maday:RB2,PLA2015162,Terragni:2012}, and their extension toward more complex models \cite{PichiDrivingBifurcatingParametrized2022a,TonicelloNonintrusiveReducedOrder2022,KhamlichModelOrderReduction2022,HessReducedBasisModel2020}.
However, when dealing with bifurcating phenomena, this decoupling is compromised by the presence of nonlinear terms. Indeed, this prevents the assembly of all the costly high-fidelity quantities during the offline stage. For this reason, Empirical Interpolation methodologies such as EIM \cite{barrault04:_empir_inter_method} and its variant DEIM \cite{Chaturantabut2010} have been developed in the last decades. These are hyper-reduction approaches, or affine-recovery techniques, which aim at interpolating the nonlinear or non-affine terms, to be efficiently pre-computed and efficiently assembled for each new instance of the parameter.  
We remark that all the numerical results have been obtained using the RBniCS \cite{rbnics} library.

The manuscript is divided as follows: Section 2 concerns a brief theoretical and methodological introduction to bifurcation theory and its numerical approximation, both at the high-fidelity and the reduced-order level. In Section 3 we review the class of hyperelastic models to study the buckling of beams, through Kirchhoff--Saint Venant and neo-Hookean materials.
Several numerical results are presented in Section 4, to investigate the bifurcation w.r.t.\ several test cases, featuring different boundary conditions, 2-D and 3-D settings and multi-parametric contexts with varying geometrical characteristics.
\section{Numerical approximation of bifurcating phenomena}
In this section, we want to provide a brief description of bifurcating behaviors and their mathematical investigation.
Physical phenomena can exhibit either gradual or sudden changes, which can be linked to quantitative or qualitative changes, respectively. 
As an example, we can think about the problem of beam under compression.
When the load applied to the beam, denoted with $\bmu$, is sufficiently small, we observe a slight deformation that corresponds to a quantitative change. Of course, this deformation depends on the magnitude of the load, and on the specific material properties of the beam. Hence, we expect that a small perturbation of $\bmu$ leads to a configuration with the same qualitative feature.
At the same time, we empirically know that a critical value for the load $\bmu^*$ exists, after which the model does not conserve the original qualitative behavior. Practically speaking, if the load is big enough, the beam can not sustain the compression, and it buckles exhibiting a  qualitative change.
Therefore, we will denote the qualitative changes that are not stable under small perturbations of the parameter $\bmu$ as \textit{bifurcation phenomena}, and the points $\bmu^*$ in which they occur as \textit{bifurcation points} \cite{seydel2009practical,kielhofer2006bifurcation,kuznetsov2004elements}. 

Here, we are interested in nonlinear and parametric PDEs, in the context of continuum mechanics, that model this kind of phenomena, and for which the solution at a given parameter $\bmu$ exist, but may not be unique \cite{Brezzi1980,caloz1997numerical}. Indeed, the local (w.r.t.\ the parameter) well-posedness of the problem \eqref{eq:strong_form} relies on the local spectral property of its Frech\'et derivative, and if it fails to fulfill the inf-sup stability assumptions \cite{Brezzi1980,pichi_phd}, the behavior of the model has to be investigated more carefully. 

\subsection{Mathematical framework}\label{se:mat}

Let $\Pa \subset \R^P$ be the parameter space, a closed and bounded subset of the space $\R^P$ with $P \geq 1$. A parameter $\bmu$ is a point in the space $\Pa$ that can contain physical or geometrical information about the system. Then, we introduce the (reference) domain $\Om \subset \R^d$, where $d = 2, 3$ is the spatial dimension.

We define the Hilbert space $\X \eqdot \X(\Omega)$ and its dual space $\X'$, as the space of linear and continuous functional over $\X$. The Hilbert space is equipped with the norm $\norm{V}_\X = (V, V)_\X^{1/2}$ for any $V  \in \X$, induced by its inner product $(\cdot, \cdot)_\X$.
Within this setting, we denote the duality pairing between $\X'$ and $\X$ by means of $$\langle G, V\rangle \eqdot \langle G, V\rangle_{{\X'}\X} \quad \forall \, G \in \X', V \in \X .$$
We denote with $G: \X \times \Pa \to \X'$ the parametrized mapping which represents the nonlinear PDE, and thus the model under investigation.
The strong form of the parametric PDE problem reads as: given $\bmu \in \Pa$, find $X(\bmu) \in \X$ such that 
\begin{equation}
\label{eq:strong_form}
G(X(\bmu); \bmu) = 0 \quad \in \X' .
\end{equation}
To study the solutions of the abstract equation above, a fundamental step is the analysis of its variational formulation, which will play a key role in the numerical approximation of the problem.
Thus, we introduce the parametrized variational form $g(\cdot, \cdot; \bmu): \X \times \X \to \mathbb{R}$ as 
\begin{equation}
\label{eq:var_form}
g(X, Y; \bmu) = \langle G(X; \bmu), Y\rangle \quad \forall \, X, Y \in \X ,
\end{equation}
and the weak formulation of equation \eqref{eq:strong_form} reads as: given $\bmu \in \Pa$, find $X(\bmu) \in \X$ such that 
\begin{equation}
\label{eq:weak_form}
g(X(\bmu), Y; \bmu) = 0 \quad \forall \, Y \in \X .
\end{equation}

In order to have a complete description of the problem formulation in \eqref{eq:strong_form}, we need some suitable \textit{boundary conditions} (BCs) on the domain $\Omega$. Here, we will cover both the cases of homogeneous/inhomogeneous Dirichlet and Neumann conditions, depending on the test case at hand.
Moreover, these BCs will be automatically embedded into \eqref{eq:weak_form} by an appropriate choice of the space where to seek the solution.

Having introduced the mathematical formulation of the problem, a starting point for the analysis is the investigation of its well-posedness. A PDE is said to be a \textit{well-posed} problem if it is characterized by the existence and the uniqueness of the solution. As we discussed, these assumptions are not verified when dealing with more realistic models, such as the one involving bifurcating phenomena.

\subsection{High fidelity approximation}\label{se:hf}
As we mentioned in the introduction, the first level of discretization is represented by high-fidelity techniques.
These can potentially involve a high number of degrees of freedom, usually causing a computational complexity which is unaffordable in a many-query scenario. 

In order to be able to reconstruct the bifurcation diagram, it is necessary to combine three well-studied methodologies: the Galerkin Finite Element method, the Newton-Kantorovich method and a Continuation method, respectively to discretize, linearize and continue each branch.

We consider the finite dimensional space $\X_\mathcal{N} \subset \X$, and performing a Galerkin projection of the weak formulation \eqref{eq:weak_form}, the Galerkin FE method reads as: given $\bmu \in \Pa$, we seek $X_\mathcal{N} \eqdot X_\mathcal{N}(\bmu) \in \X_\mathcal{N}$ that satisfies 
\begin{equation}
\label{eq:weakgal}
\langle G(X_\mathcal{N}; \bmu), Y_\mathcal{N} \rangle = g(X_\mathcal{N}, Y_\mathcal{N} ; \bmu) = 0\ , \quad \forall \, Y_\mathcal{N} \in \X_\mathcal{N}. 
\end{equation}

Since the discretized weak problem \eqref{eq:weakgal} still involves a nonlinear structure, we have to deal with a method that linearizes it, in order to be effectively treated and solved by means of the Galerkin FE method.
The nonlinear solver chosen to linearize the weak formulation in \eqref{eq:weakgal} is the well-known Newton-Kantorovich method \cite{ciarlet2013linear,quarteroniReducedBasisMethods2016}, which reads as follows: given $\bmu \in \Pa$ and chosen an initial guess $X_{\N}^0(\bmu) \in \X_{\N}$, for every $k = 0, 1, \dots$ we seek the variation $\delta X_{\N} \in \X_{\N}$ such that
\begin{equation}
\label{eq:weaknewton}
\mathrm{d}g[X_{\N}^k(\bmu)](\delta X_{\N}, Y_{\N}; \bmu) =  g(X_{\N}^k(\bmu),Y_{\N} ; \bmu) \ , \quad \forall \, Y_{\N} \in \X_{\N} ,
\end{equation}
then we update the solution at the iteration $k+1$ as
$$
X_{\N}^{k+1}(\bmu) = X_{\N}^{k}(\bmu) - \delta X_{\N}
$$
and repeat these steps until an appropriate stopping criterion is verified.

The last ingredient to discover the bifurcation diagram, and so how the stability features of a solution change with respect to the variation of the parameter $\bmu$, consists in the capability of following the solution branches.
Thus, one is interested in the \textit{continuation methods}, allowing to generate a sequence of solutions, corresponding to the selected values of the parameter, in order to construct branches of possible configurations. 

Following a branch is not an easy task, in fact, during the investigation of the parameter space, frequent issues dealing with the continuation step are: 
\begin{itemize}
\item[$\circ$] $\Delta \bmu$ is too large, such that we can possibly skip the bifurcation point, without noticing the branching behavior;
\item[$\circ$] $\Delta \bmu$ is too small, causing a waste of computational resources, mostly when approximating regions far from the critical point.
\end{itemize}

In the following, we consider the \textit{simple continuation} method, in which, as initial guess for a certain parameter $\bmu$, one assign the actual solution corresponding to the parameter $\bmu - \Delta \bmu$.

The combination of this three methodology allow us to reconstruct a branch of the bifurcation diagram, and for this reason it is called \textit{branch-wise} procedure. In this case, with a little abuse of terminology, we call branch the unique extension of the bifurcating behavior to the pre-bifurcation regime. 

Projecting the PDE in the high-fidelity space with $\mathcal N$ degrees of freedom is a computational bottleneck. Indeed, complex phenomena require a fine discretization of the domain, resulting in a possibly huge linear system that has to be solved repeatedly for each iteration of the Newton method and for each parameter in $\mathcal{P}$.

\subsection{Reduced order models}\label{se:rom}
In order to reduce the computational cost of the previously described methodology, \textit{Reduced Order Models} (ROMs) \cite{patera07:book,hesthavenCertifiedReducedBasis2015,quarteroniReducedBasisMethods2016,morepas2017,benner2017model} have been introduced a few years ago, aiming at replacing the standard high-fidelity approximations with a class of reduced techniques.

Among the others, we consider the \textit{Reduced Basis} (RB) method, which consists in a projection-based technique closely related to the Galerkin FE method.

The key feature of the RB method is the offline-online paradigm, allowing for an efficient reconstruction of the solution to \eqref{eq:weakgal}, for a given parameter $\bmu \in \Pa$, exploiting the following two steps:
\begin{itemize}
\item[$\circ$] An  \textit{offline phase}: approximated solutions to \eqref{eq:weakgal} corresponding 
to representative parameters values/system configurations are computed 
with a Finite Element method, and stored together with the assembled operators defining the problem. 
This is a computationally expensive step, and can be done through high performance computing facilities.
\item[$\circ$] An \textit{online phase}: the pre-processed information obtained during the offline phase is exploited to obtain real-time computations of the solution for each new instance of the parameter, without requiring powerful devices.
\end{itemize}

Thus, the main goal of the offline phase is the construction of a low dimensional basis for a discrete manifold $\mathbb X_N \subset \mathbb X_\mathcal{N}$, called \textit{reduced manifold}, which we assume to well approximate the high fidelity manifold $X_\mathcal{N}$.
This entails solving $N_{train}$ times the Galerkin high fidelity problem associated to
$N_{train}$ values of $\bmu$ in $\Pa$. The obtained solutions $\{X(\bmu_i)\}_{i=1}^{N_{train}}$ are usually called \textit{snapshots}, and, in this work, we exploit the \textit{Proper Orthogonal Decomposition} (POD) approach \cite{hesthavenCertifiedReducedBasis2015, patera07:book, quarteroniReducedBasisMethods2016} applied to this dataset to build the reduced basis which spans the low-dimensional manifold.
The POD is a compression strategy based on the singular value decomposition (SVD) which extracts the main features of the problems while retaining most of its energy, since it is shown to be optimal in $l^2(\R^N)$ sense.
For the sake off efficiency, one usually consider the SVD of the \textit{snapshots matrix}
\begin{equation}\label{eq:snapshots}
\mathsf S = [X_\N(\bmu_1), \dots, X_\N(\bmu_{N_{train}})] \in \R^{\N \times N_{train}},
\end{equation}
whose columns are the degrees of freedom of the $N_{train}$ snapshots, and  select as the $N$-dimensional POD basis the one spanned by the first $N$ left singular vectors.

Once built the reduced basis space, we project the weak formulation \eqref{eq:weakgal} in $\X_N$, obtaining the following problem: given $\bmu \in \Pa$, we seek $X_N(\bmu) \in \X_N$ such that
\begin{equation}
\label{eq:weakred}
g(X_N(\bmu), Y_N ; \bmu) = 0, \quad \forall \, Y_N \in \X_N .
\end{equation}
As before, we need to apply the Newton-Kantorovich method that reads as: chosen an initial guess $X^0_N(\bmu) \in \X_N$, for every $k = 0, 1, \dots$, find the variation $\delta X_N \in \X_N$ such that:
\begin{equation}
\label{eq:weakrednewton}
\mathrm{d}g[X_N^k(\bmu)](\delta X_N, Y_N; \bmu) =  g(X_N^k(\bmu), Y_N ; \bmu), \quad \forall \, Y_N \in \X_N , 
\end{equation}
and then we update the solution as $X_N^{k+1}(\bmu) = X_N^{k}(\bmu) - \delta X_N$, until an appropriate stopping criterion is verified.

At this point, the complexity reduction enabling efficiency is based on two main assumptions: 
\begin{itemize}
\item[(i)] the \textit{affine decomposition}, i.e.\ one can rewrite the weak formulation as a linear combination of $\bmu$-independent forms and $\bmu$-dependent coefficients, to assemble the system efficiently.
\item[(ii)] The discrete manifold $\mathbb X_{\mathcal{N}}$ can be approximated with a low-dimensional space $\mathbb X_N$, i.e.\ we have $N \ll \N$.
\end{itemize}

Therefore, the reduced computational cost mainly comes from avoiding to project on the large FE manifold, while relying on the small RB one.
Unfortunately, the assumption (i) is not fulfilled here, because of the $\bmu$-dependence of the solution around which we linearize the nonlinear weak formulation.
To overcome this issue, and recover a consistent speed-up, we adopt a hyper-reduction method, allowing for interpolating the nonlinear terms and recovering the affinity.

\subsection{Hyper-reduction strategies}\label{sec:hyper}
As we discussed in the previous section, the obtainment of a $\N$-independent online phase is prevented by the presence of nonlinear terms.
Indeed, even assuming that we can write the discrete counterpart of the PDE operator as
\begin{equation}
\label{eq:affine_residual_hf}
\mathsf G_\N (\mathsf X_\N; \bmu) = \sum_{q=1}^{Q_G} \theta_G^q(\bmu)\mathsf G_{\N}^q(\mathsf{X}_{\N}) ,
\end{equation}
its projection in the reduced space $\X_N$ is given by
\begin{equation}
\label{eq:affine_residual_ro}
\mathsf G_N (\mathsf X_{N}; \bmu) = \sum_{q=1}^{Q_G} \theta_G^q(\bmu) \mathsf V^{T}\mathsf G_{\N}^q(\mathsf V \mathsf{X}_{N}) ,
\end{equation}
where $\mathsf V$ encodes the reduced basis, the assembly and the projection steps still involve the degrees of freedom of the original problem. A tensor assembly of the residual can be exploited, but this is practicable only when dealing with polynomial nonlinearities, which is not the case for hyperelastic problems.
Thus, it is fundamental to exploit affine-recovery techniques to efficiently perform online evaluations. In this work, we consider a hyper-reduction approach based on the Empirical Interpolation Method (EIM) \cite{barrault04:_empir_inter_method,maday2009general} and its variant the Discrete Empirical Interpolation Method (DEIM) \cite{Chaturantabut2010}.

The hyper-reduction techniques provide an affine approximation of the reduced residual vector of the form
\begin{equation*}
\mathsf G_N(\mathsf{X}_N; \bmu) \approx \sum_{q=1}^{Q_G} c_{\theta}^q(\mathsf{X}_N; \bmu)\mathsf{V}^T h^q
\end{equation*}
where $\{h^q\}_{q=1}^{Q_G}$ represent a suitable basis and $c_{\theta}^q$ are the interpolation coefficients.

Finally, to obtain an efficient recovery of the bifurcating phenomena, we combine, as done in the high fidelity setting, the projection step in $\X_N$, the Newton method \eqref{eq:weakrednewton} and the simple continuation method for the reduced vector solution $\mathsf{X}_N$.
Given the reduced cost of each evaluation $\bmu_j \to \mathsf{X}_N(\bmu_j)$, the testing set to approximate the branches can be chosen as a much refined version of the one in the high fidelity context.
Therefore, the reduced approach allows to obtain a better investigation of the region near the critical points, capturing the dynamics of the bifurcation with smaller steps.
For what concerns the theoretical rationale behind the reduced basis method \cite{quarteroniReducedBasisMethods2016,hesthavenCertifiedReducedBasis2015}, it was proved that one of the main ingredients to have good approximation properties consists in the parametric regularity of the solutions, rather than its regularity in space.
This is an issue for bifurcation problems, where the critical points represent discontinuities in the \textit{parametric sensitivity} $\p X(\bmu)/\p\bmu$.
Hence, we expect the error analysis in the parameter space $\Pa$ to show higher peaks at bifurcation points.

\section{Continuum mechanics framework for Hyperelastic models}
Let us given $\Omega \subset \mathbb{R}^d$, with $d = 2,3$, a bounded domain as the reference configuration of a body $\mathcal{B}$ that undergoes deformation. We can capture the \textit{displacement} $\mathbf{u}$ of a material point $P$, from the reference position $X$ to the new deformed location $x$, through the deformation function $\phi: \mathbb{R}^d \to \mathbb{R}^d$ defined as $x = \phi(X)$. In fact, to capture the motion of a point, we can write $\mathbf{u}(X) = x - X$, and, from this, it is immediate to define the deformation gradient as $$F = \frac{\partial \phi(X)}{\partial X} = \nabla \mathbf{u} + I .$$ Moreover, we can also define the determinant of the deformation gradient, which encodes the volume changes, as $J = \det(F)$.

Once these quantities have been defined, we can present the standard and well-known equilibrium equation derived from the equation of motion \cite{Timoshenko,ciarlet1988three,ciarlet1997mathematical,Huynh2018,veroy03:_phd_thesis}
\begin{equation}
\label{eq:eq_1}
- \Div(P) = B  \quad \text{in} \ \Omega ,
\end{equation} 
where $B$ is an external force acting on the whole domain $\Omega$ and $P = P(F)$ is the first Piola-Kirchhoff stress tensor, related to the Cauchy stress tensor $\sigma$ by the formula $P = J \sigma F^{-T}$.

Of course, the computation of the displacement field $\mathbf{u}$ of the considered body as a function of the external loads requires complementing equation \eqref{eq:eq_1} by a suitable law. 
In particular, if we want to study the equation \eqref{eq:eq_1} in the context of elasticity problems, we have to characterize it through a relation between the stresses (forces) and the strains (displacements). 

These relationships, known as \textit{constitutive equations}, highly depend on the type of material under consideration.
In this context, we relied on the so called \textit{hyperelastic} material (or \textit{Green elastic} material) \cite{gurtin1982introduction}, which is characterized by the existence of a function which defines the Piola-Kirchhoff stress tensor.

More precisely, we can define a \textit{strain energy} function $\psi = \psi(F)$ such that 
\begin{equation}\label{eq:strain_energy}
P(F) = \frac{\partial\psi(F)}{\partial F} ,
\end{equation}
which is used to define a hyperelastic material, assuming that the stress can be obtained by taking the derivative of the energy $\psi$ with respect to the strain $F$. 

Of course, many simplifications can be adopted when dealing with specific type of materials. As an example, for a hyperelastic material which is also \textit{isotropic} (its properties are independent of the direction of examination), the strain energy function can be expressed only by means of the \textit{principal stretches}. Indeed, from these one can obtain the \textit{principal invariants} of the (left) Cauchy-Green deformation tensor $C = F^T F$, which are defined as
\begin{equation}\label{eq:invariants}
\begin{split}
I_1 &= \tr(C),\\
I_2 &= \frac{1}{2}[\tr(C)^2 - \tr(C^2)], \\
I_3 &= \det(C). 
\end{split}
\end{equation}

Moreover, we remark that these materials are characterized by the fact that the work done by stresses does not depend on the path of deformation, thus they conserve the total energy, and in contrast with linear elasticity, we do not have to require any infinitesimal assumption on the strains.
We can now proceed with the description of the models considered for the study of compressed beams.
In particular, we will analyze two different choices for the strain energy function $\psi$, investigating numerically their properties with respect to the buckling problems. 

The simplest constitutive relation is the so-called \textit{Saint Venant-Kirchhoff} (SVK) model, which is an extension of the linear elastic material model to the geometrically nonlinear regime. This model is defined through the strain energy function defined as 
\begin{equation}\label{eq:svk}
\psi(F) = \lambda_1 \mathcal{E} : \mathcal{E} + \lambda_2(\tr(\mathcal{E}))^2/2 , 
\end{equation}
 where $$\mathcal{E} = \frac{1}{2}(F^TF - I)$$ is the Green-Lagrange strain tensor (a measure for varying lengths between points), $\lambda_1$ and $\lambda_2$ are the Lam\'e constants, which are related to the material properties through the Young modulus $E$ and the Poisson ratio $\nu$ as follows $$\lambda_1 = \frac{E}{2(1 + \nu)}, \quad \lambda_2 = \frac{E\nu}{(1+\nu)(1-2\nu)}.$$
Moreover, from the definition of the hyperelastic material in \eqref{eq:strain_energy}, we obtain that the first Piola-Kirchhoff stress tensor for the SVK model is given by $$P(F) = F(2\lambda_1 \mathcal{E} + \lambda_2 \tr(\mathcal{E})I).$$
We remark that, as we said before, such model is analogous to the linear elasticity one. Indeed, one can characterize the linear elasticity model by means of the strain energy function 
\begin{equation}\label{eq:el_lin}
\psi(F) = \lambda_1 \epsilon : \epsilon + \lambda_2 (\tr(\epsilon))^2 , 
\end{equation}
where $$\epsilon = \frac{1}{2}(F + F^T - I)$$ is the infinitesimal strain tensor.

The second hyperelastic constitutive relation we will consider is the so called \textit{neo-Hookean} (NH) model, which can be expressed through the strain energy function defined as 
\begin{equation}\label{eq:nh}
\psi(F) = \frac{\lambda_1}{2}(I_1 - 3) - \lambda_1 \ln J + \frac{\lambda_2}{2}(\ln J)^2 , 
\end{equation}
where $I_1$ is the first principal invariants in \eqref{eq:invariants} and $\lambda_1$, $\lambda_2$ are the Lam\'e constants as before. The first Piola-Kirchhoff stress tensor in this case is given by $$P(F) = F[(\lambda_1(I - C^{-1}) + \lambda_2(\ln J) C^{-1}] .$$

We remark that the words buckling and bifurcation will be used here indistinctly, but we will be interested only in the approximation of the first buckling mode, properly following the post-buckling behavior as the target branch.
In the next section, we will present the mathematical formulation of the problem, focusing on its weak formulation that will allow us to apply our methodology. 

\subsection{Weak formulation}\label{se:hyp_weak_form}
Having presented the constitutive relations considered, we can now state the mathematical formulation of the problem.
Starting from equation \eqref{eq:eq_1} we can consider the boundary value problem, in the reference domain $\Om$, given by

\begin{equation}
\label{eq:hyperelastic_pde}
\begin{cases}
- \Div(P(\mathbf{u})) = B  \quad &\text{in} \ \Omega , \\
\mathbf{u} = \mathbf{u}_D \quad &\text{in} \ \Gamma_D , \\
P(\mathbf{u})n = T  \quad &\text{in} \ \Gamma_N ,
\end{cases}
\end{equation}
where $\mathbf{u}: \Om \to \R^d$ is the in-plane displacement, $B$ is the \textit{body force} per unit reference area, and $T$ is the \textit{traction force} per unit reference length.
Thus, our aim is to study the deformation of the domain $\Omega$, subjected to a prescribed displacement $\mathbf{u}_D$  on the Dirichlet boundary $\Gamma_D$ together with body and traction forces.
We remark that, due to the Dirichlet non-homogeneous BCs, the function space in which we set the problem is the space $H_{D}^1(\Om) = \left\{\mathbf{v} \in H^1(\Om) \ | \  \mathbf{v} = \mathbf{u}_D \in \Gamma_D\right\}$, where in this case we are assuming that each component of $\mathbf{u}$ is in $H_{D}^1(\Om)$.

To derive the weak formulation of the problem, we proceed as usual by taking the dot product with a test function $\mathbf{v} \in \X = H_0^1(\Om)$ and integrating over the reference domain $\Omega$, this way we obtain 
\begin{equation}
\label{eq:weak_hyp}
\int_{\Om} \Div(P(\mathbf{u})) \cdot \mathbf{v} \ d\Om + \int_{\Om}B \cdot \mathbf{v} \ d\Om = 0 \quad \forall \, \mathbf{v} \in \X .
\end{equation}  
Applying the divergence theorem and embedding the traction and displacement boundary conditions (the test function $\mathbf{v} \in \X$ satisfies homogeneous BCs on $\Gamma_D$), the weak formulation reads as: find $\mathbf{u} \in H_{D}^1(\Om)$ such that 

\begin{equation}
\label{eq:weak_hyp_2}
\int_{\Om} P(\mathbf{u}) : \nabla \mathbf{v} \ d\Om  - \int_{\Om}B \cdot \mathbf{v} \ d\Om - \int_{\Gamma_N}T \cdot \mathbf{v} \ d\Gamma = 0 \quad \forall \, \mathbf{v} \in \X .
\end{equation}  

In practice, a standard technique to set the problem in the same space is to consider a lifting function $R_D \in H^1(\Om)$ such that $R_D|_{\Gamma_D} = \mathbf{u}_D$. Then, one sets $\mathring{\mathbf{u}} = \mathbf{u} - R_D$ and, from the immediate consideration that $\mathring{\mathbf{u}} \in \X$, the symmetry in spaces between solution and test functions is restored.

Furthermore, the boundary value problem in \eqref{eq:hyperelastic_pde} for hyperelastic media can be also expressed as a minimization problem by means of the \textit{Theorem of Virtual Work} \cite{gurtin1982introduction}. 
In fact, we can define the \textit{potential energy} of the beam in terms of the strain energy function $\psi$ as 
\begin{equation}
\label{eq_2}
\Pi(\mathbf{u}) = \int_\Omega {\psi(\mathbf{u})\ d\Om} -\int_\Omega {B \cdot \mathbf{u}\ d\Om} - \int_{\Gamma_N} {T \cdot \mathbf{u}\ d\Gamma} .
\end{equation}
At minimum point of $\Pi(\mathbf{u})$, the directional derivative of $\Pi$ w.r.t.\ the change in $\mathbf{u}$ is given by 
\begin{equation}\label{eq:weak_hyp_comp}
g(\mathbf{u}, \mathbf{v}) \eqdot D_v \Pi(\mathbf{u}) = \evalat[\bigg]{\frac{\mathrm{d}\Pi(\mathbf{u} + \delta \mathbf{v})}{\mathrm{d}\delta}}{\delta=0} ,
\end{equation}
 is equal to zero for all $\mathbf{v} \in \X$, that is 
\begin{equation}
\label{eq_3}
g(\mathbf{u}, \mathbf{v}) = 0 \quad \forall \ \mathbf{v} \in \X .
\end{equation}

It is easy to observe that evaluating equation \eqref{eq:weak_hyp_comp} leads to the weak formulation in \eqref{eq:weak_hyp_2}.
The constitutive relations we chose are characterized by the fact that equation \eqref{eq_3} results in a nonlinear form w.r.t.\ $\mathbf{u}$.  In such cases we need to rely on the Jacobian of $g$, that is defined by 
\begin{equation}
\label{eq_4}
\mathrm{d}g[\mathbf{u}](z, \mathbf{v}) = \evalat[\bigg]{\frac{\mathrm{d}g(\mathbf{u} + \delta z, \mathbf{v})}{\mathrm{d}\delta}}{\delta=0} .
\end{equation}

Before moving to the numerical approximation of the problem, we want to highlight that the variational forms above were presented without any dependency from a generic (multi-)parameter $\bmu$.
The reason behind this choice is that we will consider different geometrical/physical parameter for each toy problem, specifying case by case the parametrized quantities.

As regards the numerical approximation of hyperelastic problems, the projection of the weak formulation \eqref{eq:weak_hyp_2} in the Finite Element space $\X_\N$ of dimension $\N$, reads as: given $\bmu \in \Pa$, seek $u_\mathcal{N} \eqdot u_\mathcal{N}(\bmu) \in \X_\mathcal{N}$ that satisfies 
\begin{equation}
\label{eq:weakgal_hyp}
g(u_\mathcal{N}, v_\mathcal{N} ; \bmu) = 0\   \quad \forall \, v_\mathcal{N} \in \X_\mathcal{N}, 
\end{equation}
where we introduced in the weak formulation the geometrical and/or physical multi-parameter $\bmu \in \Pa$.
The application of the Newton method leads us to solve, at the generic $k$-th step, the algebraic equation 
\begin{equation}
\label{eq:linearnewtgal_hyp}
\mathsf{J}_\N(\mathsf{u}_{\N}^k(\bmu); \bmu) \delta\mathsf{u}_{\N} = \mathsf{G}_{\N}(\mathsf{u}_{\N}^k(\bmu); \bmu) ,
\end{equation}
updating the solution as $\mathsf{u}_{\N}^{k+1} = \mathsf{u}_{\N}^k - \delta\mathsf{u}_{\N}$ until convergence.

In order to efficiently recover the solution for repeated instances of the parameter $\bmu$ during the online stage, we rely on the reduced basis projection onto $\X_N$. Thus, given $\bmu \in \Pa$, we seek $u_N \eqdot u_N(\bmu) \in \X_N$ that satisfies 
\begin{equation}
\label{eq:weakgal_hyp_red}
g(u_N, v_N ; \bmu) = 0\   \quad \forall \, v_N \in \X_N ,
\end{equation}
which form the algebraic standpoint translates into  
\begin{equation}
\label{eq:linearnewtgal_hyp_red}
\mathsf{J}_N(\mathsf{u}_{N}^k(\bmu); \bmu) \delta\mathsf{u}_{N} = \mathsf{G}_{N}(\mathsf{u}_{N}^k(\bmu); \bmu) ,
\end{equation}
in terms of the reduced residual and the reduced Jacobian matrix, respectively.

\section{Numerical results}
Here, we will present different toy problems, in which we will study the buckling behaviors of beams with different constitutive relations, boundary conditions, geometry settings and external forces. We will discuss several implications of the exploitation of ROMs approaches to obtain an efficient computation of the bifurcation diagrams.

\subsection{2-D toy problem}\label{sec:2d_toy}
In this first application, we consider a two-dimensional beam which corresponds to the domain $\Omega = [0,1] \times [0, 0.1]$. We built the Finite Element space with $\mathbb{P}_1$ linear elements, resulting in a high fidelity dimension $\N=4328$. It is well-known that, e.g. for $B = (0, 0)$ and $T = 0$, the bifurcation diagram enjoys a $\mathbb{Z}^2$-reflective symmetry, and it undergoes a series of pitchfork bifurcations as $\mu$ is increased. Despite this, here we focused only on the first branch, being this usually the most common cause of failure.
In this simple 2-D context, we will analyze different test cases in which the buckling properties will be studied in connection to different compression conditions, materials characteristics and geometries.

\subsubsection{Dirichlet compression}\label{sec:hyp_dir}
Here, we consider a beam subjected to a parametrized uniform compression imposed by means of Dirichlet boundary conditions. In particular, we will analyze the Saint Venant-Kirchhoff and the neo-Hookean models, choosing:   a null traction force $T = 0$, a Young modulus $E = 10^6$ and a Poisson ratio $\nu = 0.3$, while for the body force we will study either $B = (0, 0)$ or $B = (0, -1000)$.
In particular, we want to study the buckling of the beam subjected to a compression of magnitude $\mu$ on its right. To do so, we rewrite the Dirichlet conditions in \eqref{eq:hyperelastic_pde} as
\begin{equation}
\label{eq_5}
\begin{split}
u = (0, 0) \quad \ &\text{on} \ \Gamma^{l}_D = \lbrace{0\rbrace} \times [0, 0.1] , \\
u = (- \mu, 0) \quad \ &\text{on} \ \Gamma^{r}_D = \lbrace{1\rbrace} \times [0, 0.1] ,
\end{split}
\end{equation}
which correspond to a clamped condition on the left end of the beam, and an increasing uniform uni-axial compression on the other end.
Thus, for a fixed compression $\mu$, we consider the functional space $$\X = \lbrace{ u \in (H^1(\Omega))^2: u = (0, 0) \ \text{on}\  \Gamma^{l}_D,\ u = (-\mu, 0) \ \text{on} \  \Gamma^{r}_D \rbrace} .$$

Despite the simplicity of the models, the buckling phenomenon makes the analysis of the behavior of the solution with respect to the compression $\mu$ not straightforward, especially concerning the reduction strategies. Indeed, for the high fidelity setting we had to choose $N_{train} = 1000$ points in the parameter space $\Pa = [0, 0.2]$, in order to properly detect the critical point and follow the post-buckling branch. We applied for all the test cases a POD tolerance $\epsilon_{POD} = 10^{-8}$ and an online continuation method based on $K = 2000$ equispaced points in $\Pa$, which corresponds to a continuation step $\Delta \mu = 10^{-4}$.

\begin{figure}[h!]
\centering
     \begin{subfigure}[b]{0.49\textwidth}
         \centering
         \includegraphics[width=\textwidth]{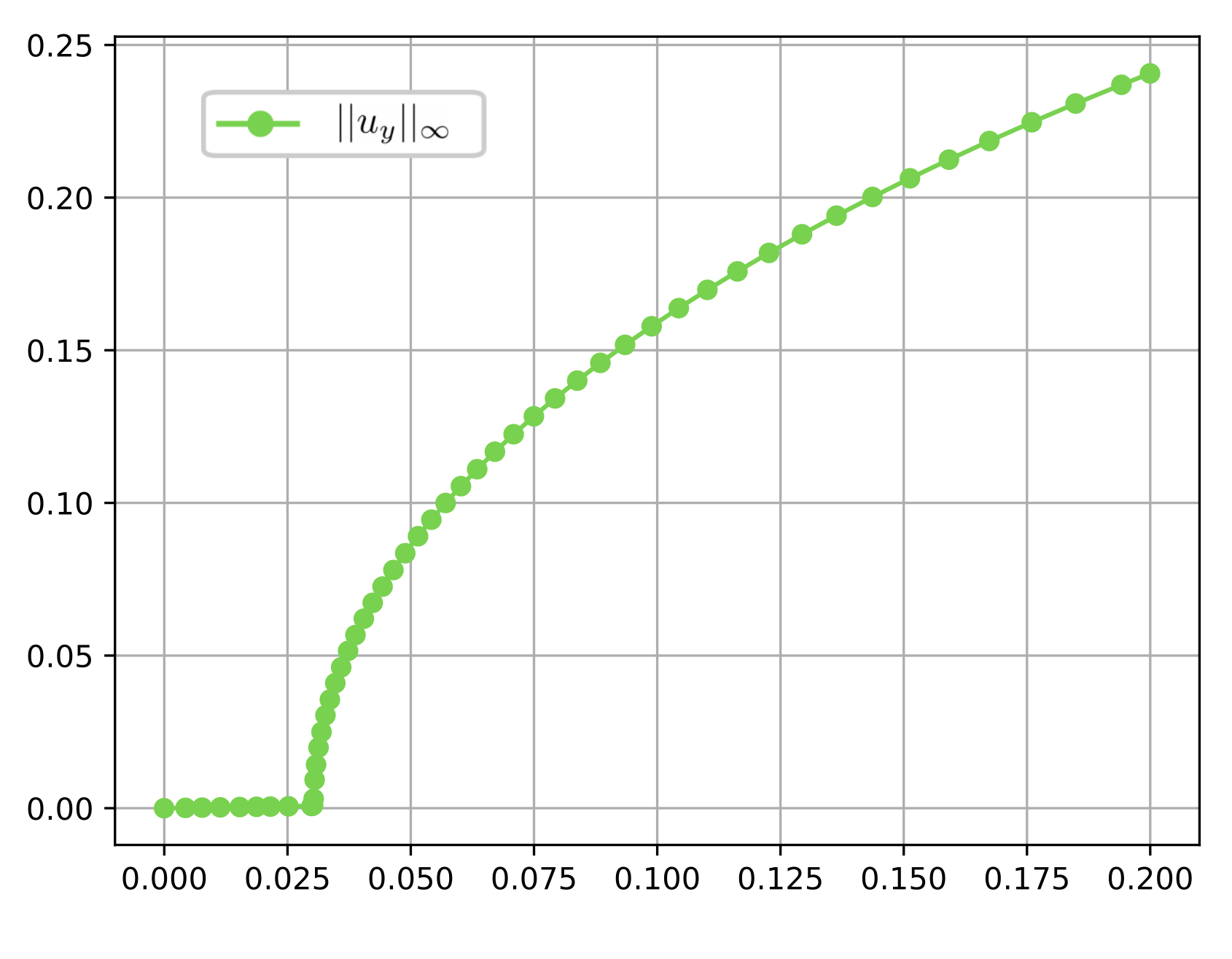}
         \put(-95,4){\makebox(0,0){$\mu$}}
         \caption{$B = (0, 0)$}
         \label{fig:9_a}
     \end{subfigure}
     \hfill
     \begin{subfigure}[b]{0.49\textwidth}
         \centering
         \includegraphics[width=\textwidth]{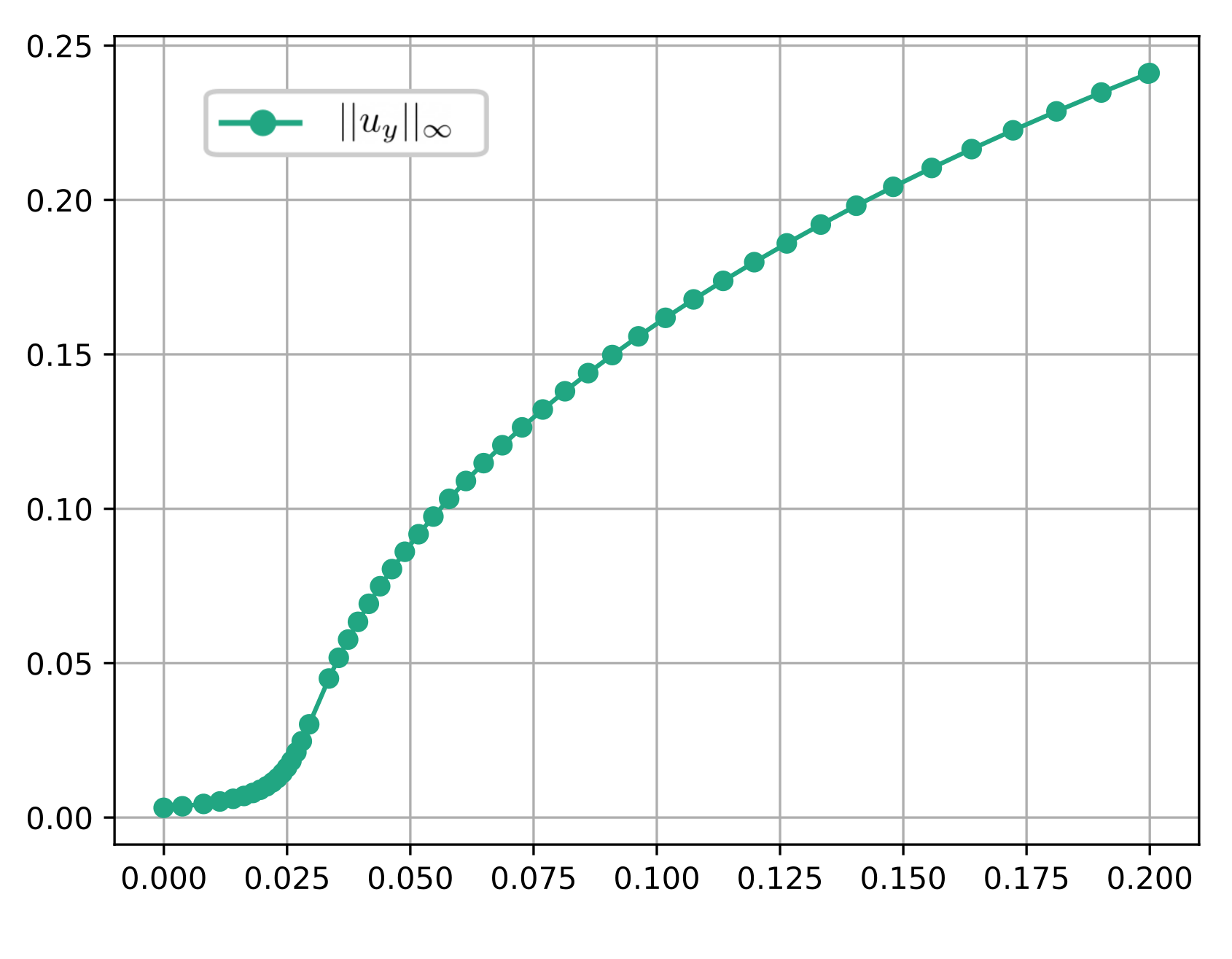}
         \put(-95,4){\makebox(0,0){$\mu$}}
         \caption{$B = (0, -1000)$}
         \label{fig:9_b}
     \end{subfigure}
	\caption{Reduced basis bifurcation diagrams for the SVK beam with different body forces.}
\label{fig_9}
\end{figure}

\begin{figure}[h!]
\centering
\includegraphics[width=0.5\textwidth]{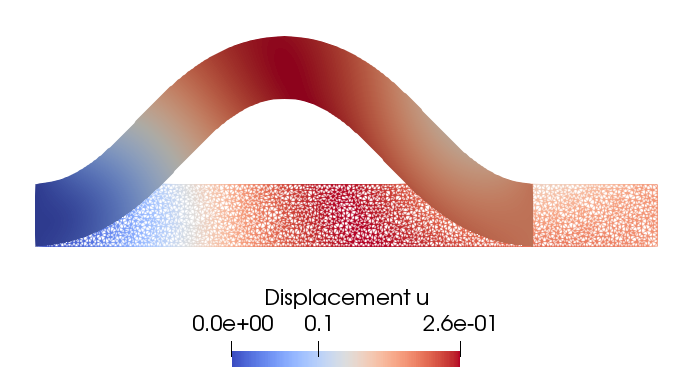}
\caption{High fidelity displacement $u$ for the SVK beam with $B = (0, 0)$ at $\mu = 0.2$.}
\label{fig_6}
\end{figure}

Let us now consider the reduced order strategy for the SVK model with null body force $B = (0, 0)$.
We obtained a reduced basis space of dimension $N = 5$, that we used for the online projection in order to reconstruct the bifurcation diagram in Figure \ref{fig:9_a}. The functional considered in this case is the infinite norm of the second component of the displacement, namely $s(u) = \norm{u_y}_{\infty}$. 
We can clearly observe that the buckling of the beam occurs for the value $\mu^* \approx 0.03$. In fact, in such point the vertical component of the displacement $u$ changes suddenly from being trivial to causing the buckling.  Furthermore, we want to highlight the sharp discontinuity in the sensitivity $\frac{\p u(\mu)}{\p \mu}$ at the buckling point.   
A representative solution of the post-buckling branch is depicted in Figure \ref{fig_6} for $\mu = 0.2$ with respect to the original un-deformed configuration (mesh wireframe). We remark that here the beam buckles upwards.

\begin{figure}[h!]
\centering
     \begin{subfigure}[b]{0.48\textwidth}
         \centering
         \includegraphics[width=\textwidth]{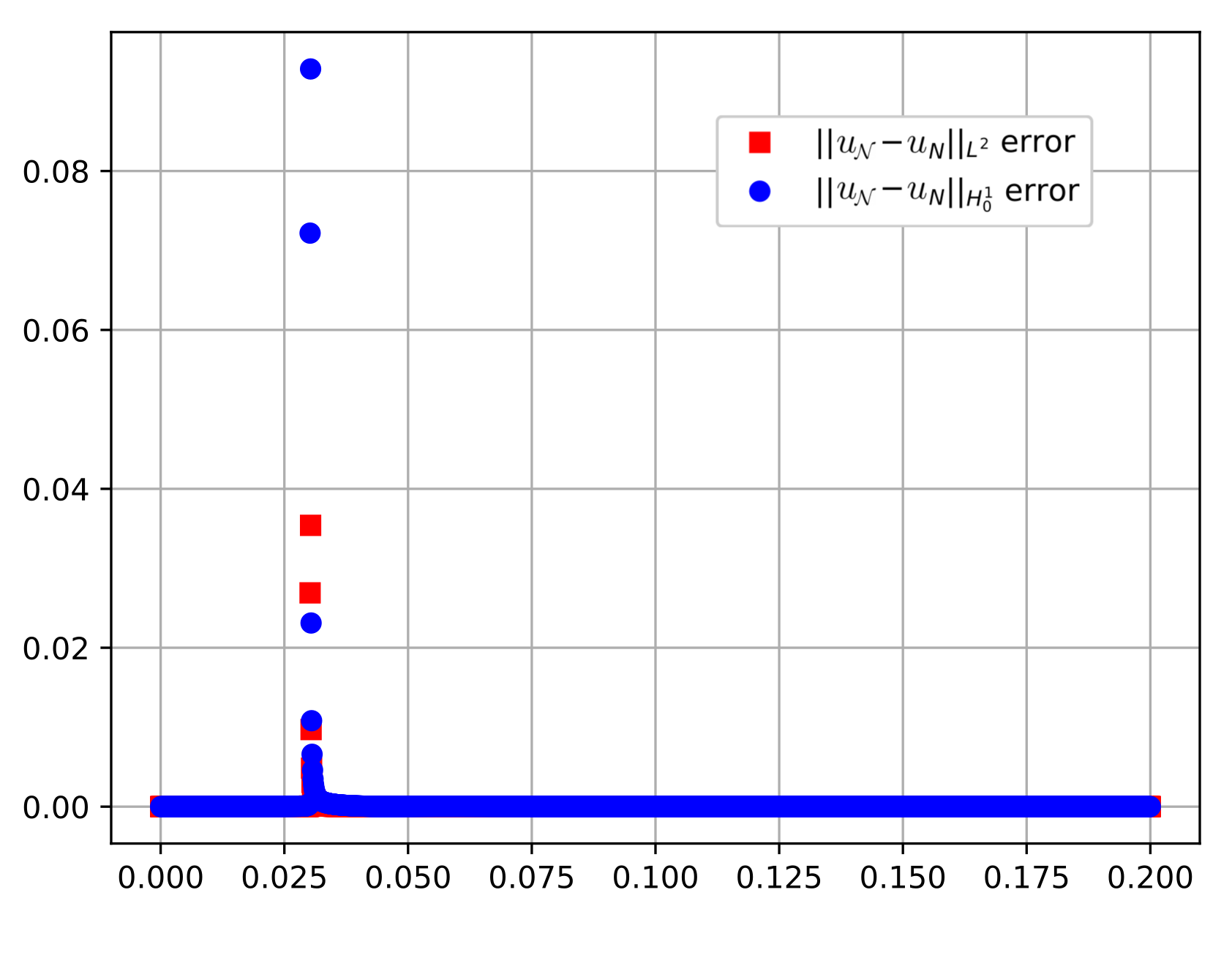}
          \put(-95,4){\makebox(0,0){$\mu$}}
         \caption{$B = (0, 0)$}
         \label{fig:10_a}
     \end{subfigure}
     \hfill
     \begin{subfigure}[b]{0.49\textwidth}
         \centering
         \includegraphics[width=\textwidth]{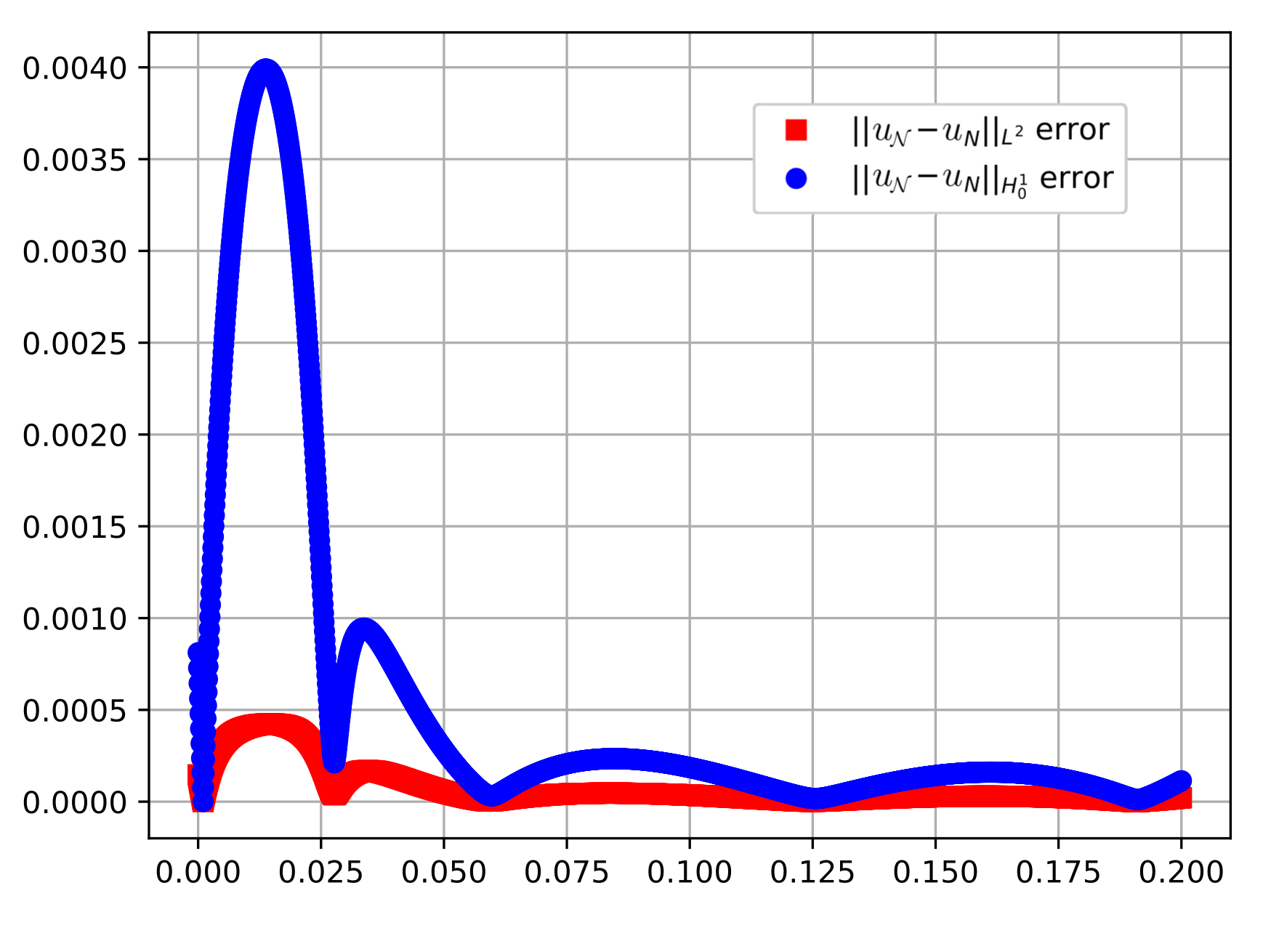}
         \put(-93,4){\makebox(0,0){$\mu$}}
         \caption{$B = (0, -1000)$}
         \label{fig:10_b}
     \end{subfigure}
	\caption{Reduced basis errors for the SVK beam with respect to different body forces.}
\label{fig_10}
\end{figure}

As we can see from Figure \ref{fig:10_a}, the POD with only $N=5$ basis functions is able to recover a good approximation of the bifurcation diagram, with an average error on $\Pa$ of order $10^{-6}$. Despite this, we can observe also here, even more clearly, that the RB error has its maximum at the buckling point. Moreover, here the difference between the maximum and average RB error is of almost 4 orders of magnitude, denoting the great difficulty of the reduced manifold at reproducing the bifurcation.
We also want to remark that, when Dirichlet BCs are considered in the PDE, a lifting function has to be added to the POD basis, in order to encode the value at the boundary.
The speed-up, equals to  1.26, is quite low, in fact to plot the high fidelity version of the bifurcation diagram in Figure \ref{fig_9} we spent $t_{HF} = 454$(s) while the reduced order one required $t_{RB} = 358$(s). This is due to the fact that: (i) the number of degrees of freedom is kept low by the linear elements, (ii) empirical interpolation strategies were not applied, (iii) an $\mathcal{N}$-independent output functional could have been chosen. 

\begin{figure}[h!]
\centering
\includegraphics[width=0.5\textwidth]{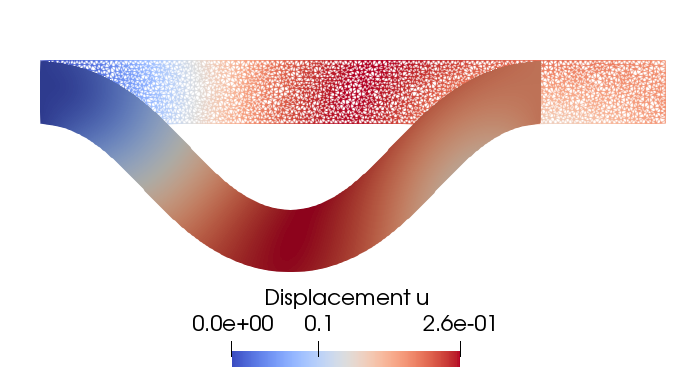}
\caption{High fidelity displacement $u$ for the SVK beam with $B = (0, -1000)$ at $\mu = 0.2$.}
\label{fig_6_B}
\end{figure}

Now we want to study how a perpendicular body force acting on the beam changes the solution properties.
For this reason we consider a non-trivial $B = (0, -1000)$ and within the same reduced setting we obtained $N = 4$ reduced basis dimension upon which we built the buckling diagram in Figure \ref{fig:9_b}. Is it immediate to observe that the sharp discontinuity in the sensitivity was a bit smoothed by the force $B$ which therefore produces a more gradually buckling behavior. Moreover, as it is possible to observe from Figure \ref{fig_6_B}, the action of a ``gravity"-like force causes the branching behavior to be characterized by a downwards buckling.
Hence, the reflective symmetry has been broken by imposing a gravitational body force and the bifurcation diagram behaves more smoothly near the buckling point. Furthermore, due to the elastic properties chosen, the body force is only able to qualitatively change the buckling without consistently affecting the quantitative displacement.

The reduced error plot in Figure \ref{fig:10_b} shows an effect of this smoothness, in fact now the maximum and average error are of the same order. We also notice that here the error is maximum in the region of the unbuckled states. This is an expected behavior given by the greater amount of buckled snapshots computed during the offline phase w.r.t.\ the unbuckled counterpart. Moreover, we remark that the speed-up is the same as the unforced case.
Finally, we highlight that in general a regularization action for the bifurcation diagram could be the key for a good RB strategy, when only a qualitative understanding on the system is required.

\begin{figure}[h!]
\centering
     \begin{subfigure}[b]{0.49\textwidth}
         \centering
         \includegraphics[width=\textwidth]{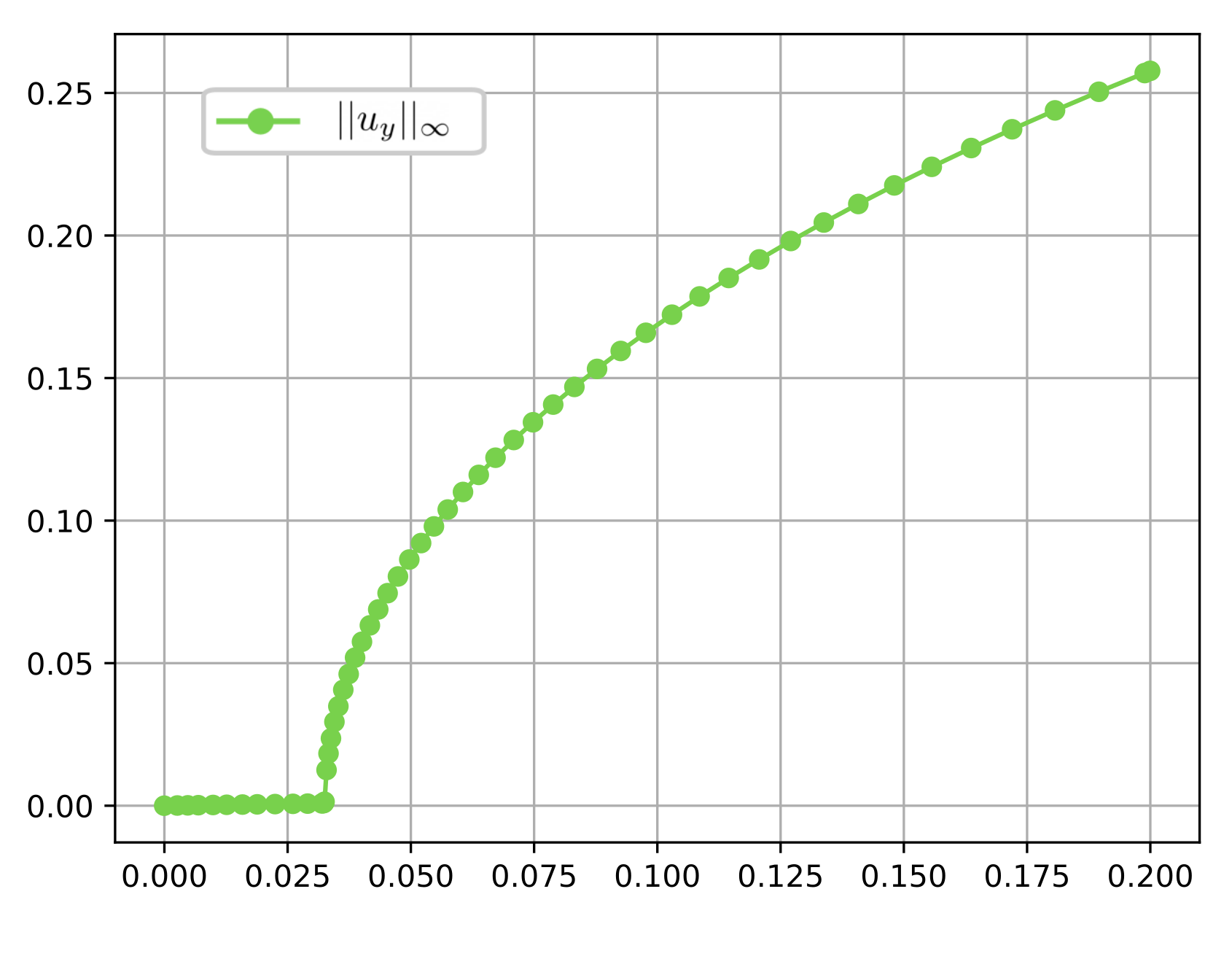}
         \put(-95,4){\makebox(0,0){$\mu$}}
         \caption{$B = (0, 0)$}
         \label{fig:1_a}
     \end{subfigure}
     \hfill
     \begin{subfigure}[b]{0.49\textwidth}
         \centering
         \includegraphics[width=\textwidth]{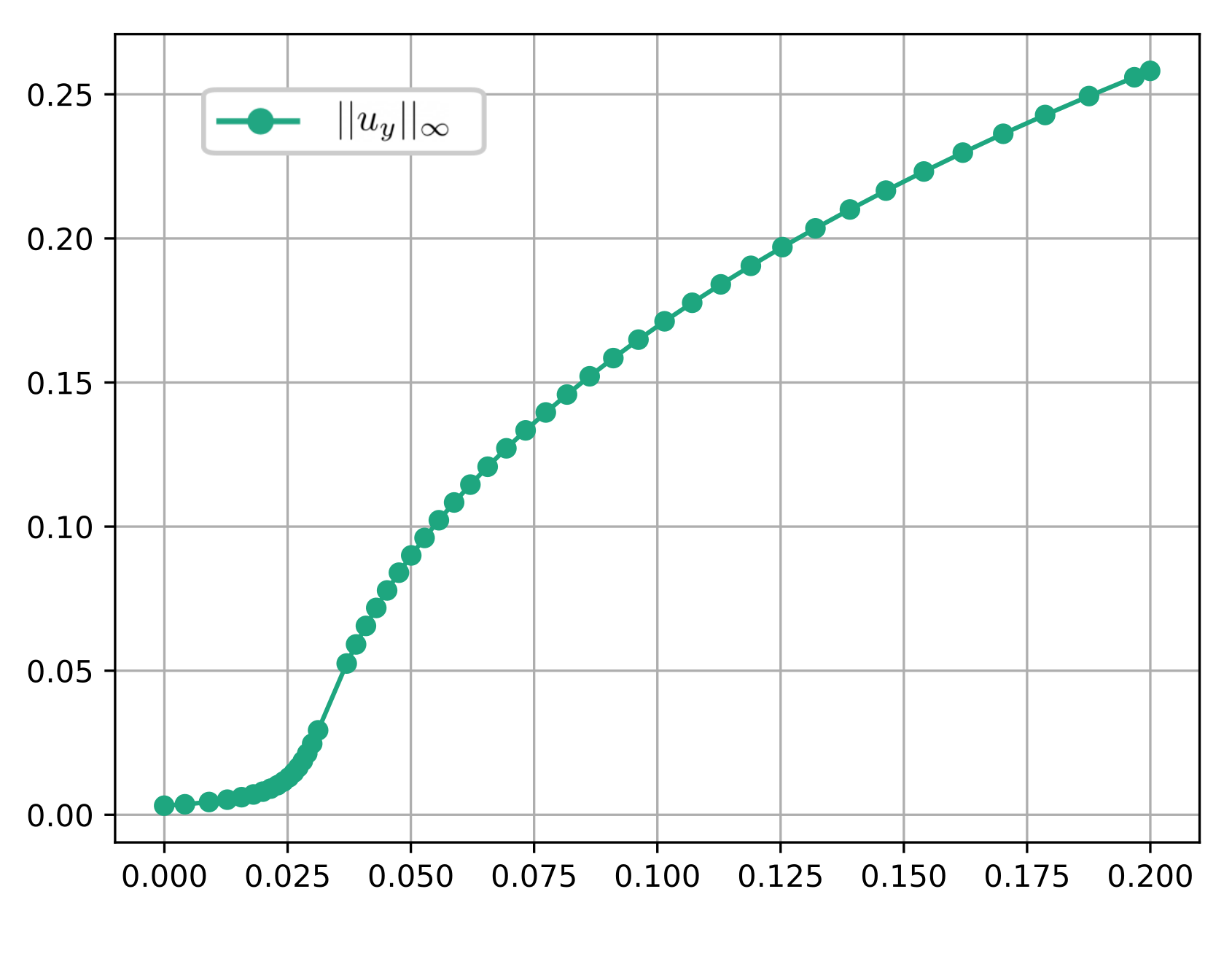}
         \put(-95,4){\makebox(0,0){$\mu$}}
         \caption{$B = (0, -1000)$}
         \label{fig:1_b}
     \end{subfigure}
	\caption{Reduced basis bifurcation diagrams for the NH beam with different body forces.}
\label{fig_1}
\end{figure}

We can now move to the investigation of the neo-Hookean constitutive relation in \eqref{eq:nh}.
Using the same offline setting, we obtained a reduced basis space of dimension $N = 4$ for both body forces test cases.
We can show the reduced basis bifurcation diagrams in Figure \ref{fig_1}, from which it can be observed the same smoothing effect of the body force near the buckling. Even though it is well-known from the literature that the neo-Hookean model is more accurate far from the small displacement regime, we can notice that the detection of buckling point is consistent with the one predicted by the Saint Venant-Kirchhoff model. 
The results do not differ too much form the previous ones, only a small increment of the maximum displacement can be observed in Figure \ref{fig_1}.

\begin{figure}[h]
\centering
\includegraphics[width=0.5\textwidth]{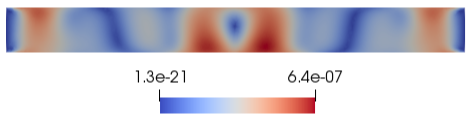}
\caption{Reduced basis error of the displacement $u$ for the NH beam with $B = (0, 0)$ at $\mu = 0.2$.}
\label{fig_err_sol_nh}
\end{figure}

\begin{figure}[h!]
\centering
\includegraphics[width=0.5\textwidth]{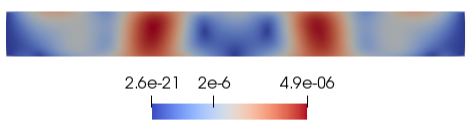}
\caption{Reduced basis error of the displacement $u$ for the NH beam with $B = (0, -1000)$ at $\mu = 0.2$.}
\label{fig_err_sol_B_nh}
\end{figure}

Representative reduced error plots on the displacement of the post-buckling branch for both $B =(0, 0)$ and $B =(0, -1000)$ case are depicted in Figures \ref{fig_err_sol_nh} and \ref{fig_err_sol_B_nh} for $\mu = 0.2$, showing the good approximation accuracy reached by the POD basis for the post-buckling branch.

\begin{figure}[h!]
\centering
     \begin{subfigure}[b]{0.49\textwidth}
         \centering
         \includegraphics[width=\textwidth]{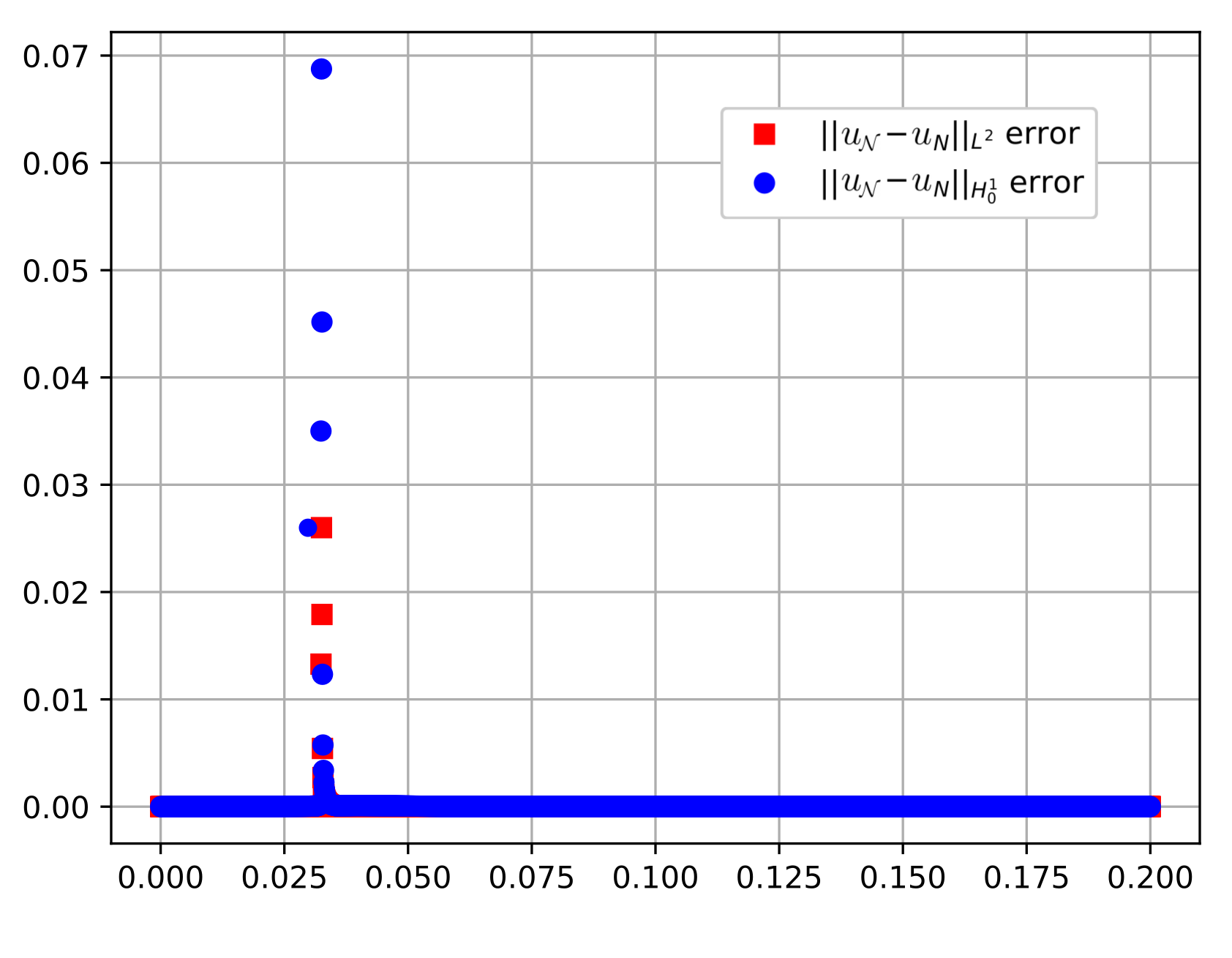}
                  \put(-96,4){\makebox(0,0){$\mu$}}
         \caption{$B = (0, 0)$}
         \label{fig:7_a}
     \end{subfigure}
     \hfill
     \begin{subfigure}[b]{0.49\textwidth}
         \centering
         \includegraphics[width=\textwidth]{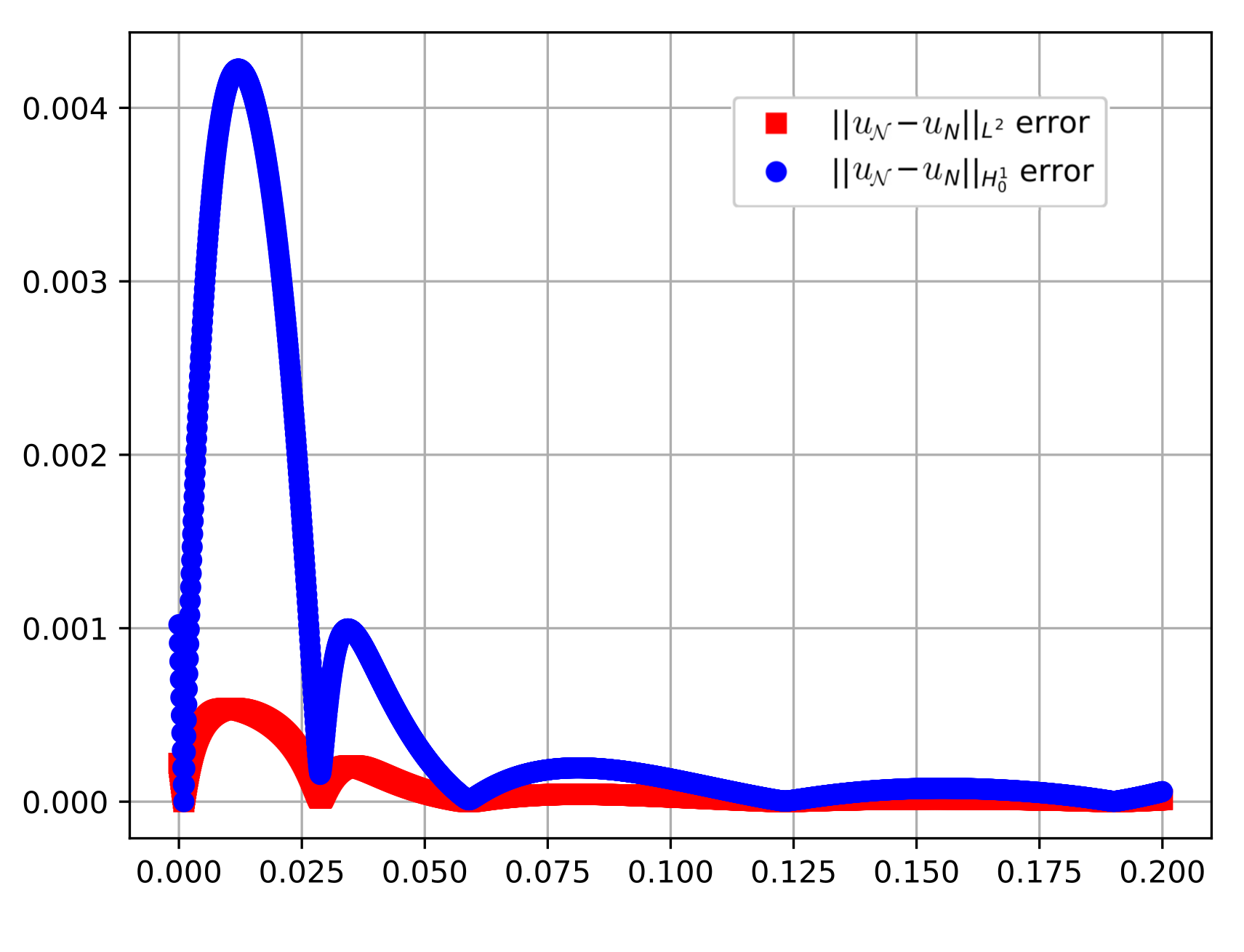}
                  \put(-95,4){\makebox(0,0){$\mu$}}
         \caption{$B = (0, -1000)$}
         \label{fig:7_b}
     \end{subfigure}
	\caption{Reduced basis errors for the NH beam with respect to different body forces.}
\label{fig_7}
\end{figure}

As before, we present the plots of the reduced basis error with respect to each $\mu$ in the parameter domain $\Pa$ in Figure \ref{fig_7}. We can draw the same conclusions about the accuracy of the reduction strategy when both body forces were applied.

\begin{figure}[h!]
\centering
\includegraphics[width=0.5\textwidth]{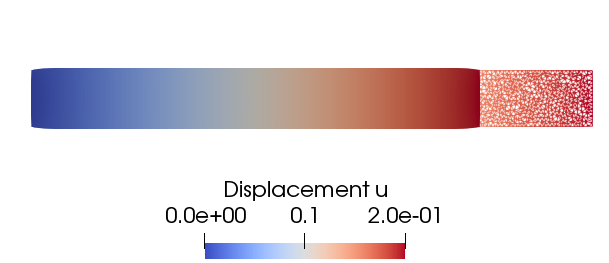}
\caption{High fidelity $0$-th mode displacement $u$ for the SVK beam with $B = (0, 0)$ at $\mu = 0.2$.}
\label{fig_6_0}
\end{figure}

\begin{figure}[h!]
\centering
\includegraphics[width=0.5\textwidth]{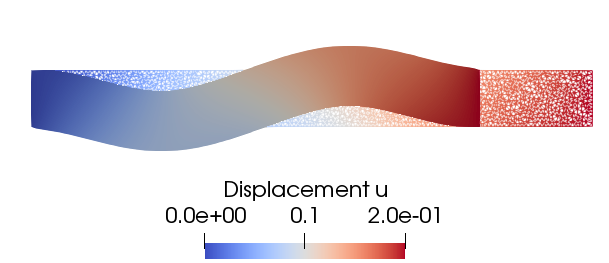}
\caption{High fidelity $2$-th mode displacement $u$ for the SVK beam with $B = (0, 0)$ at $\mu = 0.2$.}
\label{fig_6_2}
\end{figure}

The speed-ups are similar in both forced and unforced cases, of order $1.22$, due to the same reasons. The complexity of the model makes the computation of the bifurcation diagrams more expensive, in fact we spent almost $t_{HF} = 763$(s) for the high fidelity version and $t_{RB} = 624$(s) for the reduced order one. 
Before ending the analysis of this test case, we remark that although we were only interested in the first buckling, as we said previously many configurations can coexist for the same values of the compression parameter $\mu$.
As an example, in Figures \ref{fig_6_0} and \ref{fig_6_2} we plot the high fidelity displacements $u$ at $\mu = 0.2$, corresponding to the $0$-th and $2$-th buckling modes.


\subsubsection{Neumann compression}\label{sec:hyp_neu}

In this section, we will consider a different type of compression, in which we are no longer fixing the displacement at the boundary, rather we parametrize a traction force T. In literature this is usually called the cantilever beam test case.

In particular, we want to investigate the Saint Venant-Kirchhoff model with a Young modulus $E = 10^6$ and a Poisson ratio $\nu = 0.3$, with the trivial body force $B = (0, 0)$. Hence, the compression is modelled through the traction $T = (-\mu, 0)$ on $\Gamma_N^r = \{1\} \times [0, 0.1]$, while the beam is clamped, homogeneous Dirichlet condition $u = (0, 0)$, at the opposite edge $\Gamma_D^l = \{0\} \times [0, 0.1]$.

Of course, given the changed boundary conditions, we expect to obtain a different buckled configuration, since now the beam, thanks to the Neumann condition, has more freedom to move.
For what concerns the high fidelity setting, here we were able to reconstruct the post-buckling branch with much fewer points, indeed we chose $N_{train} = 50$ equispaced points in the parameter space $\Pa = [2200, 2400]$ which represent the magnitude of the traction force $T$. For this model we experimented also the POD approximation accuracy with respect to the number of basis functions involved. Hence, we considered two reduced models, the first with tolerance $\epsilon_{POD} = 10^{-8}$ corresponding to a reduced basis space dimension $N_1 = 4$, and the second one with tolerance $\epsilon_{POD} = 10^{-10}$ corresponding to a reduced basis space dimension $N_2 = 5$. For the online phase we built the branch over $K = 400$ equispaced points in $\Pa$, which corresponds to a continuation step $\Delta \mu =  5 \cdot 10^{-1}$.

\begin{figure}[h!]
\centering
\includegraphics[width=0.6\textwidth]{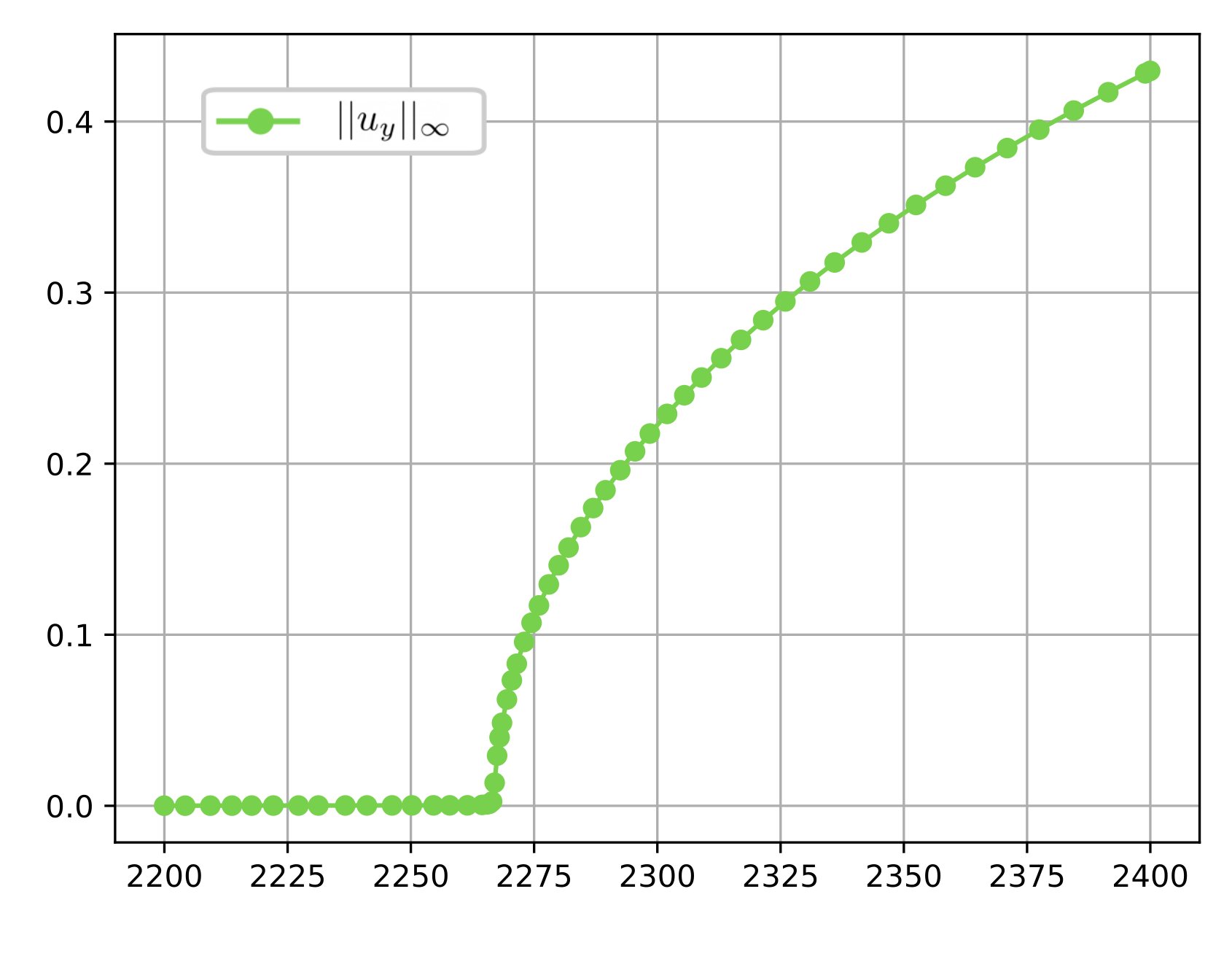}
\put(-115,0){\makebox(0,0){$\mu$}}
\caption{Reduced basis bifurcation diagrams for the SVK beam with Neumann compression and $B = (0, 0)$.}
\label{fig_1_force}
\end{figure}

We considered again the functional $s(u) = \norm{u_y}_\infty$ to plot the reduced basis bifurcation diagram in Figure \ref{fig_1_force}. Since the body force is trivial, consistently with the previous analysis we have a sharp sensitivity near the buckling point $\mu^* \approx 2267$. 

\begin{figure}[h]
\centering
\includegraphics[width=0.4\textwidth]{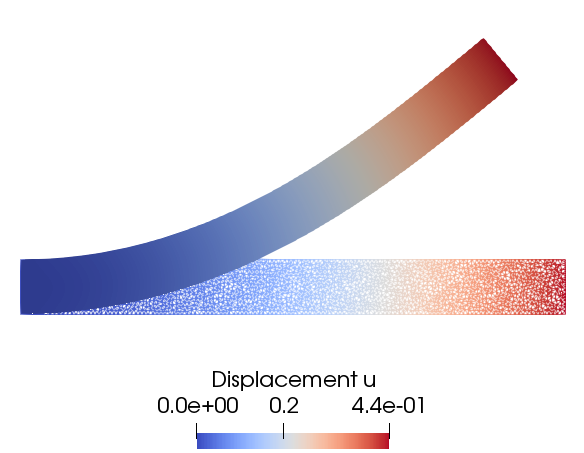}
\caption{High fidelity displacement $u$ for the SVK beam with Neumann compression and $B = (0, 0)$ at $\mu = 0.2$.}
\label{fig_err_sol_svk_force}
\end{figure}

As we expected, the displacement solution shows an upwards buckling with a qualitative different behavior with respect to the Dirichlet compression test case. We plot in Figure \ref{fig_err_sol_svk_force}  a representative solution of the post-buckling branch for the Neumann compression at $\mu = 0.2$.

\begin{figure}[h!]
\centering
     \begin{subfigure}[b]{0.48\textwidth}
         \centering
         \includegraphics[width=\textwidth]{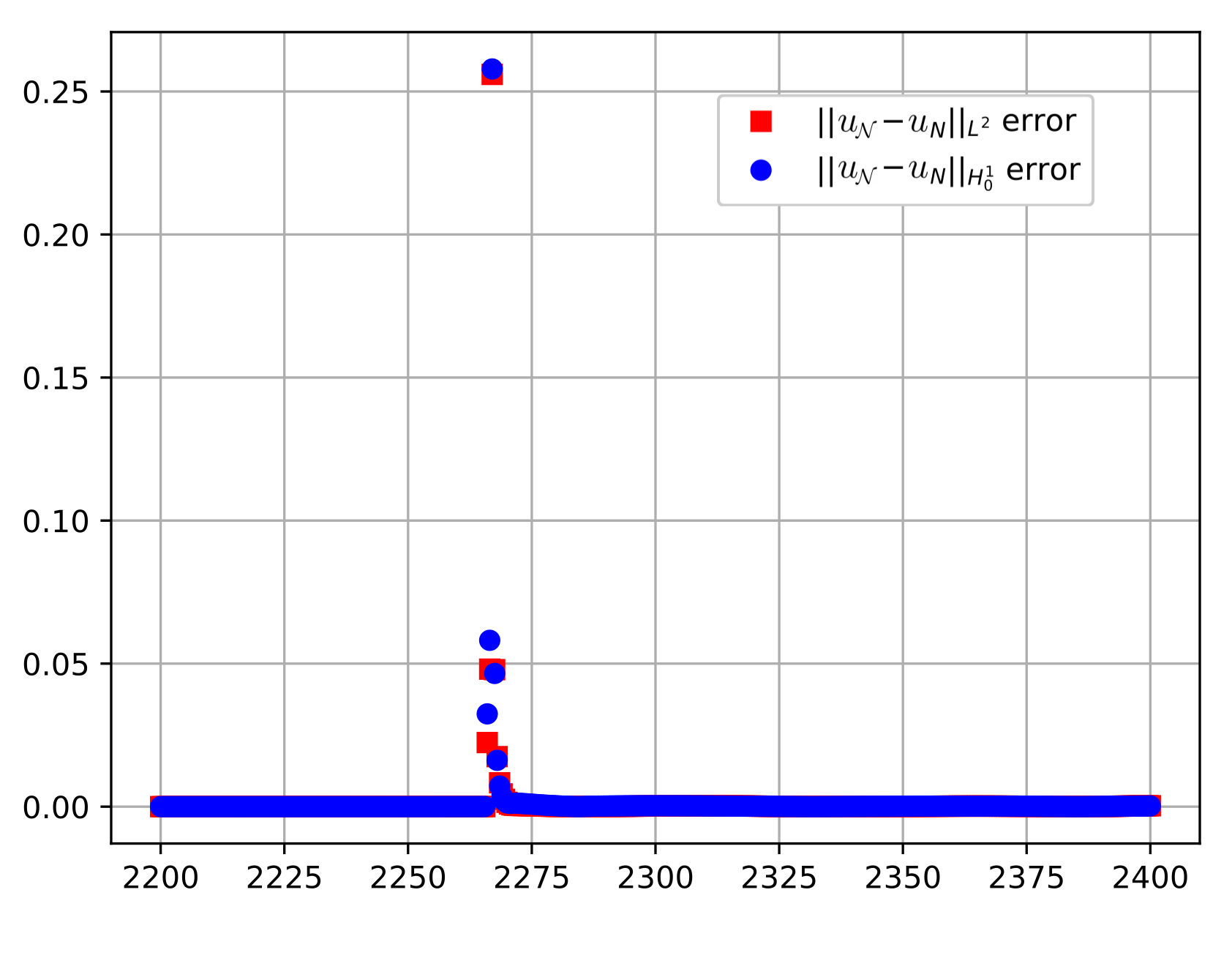}
         \put(-95,4){\makebox(0,0){$\mu$}}
         \caption{$N_1 = 4$}
         \label{fig:10_comp_a}
     \end{subfigure}
     \hfill
     \begin{subfigure}[b]{0.49\textwidth}
         \centering
         \includegraphics[width=\textwidth]{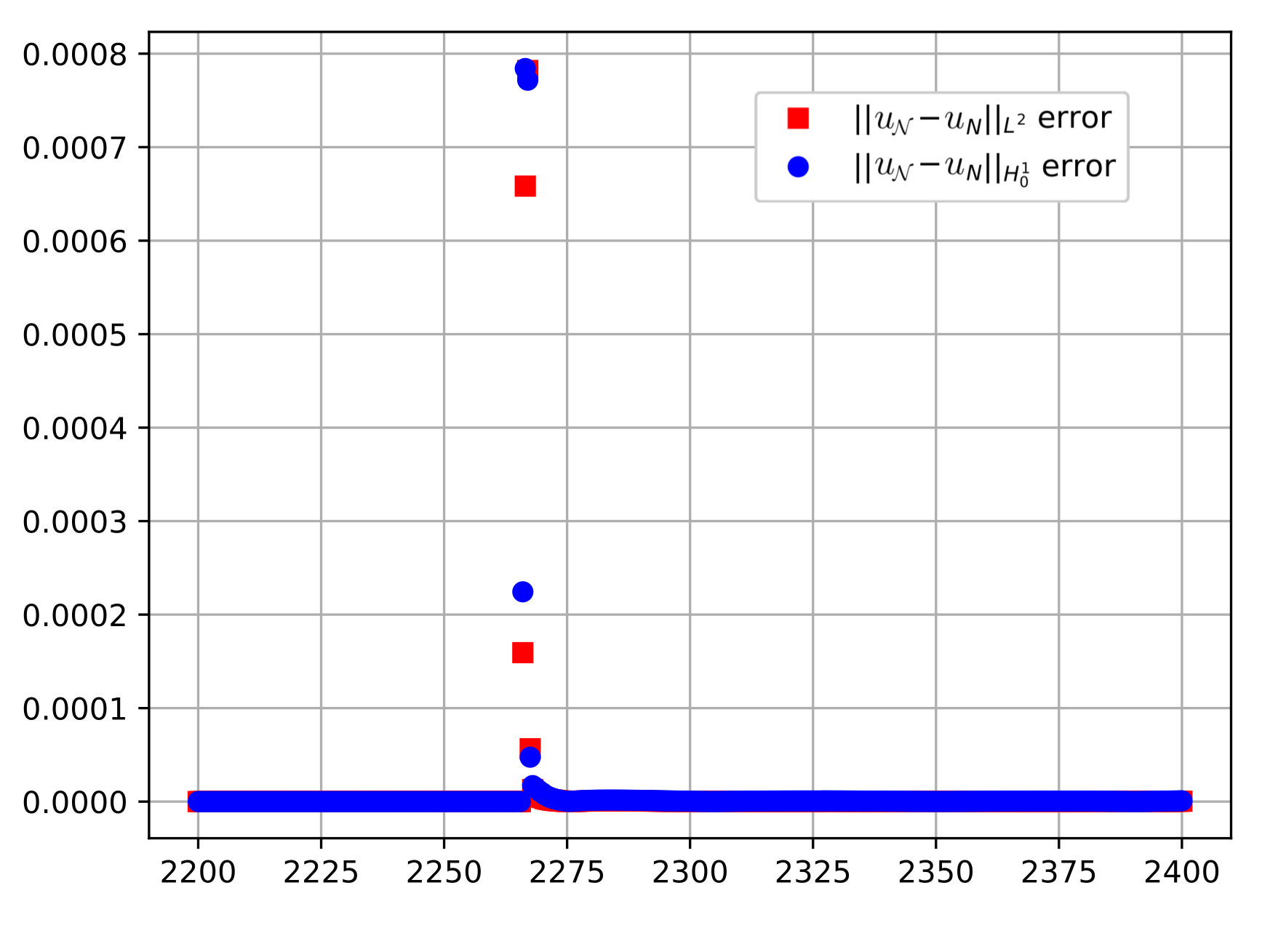}
          \put(-93,4){\makebox(0,0){$\mu$}}
         \caption{$N_2 = 5$}
         \label{fig:10_comp_b}
     \end{subfigure}
	\caption{Reduced basis errors for the SVK beam with a different number of basis functions.}
\label{fig:10_comp}
\end{figure}

Now we come back to the comparison between the two ROMs, employing respectively $N_1 = 4$ and $N_2 = 5$ basis. As we can see from Figure \ref{fig:10_comp}, the POD of dimension $N_1=4$ (left) produces a not completely satisfactory reduced error at the buckling point $\mu^*$ of order $2.e-1$, although it corresponds to a quite low tolerance $\epsilon_{POD} = 10^{-8}$. Due to the expected exponential convergence of the RB method, adding a single basis function, selected by the tolerance $\epsilon_{POD} = 10^{-10}$, the maximum error diminishes of almost three orders. Of course, in both cases, the average errors of order below $10^{-7}$ on $\Pa$ confirm the good approximation property of the POD.
Once again, the peaks of the error are reached where the sensitivity loses its differentiability.

The speed-up is still very low, almost $1.2$, but also the computation of high fidelity version of the bifurcation diagram in Figure \ref{fig_1_force} was less costly, indeed to build it we spent $t_{HF} = 180$(s), while the reduced order one required $t_{RB} = 150$(s). 

\subsubsection{Multi-parameter test case}

Here, we want to extend the investigation done in Section \ref{sec:hyp_dir} by modelling different materials through the parametrization of the elasticity constants. In particular, staying within the physical parametrization (we will see in the next section an example of the geometrical one), we want to understand how different materials behave with respect to an increasing compression imposed by means of Dirichlet boundary condition. Hence, we can define the multi-parameter $\bmu \in \Pa \subset \mathbb{R}^3$  as the triplet $\bmu = (\mu, E, \nu)$, where $\mu$  is the bifurcating parameter which controls the Dirichlet compression and $E$, $\nu$ are respectively the Young modulus and the Poisson ratio. 
We considered the SVK model with gravitational body force $B = (0, -1000)$, and $s(u) = \norm{u_y}_{\infty}$ as the output functional for the buckling detection.

For the analysis, we chose a parameter space given by $\Pa = [0.0, 0.2] \times [10^5, 10^7] \times [0.25, 0.42]$.
For what concerns the offline phase, we sampled $N_{train} = 1000$ snapshots for equispaced values of $\mu$, for each one of the 4 vertexes of the bi-dimensional parameter space of the elasticity constants, i.e.\ $[10^5, 10^7] \times [0.25, 0.42]$. 

Hence, as in the multi-parameter test cases in the previous chapters, also here we have sampled the bifurcation diagrams for different physical configuration, and then we have used the POD compression to extract the modes needed to recover the buckling of beams made by different materials.  

\begin{figure}[h]
\centering
\includegraphics[width=0.7\textwidth]{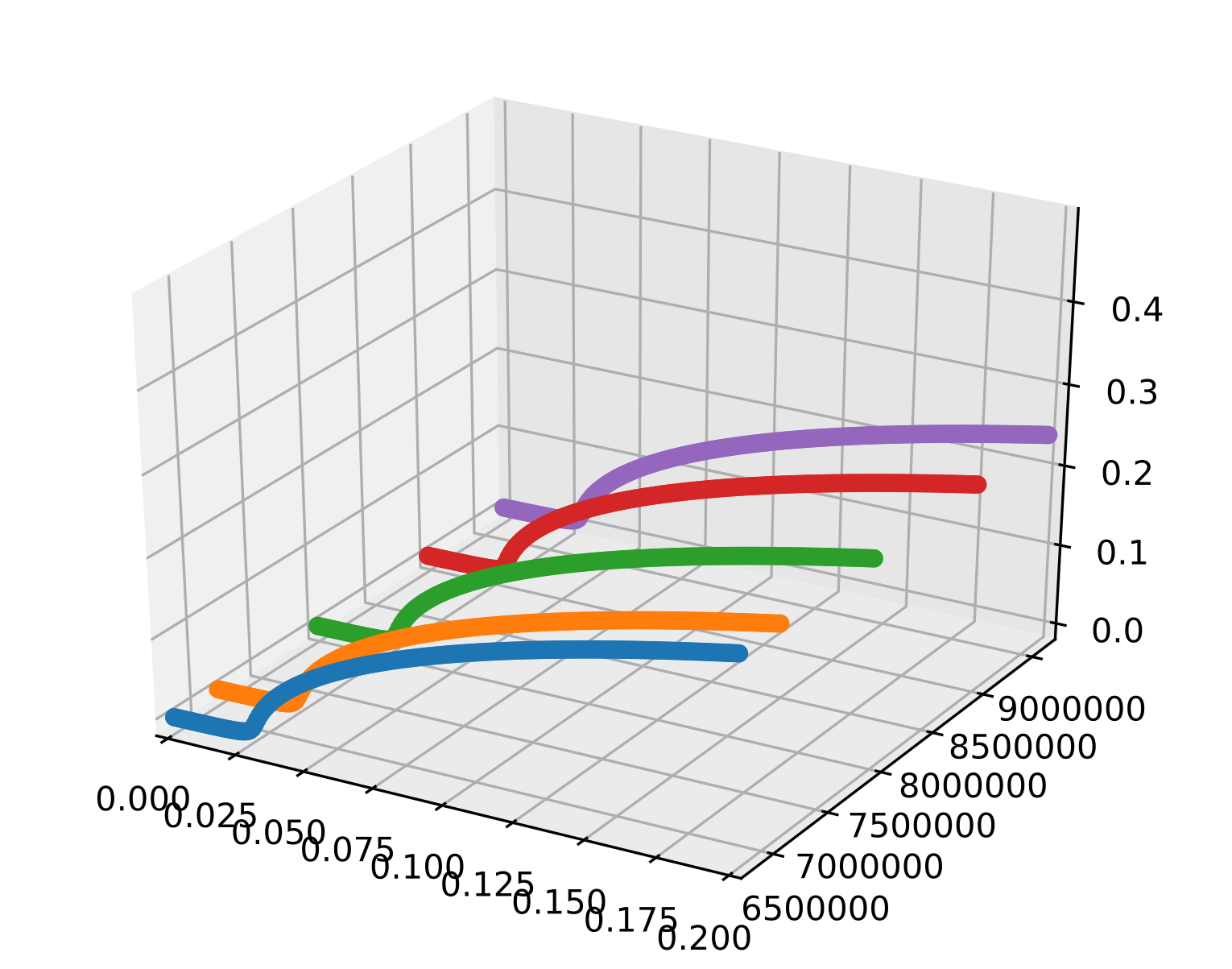}
    \put(-5,115){$\norm{u_y}_{\infty}$}
    \put(-45,23){$E$}
    \put(-220,12){$\mu$}
\caption{3D reduced bifurcation plot  for SVK beam with $B = (0, -1000)$ for five random pairs $(E, \nu) \in [10^5, 10^7]  \times [0.25, 0.42]$.}
\label{fig_20}
\end{figure}

Since the parameter space is 3-D, we expected that a much greater number of basis functions are needed, indeed using a POD tolerance $\epsilon_{POD} = 10^{-8}$ we obtained a reduced basis space of dimension $N = 43$.
The online continuation method was based on $K = 201$ equispaced $\mu$ values in $\Pa$, which corresponds to a continuation step $\Delta \mu = 10^{-3}$, and the bifurcation diagram was depicted for 5 random pairs of the elasticity constants $(E, \nu) \in [10^5, 10^7]  \times [0.25, 0.42]$.
Therefore, let us show in Figure \ref{fig_20} the reduced basis bifurcation diagrams with respect to the Young modulus.

\begin{figure}[h!]
\centering
     \begin{subfigure}[b]{0.49\textwidth}
         \centering
         \includegraphics[width=\textwidth]{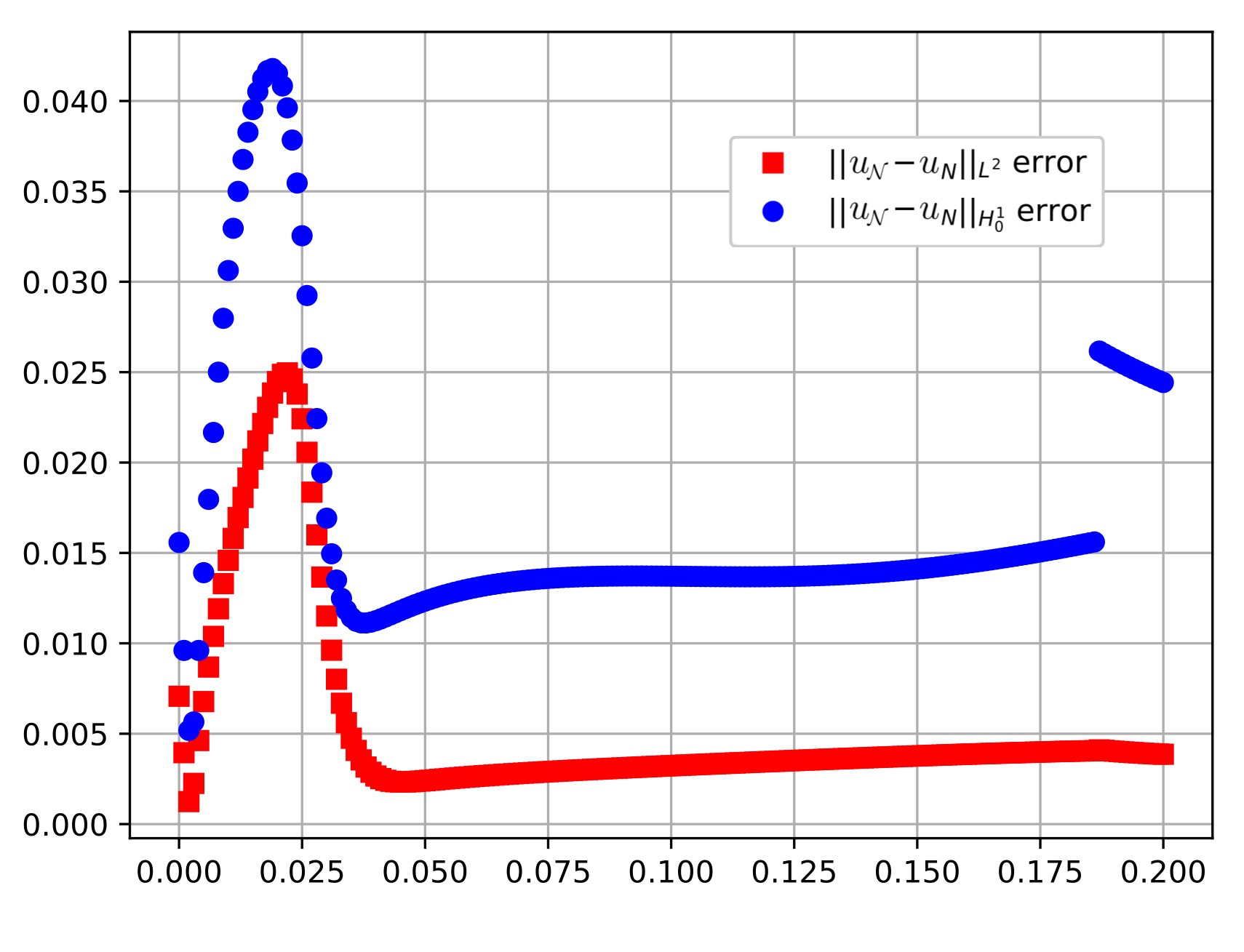}
                  \put(-95,4){\makebox(0,0){$\mu$}}
         \caption{$(E_1, \nu_1) = (5.2\cdot 10^5, 0.352)$}
         \label{fig:10_a_multi}
     \end{subfigure}
     \hfill
     \begin{subfigure}[b]{0.49\textwidth}
         \centering
         \includegraphics[width=\textwidth]{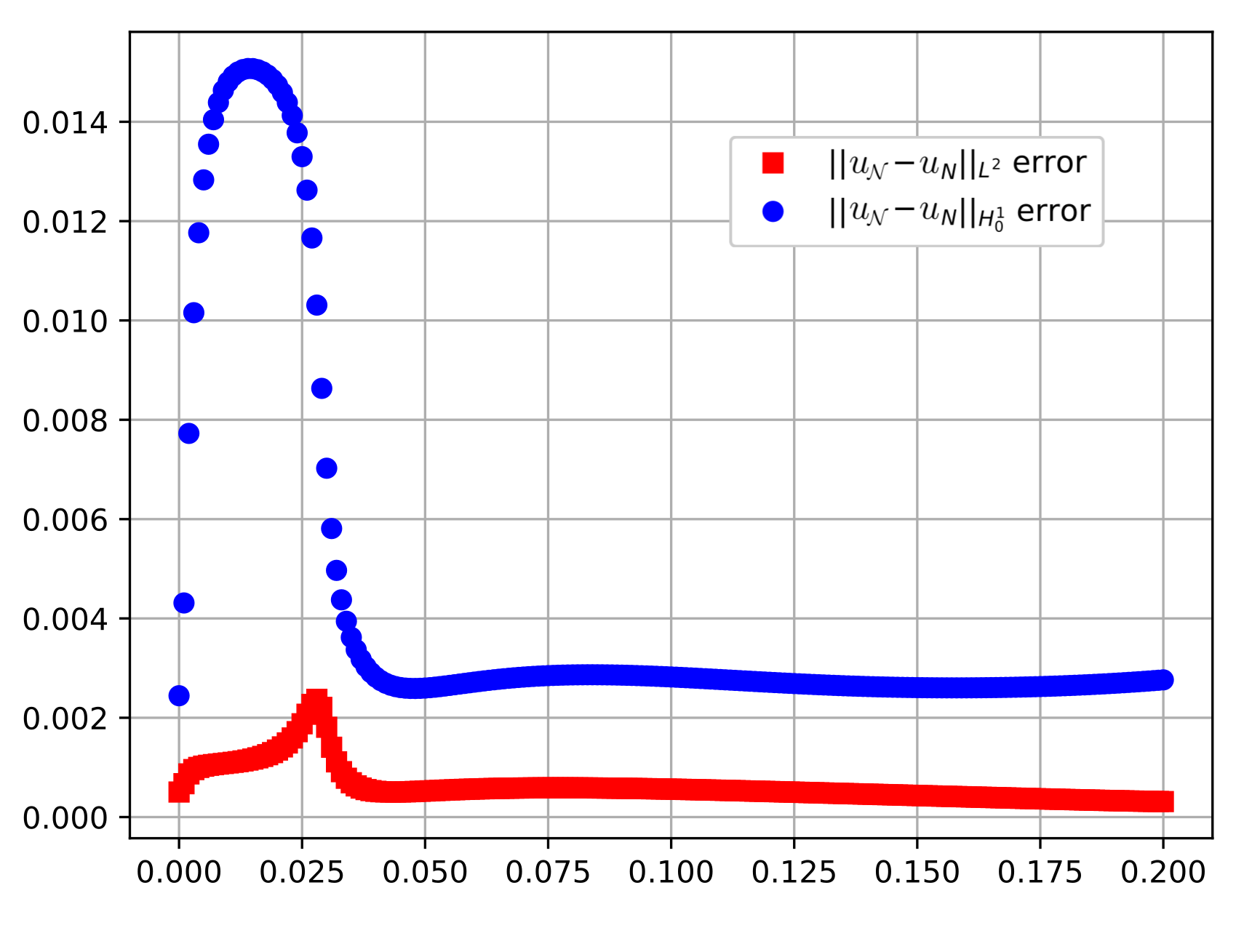}
                  \put(-95,4){\makebox(0,0){$\mu$}}
         \caption{$(E_2, \nu_2) = (2.6 \cdot 10^6, 0.272)$}
         \label{fig:10_b_multi}
     \end{subfigure}
	\caption{Reduced basis errors for the SVK beam with $B = (0, -1000)$ for fixed pairs $(E, \nu)$.}
\label{fig_10_multi}
\end{figure}

As we can see the reduced model was able to reconstruct the post-buckling behavior for all the pairs. As an example we show in Figure \ref{fig_10_multi} the reduced errors for two other pairs $(E_1, \nu_1) = (5.2\cdot 10^5, 0.352)$ and $(E_2, \nu_2) = (2.6 \cdot 10^6, 0.272)$.
We remark that the choice of the non-trivial body force was takes in order to avoid sharp gradients in the bifurcation diagrams and thus larger errors in the reduced approximation.

\begin{figure}[h]
\centering
\includegraphics[width=0.7\textwidth]{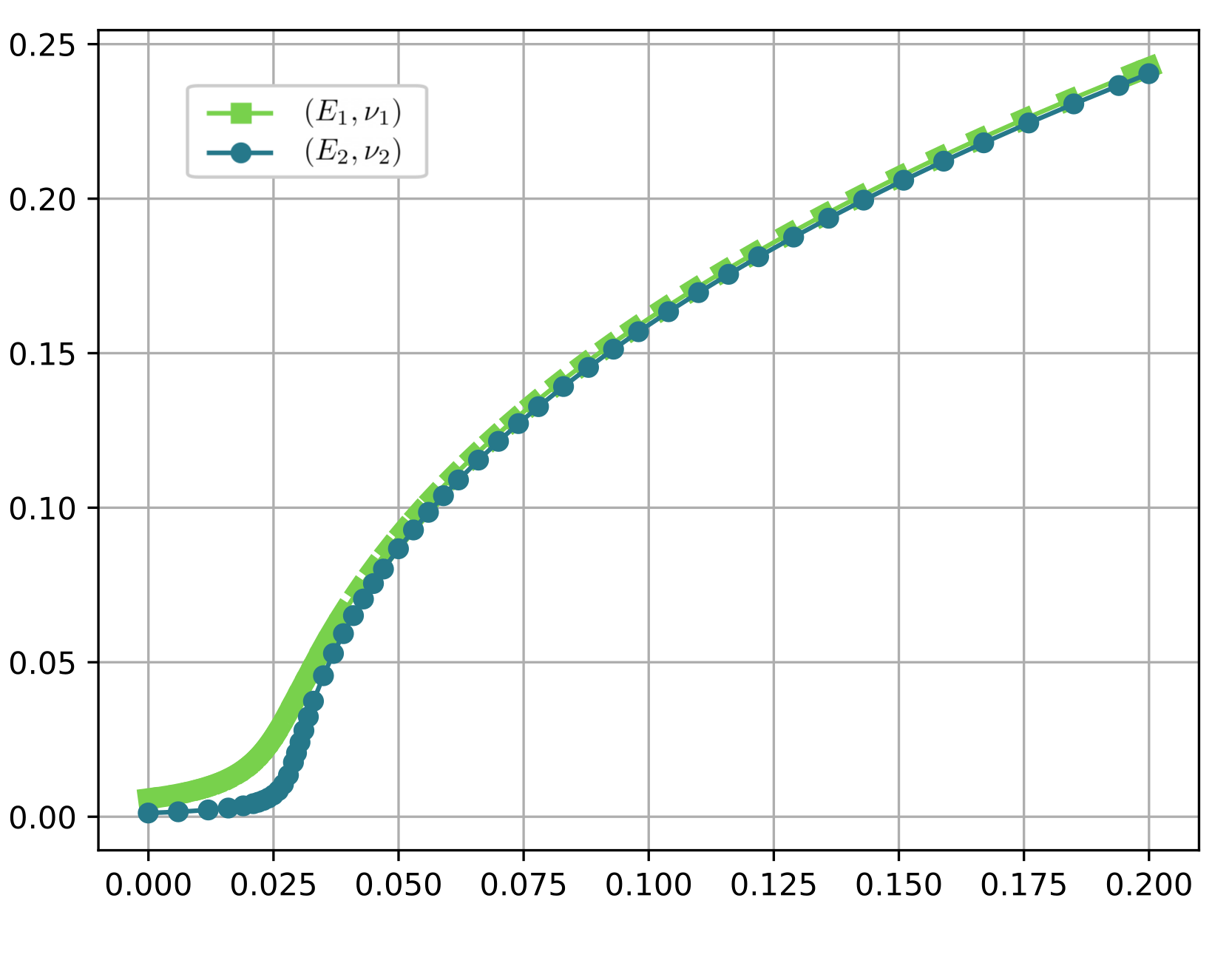}
\put(-140,2){\makebox(0,0){$\mu$}}
\put(-330,125){$\norm{u_y}_{\infty}$}
\caption{Reduced bifurcation diagram for SVK beam with $B = (0, -1000)$ for $(E_1, \nu_1)  =  (5.2\cdot 10^5, 0.352)$ and $(E_2, \nu_2)  =  (2.6 \cdot 10^6, 0.272)$.}
\label{fig_20_comp}
\end{figure}

Despite the displacement solution for the materials considered during the online phase shows similar properties in the considered parameter range, we can observe from Figure \ref{fig_20_comp} that a slight modification of the buckling point can occur. Moreover, we remark that for some choices of the elasticity constants, the body force can become irrelevant leading again to sharp discontinuity in the sensitivity. 

As regard the computational speed-up, given the increased dimension of the reduced space, the plot of the reduced diagram costs exactly as the high fidelity version $t_{HF} \approx t_{RB} = 900$(s).
Hence, the need for an empirical interpolation approach capable of decoupling the online phase from the high fidelity degrees of freedom is still more evident. Despite this, the interesting point here is that, even through a naive approach for the sampling of the elasticity parameters, we were able to detect the buckling and the related post-buckling behavior for a wide range of a 3-D parameter space, by means of $N = 43$ basis functions.

\subsubsection{Geometrical parametrization}\label{ssec:geom_param}
In this test case, we will consider for the first time a parametrized geometry, trying to understand how this influences the buckling property of the beam. 
Thus, let us gives brief introduction to the geometrical parametrization test case \cite{quarteroniReducedBasisMethods2016}, which we will use also in subsequent sections.

Given a parameter $\bmu \in \Pa$, we can distinguish between the physical quantities $\bmu_p \in \Pa_p$ (compression, forces, viscosity, trap strength) and geometrical quantities $\bmu_g \in \Pa_g$ (lengths, angles).  Hence, we now consider the parameter dependent domain $\widetilde{\Omega}(\bmu_g)$ and its reference configuration $\Omega$, related by the transformation map 
\begin{equation}\label{transformation_map}
\Phi : \Om \times \Pa_g \to \R^d \quad \text{such that} \quad \widetilde{\Omega}(\bmu_g) = \Phi(\Om; \bmu_g),  \quad \forall \bmu_g \in \Pa_g.
\end{equation}
The aim of such map is allow us to write a generic weak formulation posed on the parameter dependent domain $\widetilde{\Omega}(\bmu_g)$, in the reference configuration $\Omega$, in order to guarantee the assembly of the high fidelity $\bmu$-independent quantities during the offline phase.

The key point is the formula for the change of variables, which for any integrable functions $\widetilde{f}: \widetilde{\Om} \to \R^d$ is given by 
\begin{equation}
\label{eq:change_var}
\int_{\widetilde{\Om}(\bmu_g)}\widetilde{f} \ d\widetilde{\Om} = \int_{\Om} f \det(\mathsf{J}_{\Phi}) \ d\Om ,
\end{equation}
where $f = \widetilde{f} \circ \Phi$  and $\mathsf{J}_{\Phi}$ is the Jacobian of the transformation map $\Phi$.
When the forms involve spatial derivatives, one relies on the chain rule to obtain a formula that encodes the parametrization dependence. As an example we can consider the standard $H^1$ semi-norm, which can be transformed as

\begin{equation}
\label{eq:change_var_semi}
\int_{\widetilde{\Om}(\bmu_g)}\nabla_{\widetilde{x}}\widetilde{f} : \nabla_{\widetilde{x}}\widetilde{h}  \ d\widetilde{\Om} = \int_{\Om} (\nabla_x f) \mathsf{K} : \nabla_x h \ d\Om ,
\end{equation}
where the parametrized tensor $\mathsf{K}: \R^d \times \Pa \to \R^{d \times d}$ is defined as $$\mathsf{K}(x; \bmu) = \mathsf{J}^{-1}_{\Phi}(x;\bmu)\mathsf{J}^{-T}_{\Phi}(x;\bmu)\det(\mathsf{J}_{\Phi}(x;\bmu)) .$$ 

\begin{figure}[h]
\centering
\includegraphics[width=0.7\textwidth]{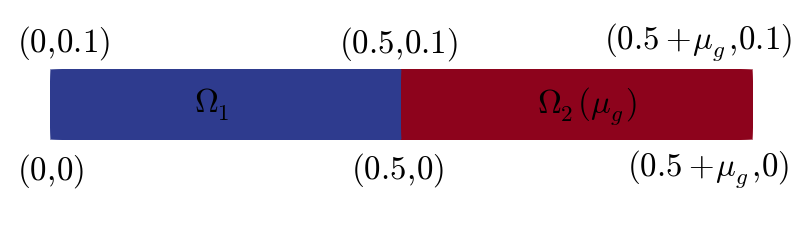}
\caption{2-D beam with parametrized geometry and $\mu_g \in [0.5, 1]$.}
\label{fig_17_domain}
\end{figure}

Having presented the basic notions on geometrical parametrization, we can go back to our 2-D toy problem to investigate how the length of the beam is related to the buckling point.
For this reason we now consider the geometrically parametrized beam, where its semi-length $\mu_g$ is added as new parameter.
Therefore, the domain depicted in Figure \ref{fig_17_domain} can be expressed as $\widetilde{\Omega}(\mu_g) = \widetilde{\Omega}_1 \cup \widetilde{\Omega}_2(\mu_g)$, where $\widetilde{\Omega}_1 = [0, 0.5] \times [0, 0.1]$ and $\widetilde{\Omega}_2(\mu_g) = [0.5, 0.5 + \mu_g] \times [0, 0.1]$.
In this case the transformation map is simply given by the affine function $$\Phi(x; \mu_g) =  \begin{bmatrix} 2\mu_g(x-0.5) + 0.5 \\ y \end{bmatrix} \quad \text{for } x \in \Omega_2 =  \widetilde{\Omega}_2(0.5) .$$

As a consequence of the geometrical parametrization, a consistent number of terms is involved in the affine (w.r.t.\ the parameter) decomposition of the weak formulation. As an example, for the SVK model one obtains the split of \eqref{eq:weakgal_hyp} in 3 linear and 5 nonlinear terms, already for this simple geometry change.

\begin{figure}[h]
\centering
\includegraphics[width=0.7\textwidth]{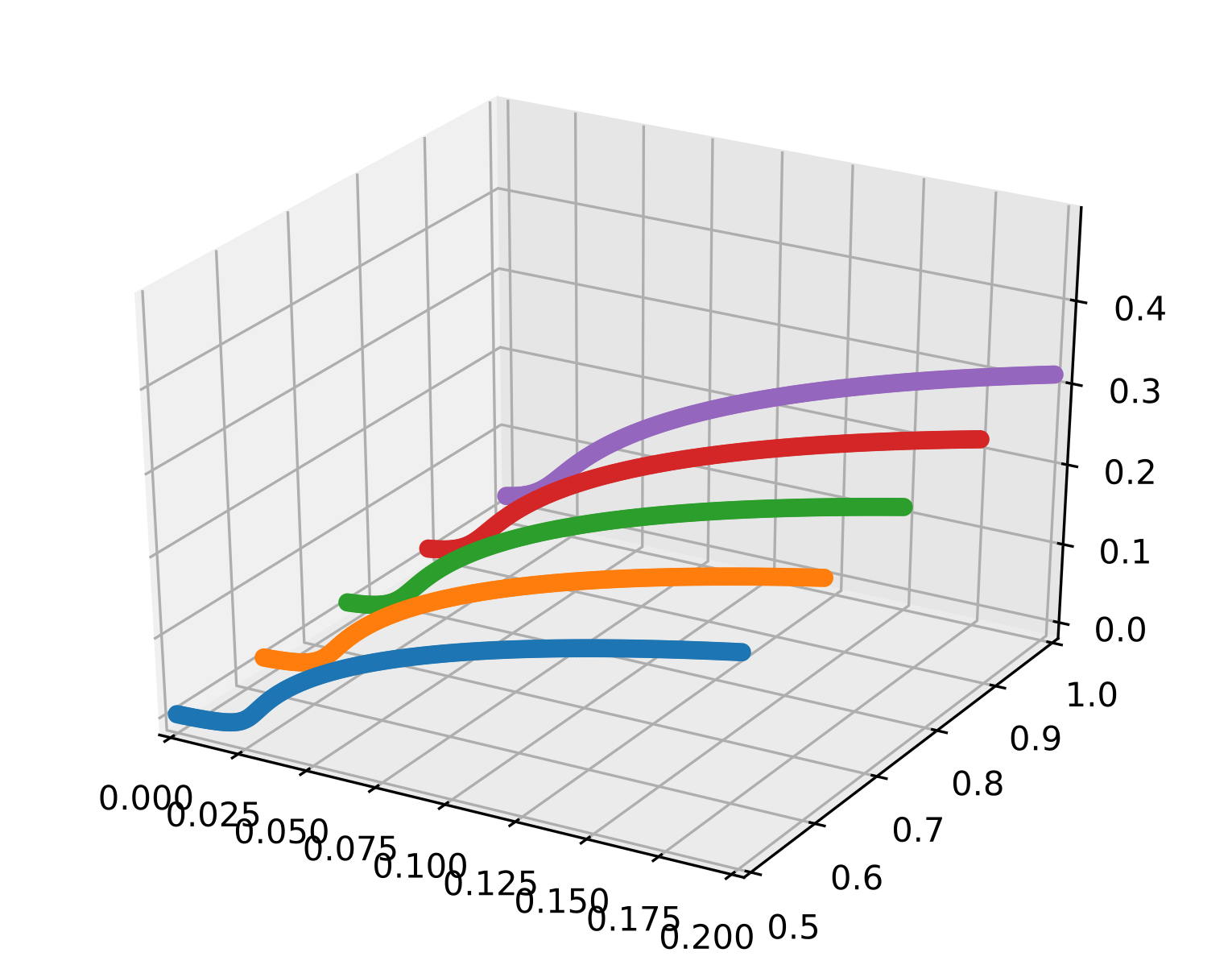}
    \put(-5,115){$\norm{u_y}_{\infty}$}
    \put(-50,28){$\mu_g$}
    \put(-220,12){$\mu$}
\caption{3D bifurcation plot for SVK beam with $B = (0, -1000)$ with $\mu_g \in \Pa_g$.}
\label{fig_17}
\end{figure}

Regarding the physical setting, we consider the SVK model compressed through Dirichlet boundary conditions with null traction force $T = 0$, Young modulus $E = 10^6$, Poisson ratio $\nu = 0.3$ and gravitational body force $B = (0, -1000)$.
During the offline phase we sampled $N_{train} = 1000$ points in the parameter space $\Pa_p = [0, 0.2]$, for each one of the 3 equispaced points in $\Pa_g = [0.5, 1]$. Using, as always, a POD global compression with tolerance $\epsilon_{POD} = 10^{-8}$ we obtained a reduced basis space of dimension $N = 12$. 
On the contrary, the online continuation method to reconstruct the 3-D bifurcation diagram in Figure \ref{fig_17} we chose $K = 201$ equispaced points in $\Pa_p$, for 5 equispaced values of the semi-length $\mu_g \in \Pa_g$. 

\begin{figure}[h!]
\centering
    \begin{subfigure}[b]{0.49\textwidth}
         \centering
         \includegraphics[width=\textwidth]{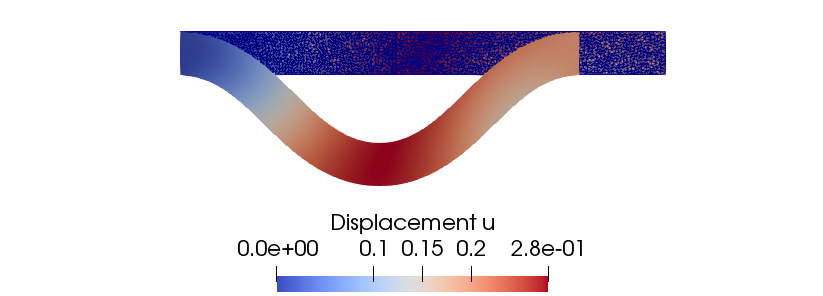}
         \caption{$\mu_g = 0.625$}
         \label{fig:6_a_geom_1}
     \end{subfigure}
     \hfill
    \begin{subfigure}[b]{0.49\textwidth}
         \centering
         \includegraphics[width=\textwidth]{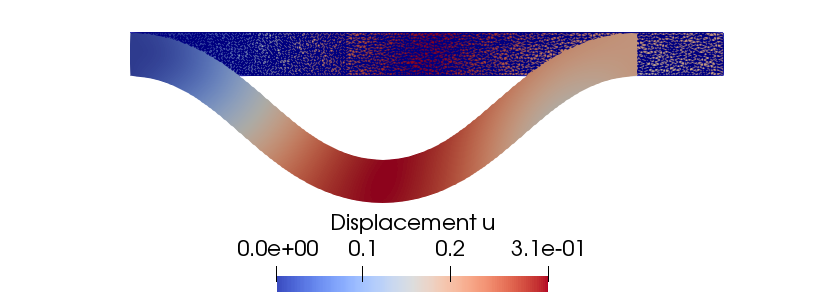}
         \caption{$\mu_g = 0.875$}
         \label{fig:6_b_geom_1}
     \end{subfigure}
\caption{High fidelity displacement $u$ for the SVK beam with $B = (0, -1000)$ at $\mu = 0.2$ for different geometries.}
\label{fig_6_geom}
\end{figure}

Two representative solutions of the post-buckling branches for $\mu_g = 0.625$ and $\mu_g = 0.875$ are depicted in Figure \ref{fig_6_geom}, respectively left and right, at $\mu = 0.2$.

\begin{figure}[h!]
\centering
     \begin{subfigure}[b]{0.49\textwidth}
         \centering
         \includegraphics[width=\textwidth]{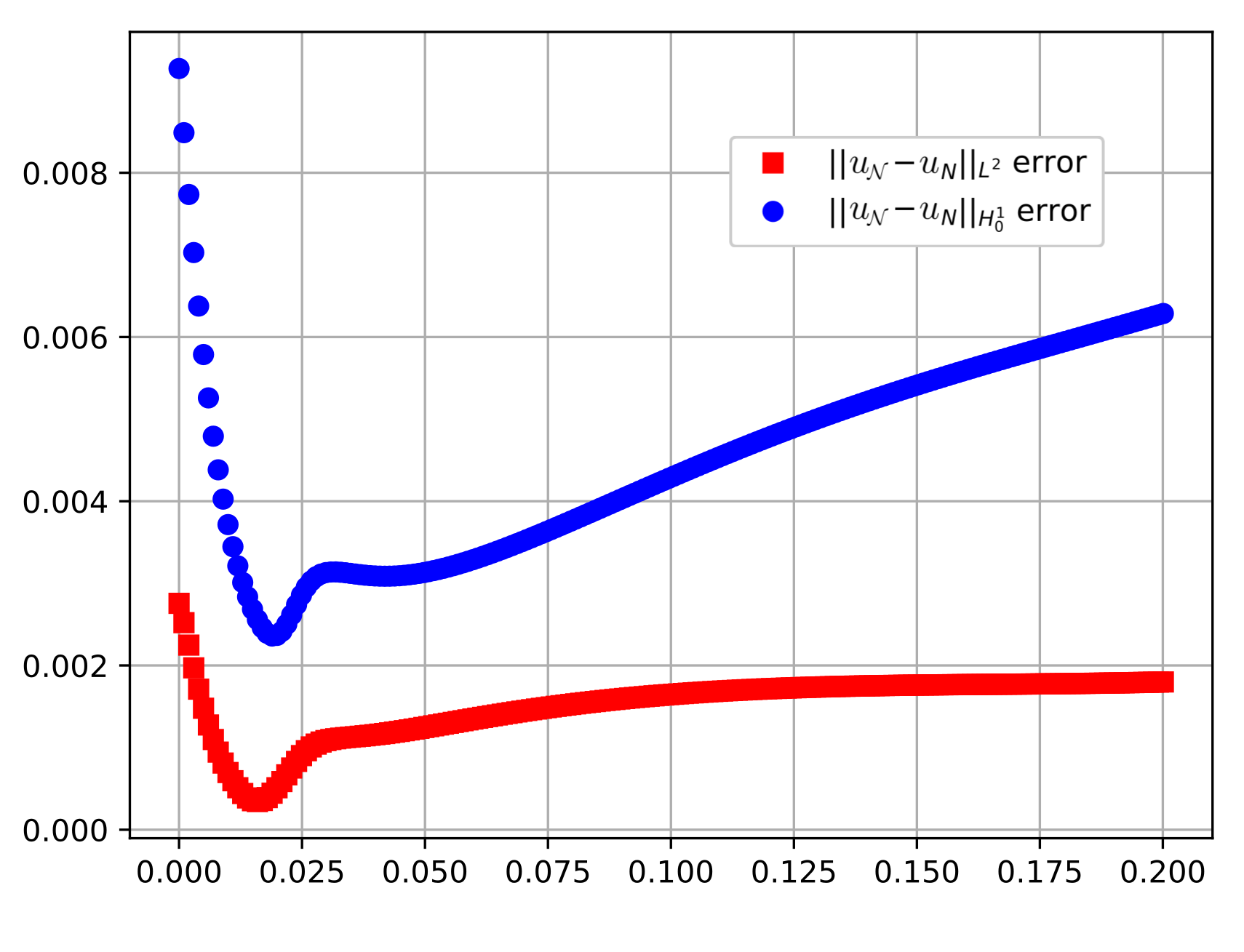}
                           \put(-95,4){\makebox(0,0){$\mu$}}
         \caption{$\mu_g = 0.625$}
         \label{fig:10_a_geom}
     \end{subfigure}
     \hfill
     \begin{subfigure}[b]{0.49\textwidth}
         \centering
         \includegraphics[width=\textwidth]{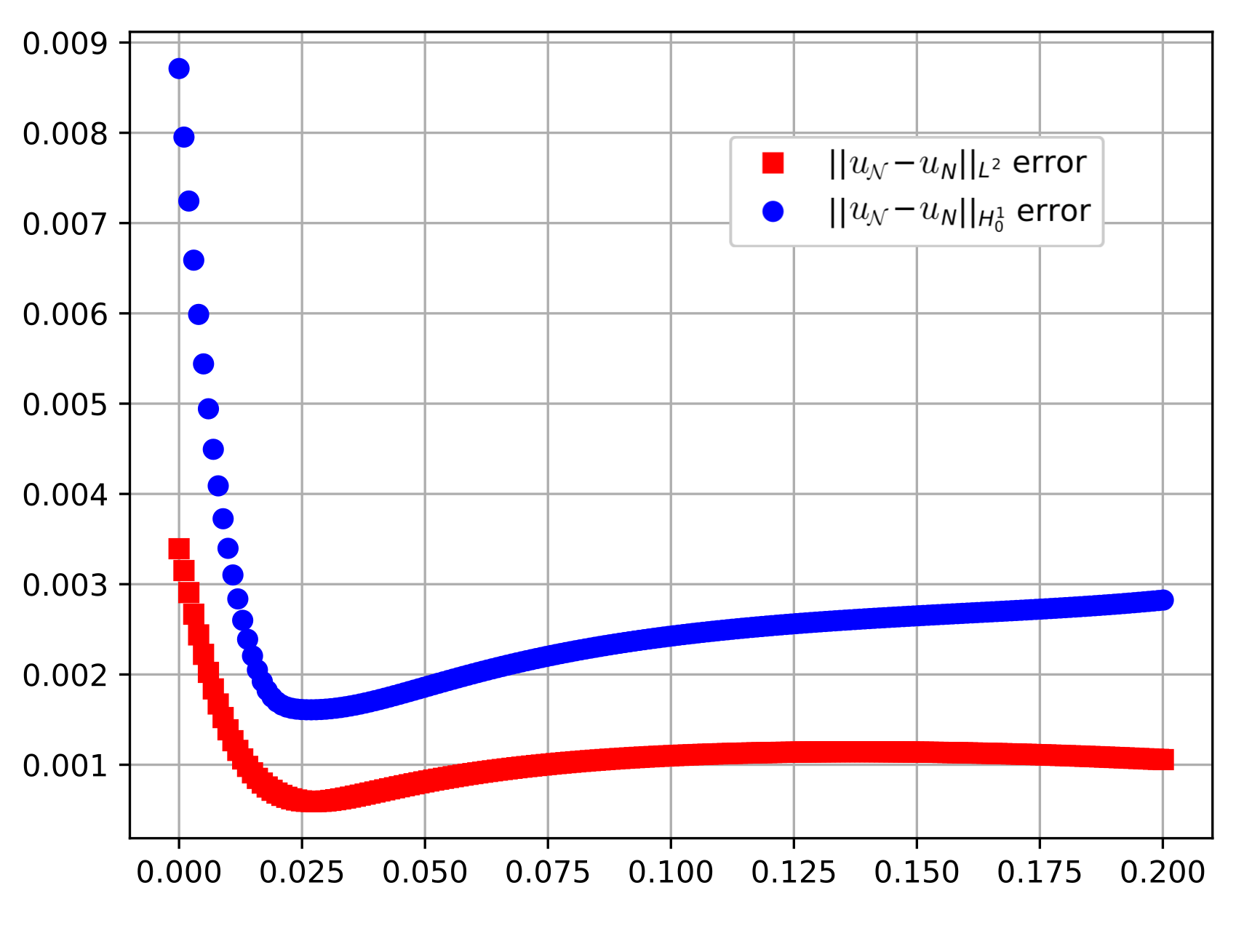}
                           \put(-95,4){\makebox(0,0){$\mu$}}
         \caption{$\mu_g = 0.875$}
         \label{fig:10_b_geom}
     \end{subfigure}
	\caption{Reduced basis errors for the SVK beam with $B = (0, -1000)$ for fixed geometries.}
\label{fig_10_geom}
\end{figure}

As we can see from Figure \ref{fig_10_geom} the reduced manifold was able to approximate the buckling also for not sampled geometries with good approximation accuracy. 
This way, we can effectively study the evolution of the buckling varying the length of the beam.

\begin{figure}[h]
\centering
\includegraphics[width=0.7\textwidth]{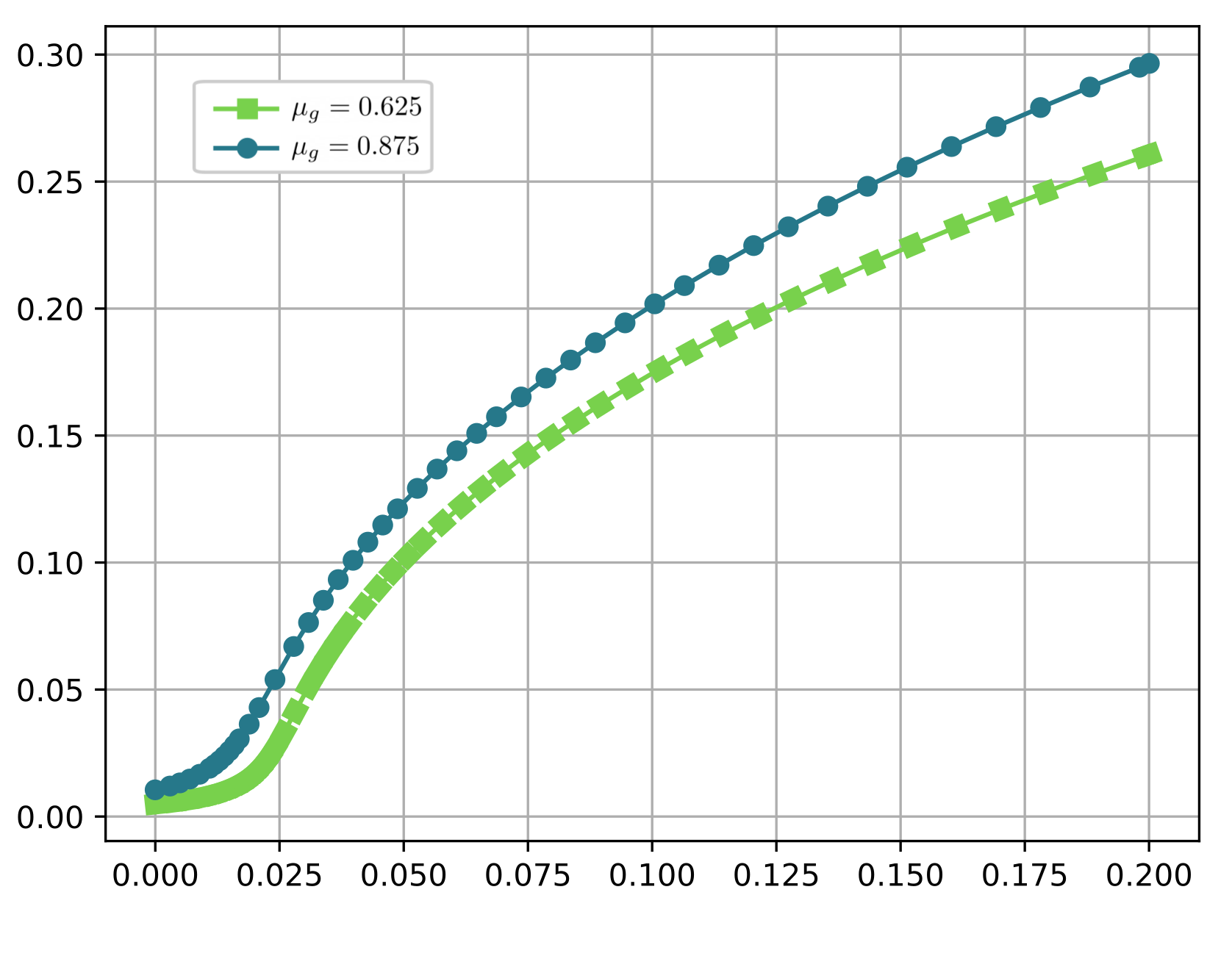}
\put(-137,4){\makebox(0,0){$\mu$}}
\put(-330,125){$\norm{u_y}_{\infty}$}
\caption{Reduced bifurcation diagram for the SVK beam with $B = (0, -1000)$ for $\mu_g = 0.625$ and $\mu_g = 0.875$.}
\label{fig_20_comp_geom}
\end{figure}

As we expected when the length of the beam varies, also the buckling point changes its position. Indeed, as we can observe from Figure \ref{fig_20_comp_geom} the longer is the beam the sooner it buckles. This is still more evident in the 2-D projected diagram where the branches correspondent to $\mu_g = 0.625$ and $\mu_g = 0.875$ are plotted.
Indeed, the buckling point for the beam with corresponding semi-length $\mu_g = 0.875$ buckles for a value of $\mu^* = 0.016$, while for the configuration with $\mu_g = 0.625$ it occurs at $\mu^* = 0.022$.

Once again the speed-up is essentially null with $t_{HF} \approx t_{RB} = 3000$(s), in fact the lower number of basis function is balanced by the increased number of terms involved in the weak formulation, which due to the lack of empirical interpolation strategies require a costly projection on the high fidelity dimension space.

\subsection{3-D toy problem}\label{sec:3d_toy}
Having analyzed in the previous section a variety of settings for the buckling of a two-dimensional beam, the main step towards real applications is to consider a 3-D geometry, as an extension of the one studied previously. 
Following the test case investigated in \cite{zanon:_phd_thesis}, we fix the domain as $\Omega = [0, 1] \times [0, 0.2]
\times [0, 0.079]$, we choose the body forces as either $B = (0, 0, 0)$ or $B = (0, 0, -1000)$ with trivial traction $T = 0$, and fixed material properties $E = 10^6$ and $\nu = 0.3$. 

Hence, we go back to the one parameter case, in which the compression parameter is imposed through Dirichlet boundary conditions. For this reason, we seek the solution displacement $u$ in the functional space $$\X = \lbrace{ u \in (H^1(\Omega))^3: u = (0, 0, 0) \ \text{on}\  \Gamma^{l}_D,\ u = (-\mu, 0, 0) \ \text{on} \  \Gamma^{r}_D \rbrace} ,$$ where $\Gamma^{l}_D = \lbrace{0\rbrace} \times [0, 0.2]
\times [0, 0.079]$ and $\Gamma^{r}_D = \lbrace{1\rbrace} \times [0, 0.2]
\times [0, 0.079]$.

We built the Finite Element space with $\mathbb{P}_1$ linear elements on a tetrahedral mesh, resulting in a high fidelity dimension $\N = 1734$.
As before, the parameter space is given by $\Pa = [0.0, 0.2]$, and its exploration is performed through the simple continuation method with fixed step $\Delta\mu = 2 \cdot 10^{-4}$, resulting in $N_{train} = 1000$ snapshots.
The reduced manifold was built choosing a tolerance $\epsilon_{POD} = 10^{-8}$ and the bifurcation diagram is recovered with an online continuation method based on $K = 2000$ equispaced points in $\Pa$.

We remark that since here we added a space dimension, the buckling for the rectangular cross-section beam can happen in both the directions individuated by the cross-section of the beam perpendicular to the compression axis. Despite this, the first buckling usually occurs in the direction of minimum length, hence in this case the $z$-axis.
For this reason, in order to detect the buckling behavior, we consider as output functional the infinite norm of the $z$-component of the displacement $\norm{u_z}_{\infty}$. 
  
For the three-dimensional beam, both constitutive relations, SVK and NH, have been investigated.
For what concerns the SVK model, we obtained a reduced basis space of dimension $N = 9$, for both the choices of the body force $B$. 

\begin{figure}[h!]
\centering
     \begin{subfigure}[b]{0.49\textwidth}
         \centering
         \includegraphics[width=\textwidth]{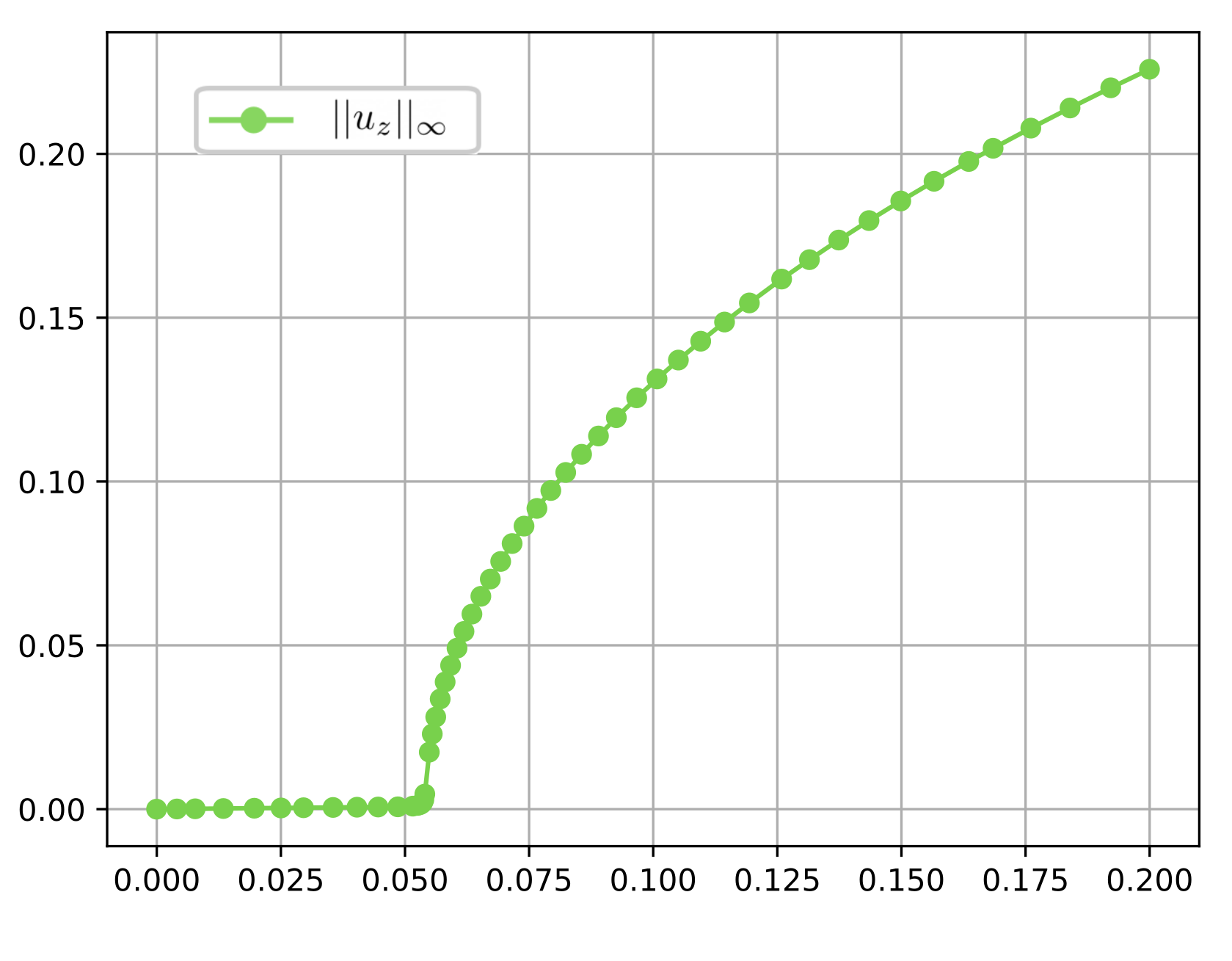}
                  \put(-95,4){\makebox(0,0){$\mu$}}
         \caption{$B = (0, 0, 0)$}
         \label{fig:15_a_3d}
     \end{subfigure}
     \hfill
     \begin{subfigure}[b]{0.49\textwidth}
         \centering
         \includegraphics[width=\textwidth]{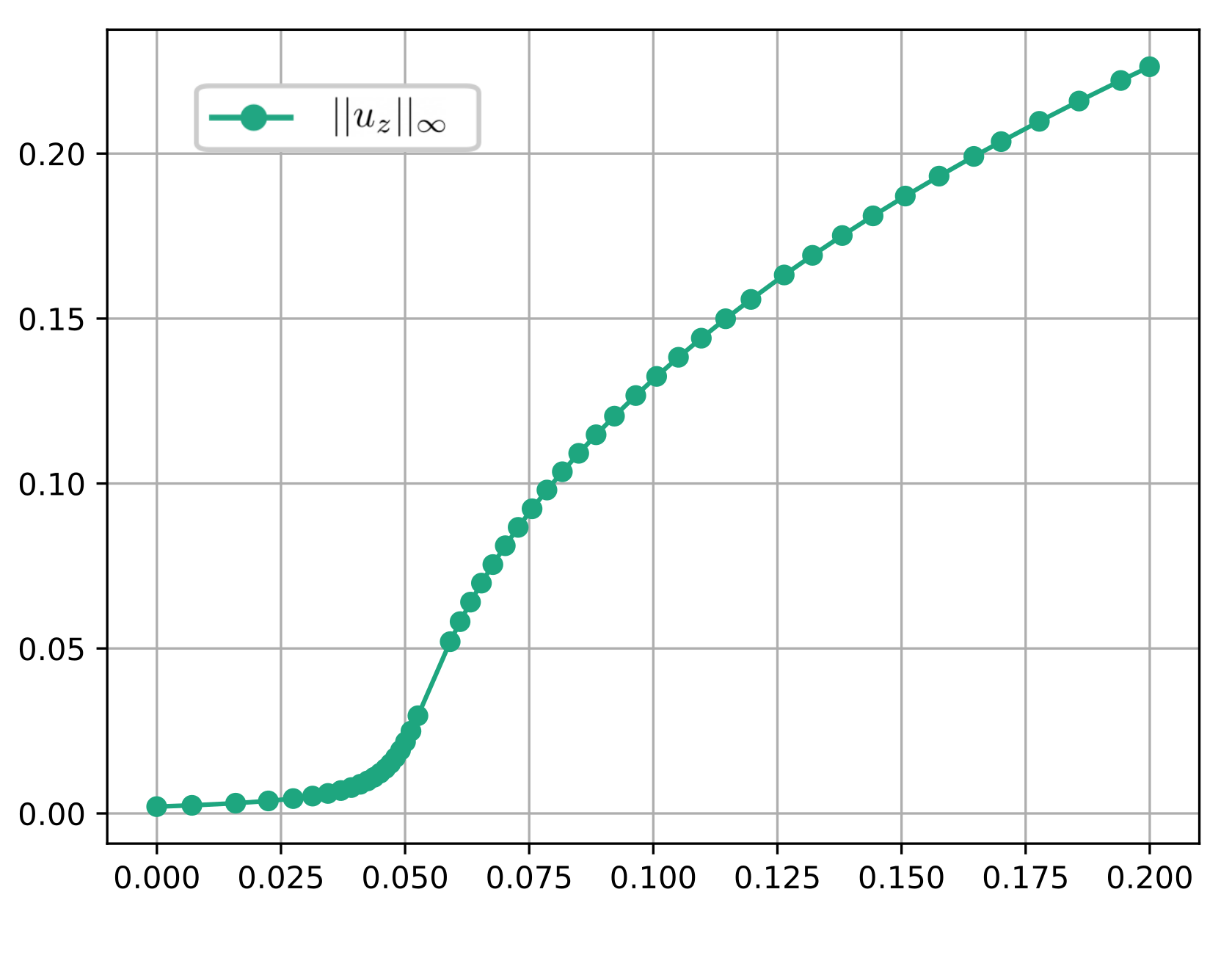}
                  \put(-95,4){\makebox(0,0){$\mu$}}
         \caption{$B = (0, 0, -1000)$}
         \label{fig:15_b_3d}
     \end{subfigure}
	\caption{Reduced basis bifurcation diagrams for the 3D SVK beam with respect to different body forces.}
\label{fig_15_3d}
\end{figure}

\begin{figure}[h]
\centering
\includegraphics[width=0.7\textwidth]{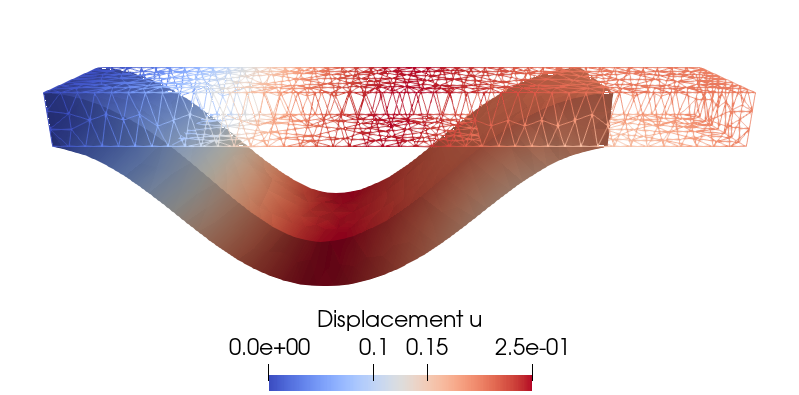}
\caption{High fidelity displacement $u$ for the 3D SVK beam with $B = (0, 0, 0)$ at $\mu = 0.2$.}
\label{fig_6_3d}
\end{figure}

In Figure \ref{fig_15_3d} we can see the bifurcation plot for trivial $B = (0, 0, 0)$ and gravitational $B = (0, 0, -1000)$ body forces. As we can see, also in this case the sharp gradient in the sensitivity was smoothed by the external force, but in Figure \ref{fig:15_a_3d} we can clearly observe that the buckling of the beam occurs for the value $\mu^* \approx 0.053$.
As predicted before the buckling occurs along the $z$-direction and a representative solution of the post-buckling branch is depicted in Figure \ref{fig_6_3d} for $\mu = 0.2$ with respect to the original un-deformed configuration (mesh wireframe). 

\begin{figure}[h!]
\centering
     \begin{subfigure}[b]{0.49\textwidth}
         \centering
         \includegraphics[width=\textwidth]{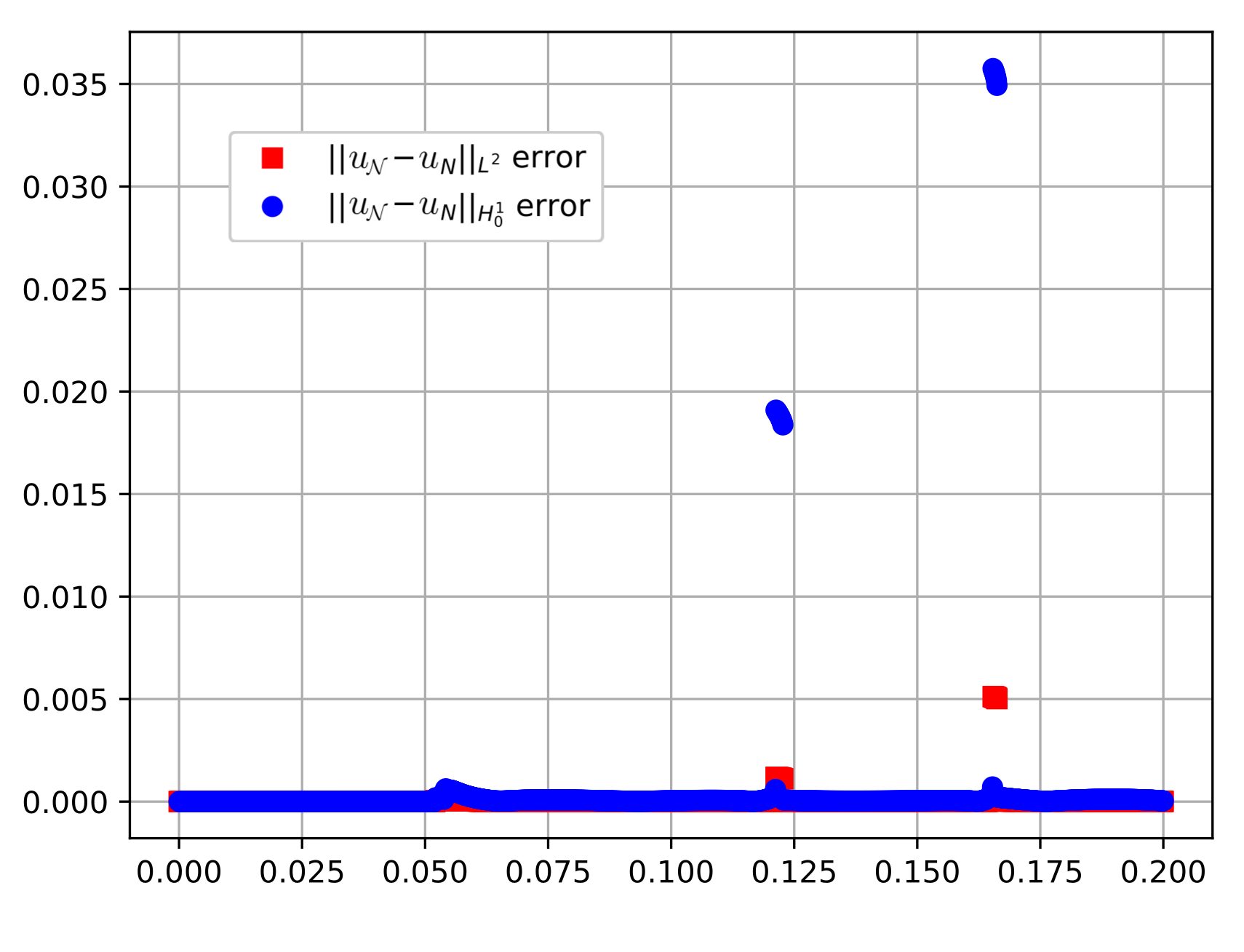}
                           \put(-95,4){\makebox(0,0){$\mu$}}
         \caption{$B = (0, 0, 0)$}
         \label{fig:10_a_3d}
     \end{subfigure}
     \hfill
     \begin{subfigure}[b]{0.49\textwidth}
         \centering
         \includegraphics[width=\textwidth]{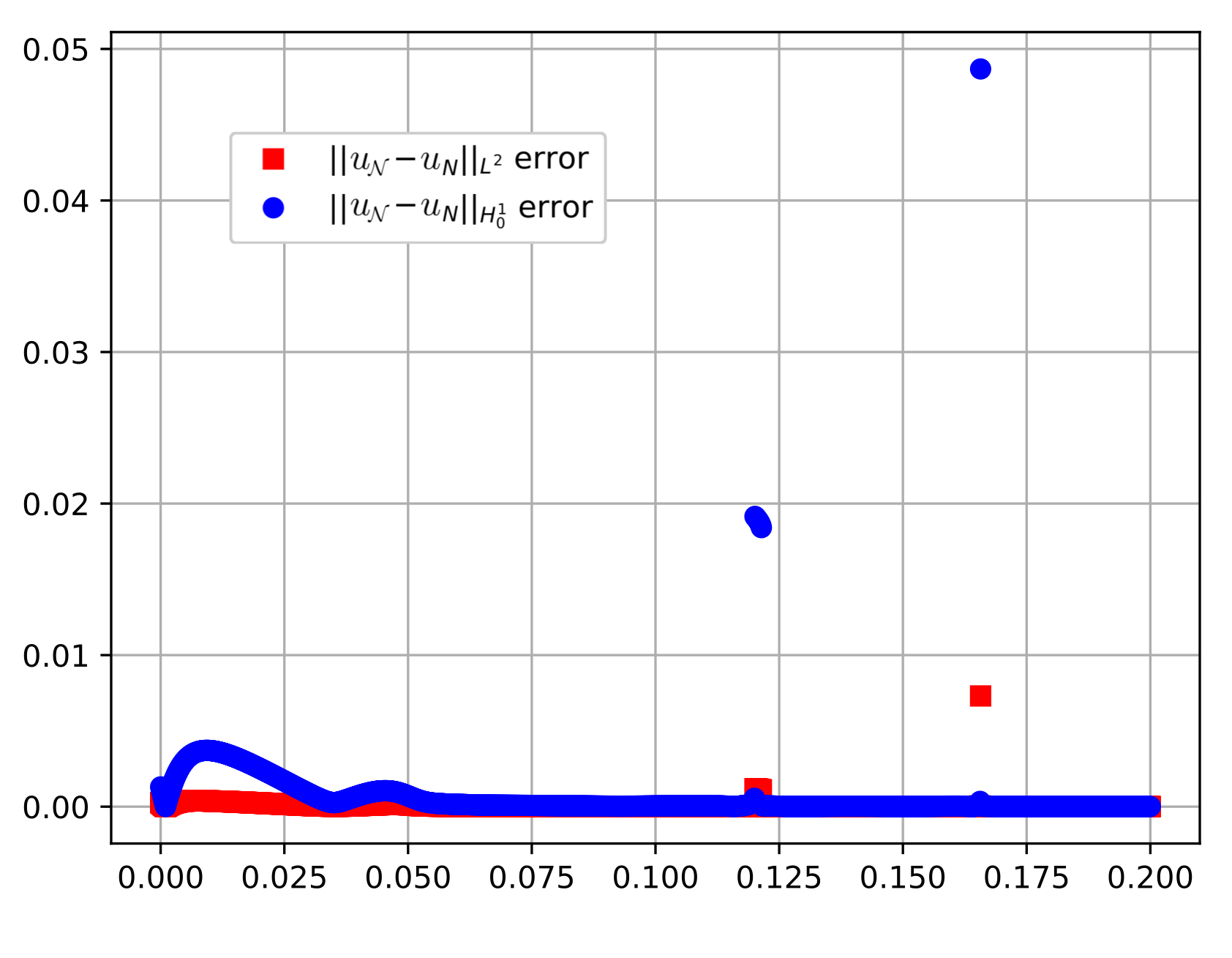}
                           \put(-95,4){\makebox(0,0){$\mu$}}
         \caption{$B = (0, 0, -1000)$}
         \label{fig:10_b_3d}
     \end{subfigure}
\caption{Reduced basis errors for the 3D SVK beam with respect to different body forces.}
\label{fig_10_3d}
\end{figure}

As we can see from Figure \ref{fig_10_3d}, the reduced order model has overall good approximation properties, since it is capable to reconstruct the post-buckling behavior with a maximum error of order $10^{-2}$. We remark that while the peaks in correspondence to the bifurcation point agree with the previous analysis considering different body forces, in this case the maximum errors were obtained for a bigger value of the compression. We did not investigate further the buckling, but from the considerations about the performances of the RB at bifurcation points, we guess that other buckling points are located near $\mu = 0.122$ and $\mu = 0.166$.

Once again, we observe the consistent difference between the maximum error and the average one of order $10^{-8}$, remarking that this is caused by the lack of regularity of the solution with respect to the parameter, due to the bifurcation phenomena. Same results as before hold for the speed-up.

\begin{figure}[h!]
\centering
     \begin{subfigure}[b]{0.49\textwidth}
         \centering
         \includegraphics[width=\textwidth]{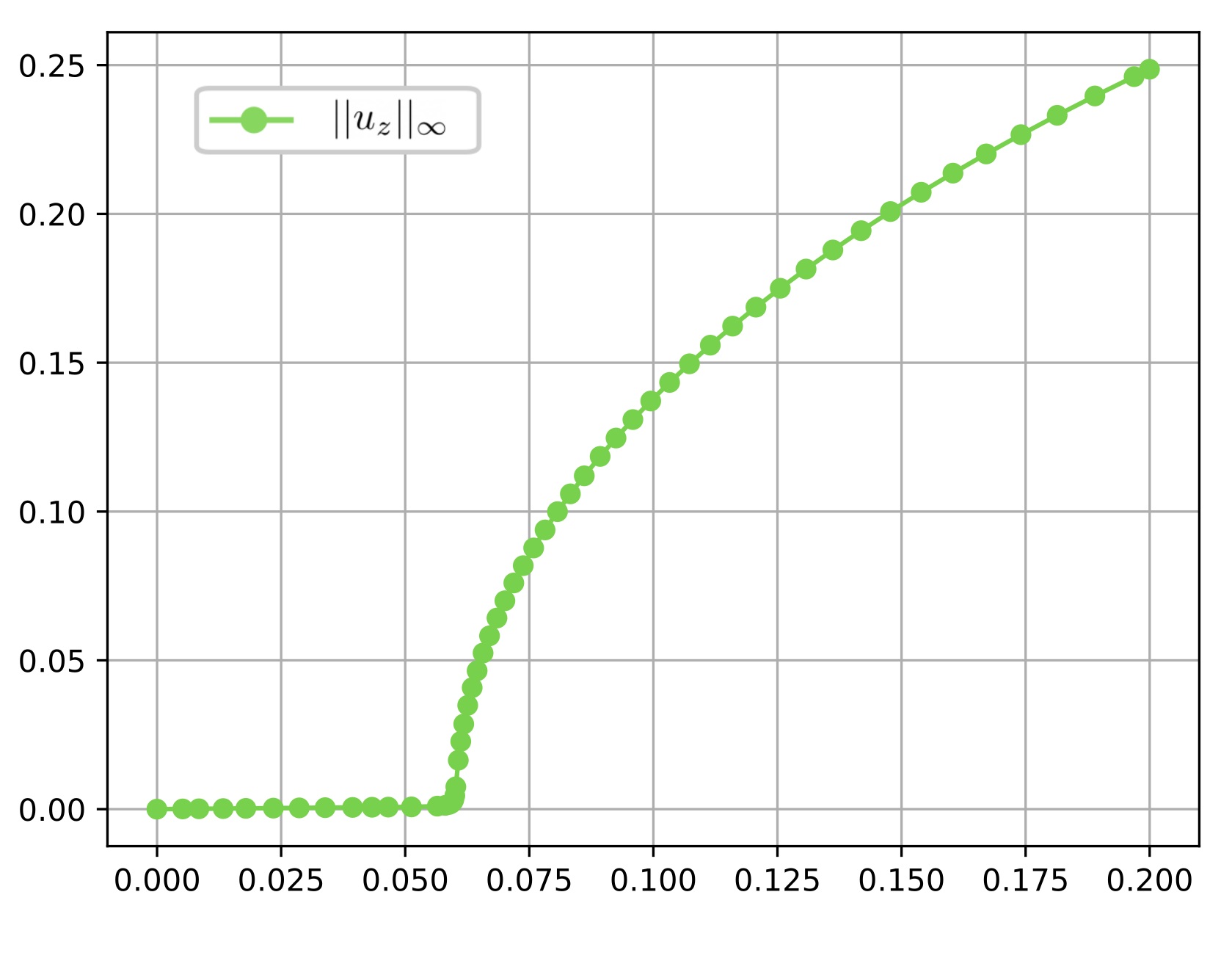}
                           \put(-95,4){\makebox(0,0){$\mu$}}
         \caption{$B = (0, 0, 0)$}
         \label{fig:15_a_3d_nh}
     \end{subfigure}
     \hfill
     \begin{subfigure}[b]{0.49\textwidth}
         \centering
         \includegraphics[width=\textwidth]{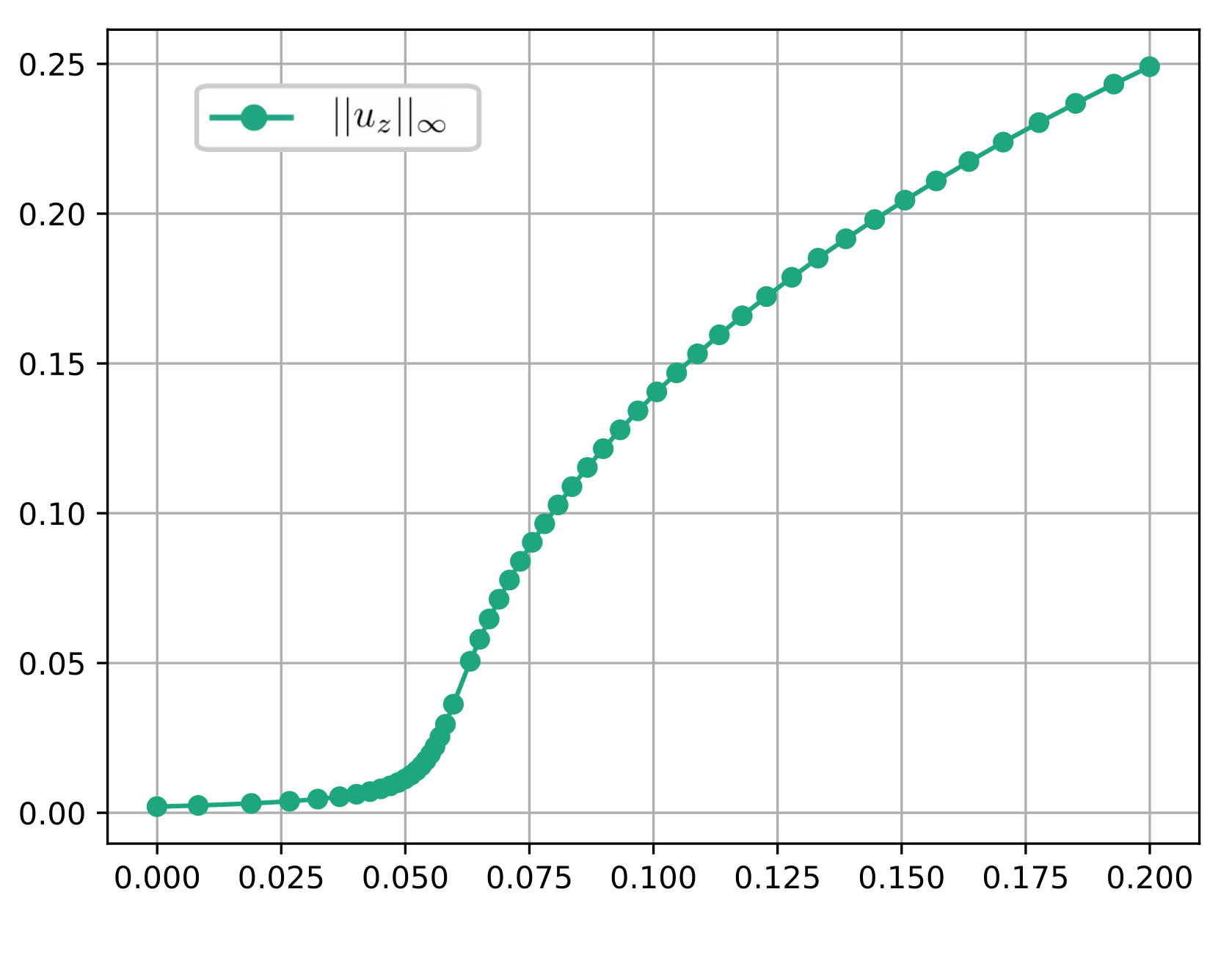}
                           \put(-95,4){\makebox(0,0){$\mu$}}
         \caption{$B = (0, 0, -1000)$}
         \label{fig:15_b_3d_nh}
     \end{subfigure}
	\caption{Reduced basis bifurcation diagrams for the 3D NH beam with respect to different body forces.}
\label{fig_15_3d_nh}
\end{figure}

For the sake of comparison let us consider also the NH model within the same setting as before.
The same POD tolerance was reached by a fewer number of basis function, namely $N = 5$.
In Figure \ref{fig_15_3d_nh} we can see the bifurcation plot for trivial $B = (0, 0, 0)$ and gravitational $B = (0, 0, -1000)$ body forces. For the three-dimensional beam we can observe that the NH constitutive relation actually predict the buckling in a slightly different location, indeed from Figure \ref{fig:15_a_3d_nh} we can detect the buckling occurring for the value $\mu^* \approx 0.059$ (compare with Figure \ref{fig:15_a_3d}).

\begin{figure}[h!]
\centering
     \begin{subfigure}[b]{0.49\textwidth}
         \centering
         \includegraphics[width=\textwidth]{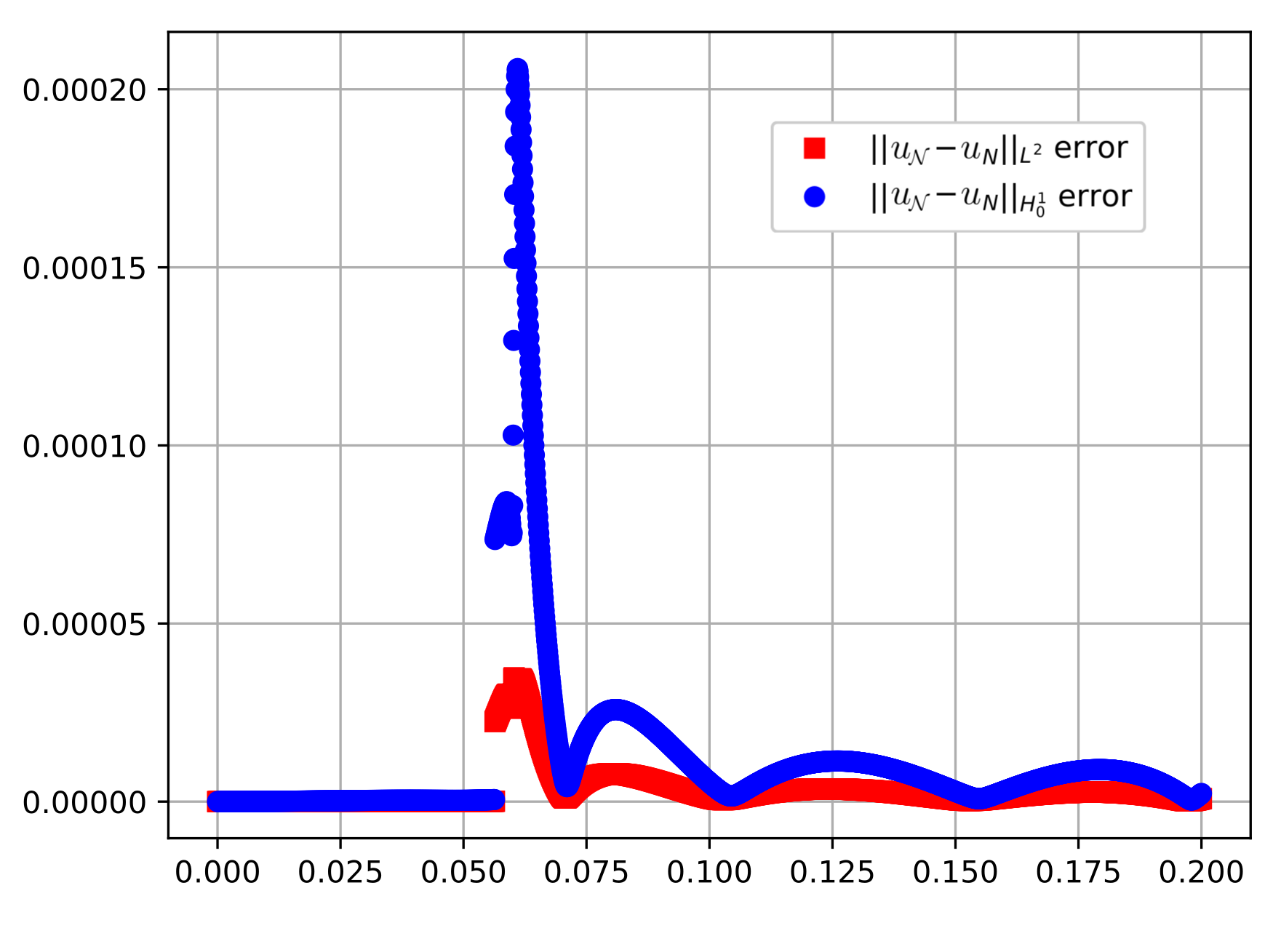}
                           \put(-92,4){\makebox(0,0){$\mu$}}
         \caption{$B = (0, 0, 0)$}
         \label{fig:10_a_3d_nh}
     \end{subfigure}
     \hfill
     \begin{subfigure}[b]{0.49\textwidth}
         \centering
         \includegraphics[width=\textwidth]{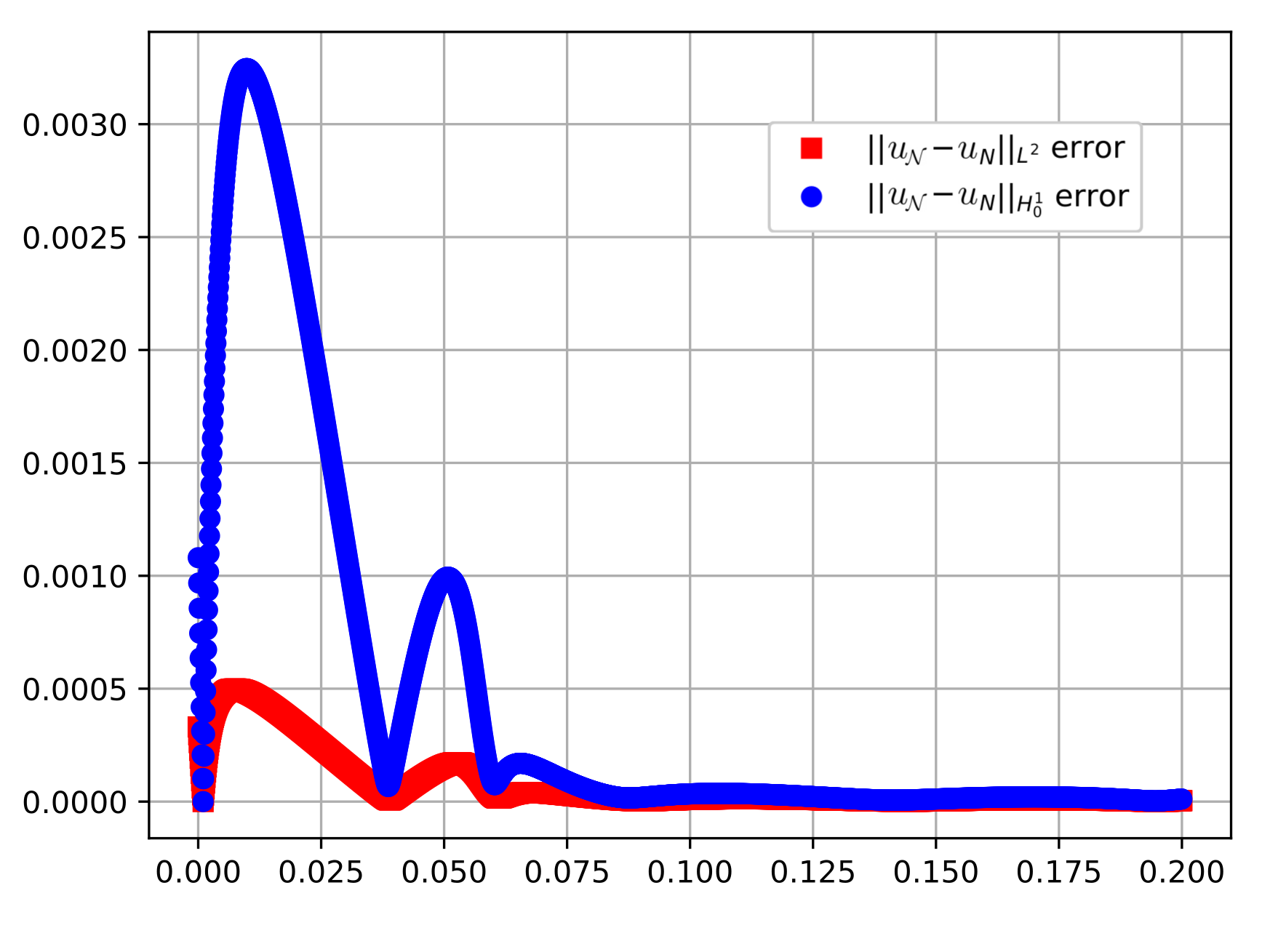}
                           \put(-92,4){\makebox(0,0){$\mu$}}
         \caption{$B = (0, 0, -1000)$}
         \label{fig:10_b_3d_nh}
     \end{subfigure}
\caption{Reduced basis errors for the 3D NH beam with respect to different body forces.}
\label{fig_10_3d_nh}
\end{figure}

\begin{figure}[h!]
\centering
\includegraphics[width=0.7\textwidth]{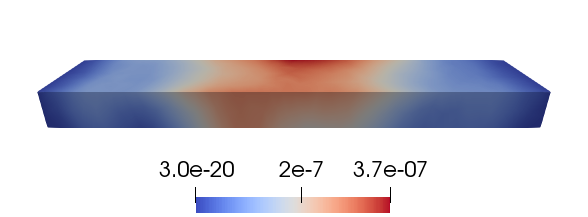}
\caption{Reduced basis error plot of the displacement $u$ for the 3D NH beam with $B = (0, 0, -1000)$ at $\mu = 0.2$.}
\label{fig_14}
\end{figure}

Furthermore, as we can see from Figure \ref{fig_10_3d_nh}, the NH model did not encounter the same accuracy issues as for the SVK model in Figure \ref{fig_10_3d}, thus showing a better approximation over $\Pa$. As an example we show in Figure \ref{fig_14} the reduced error for the post-buckling displacement solution at $\mu = 0.2$ with gravitational body force.

\begin{figure}[h]
\centering
\includegraphics[width=0.5\textwidth]{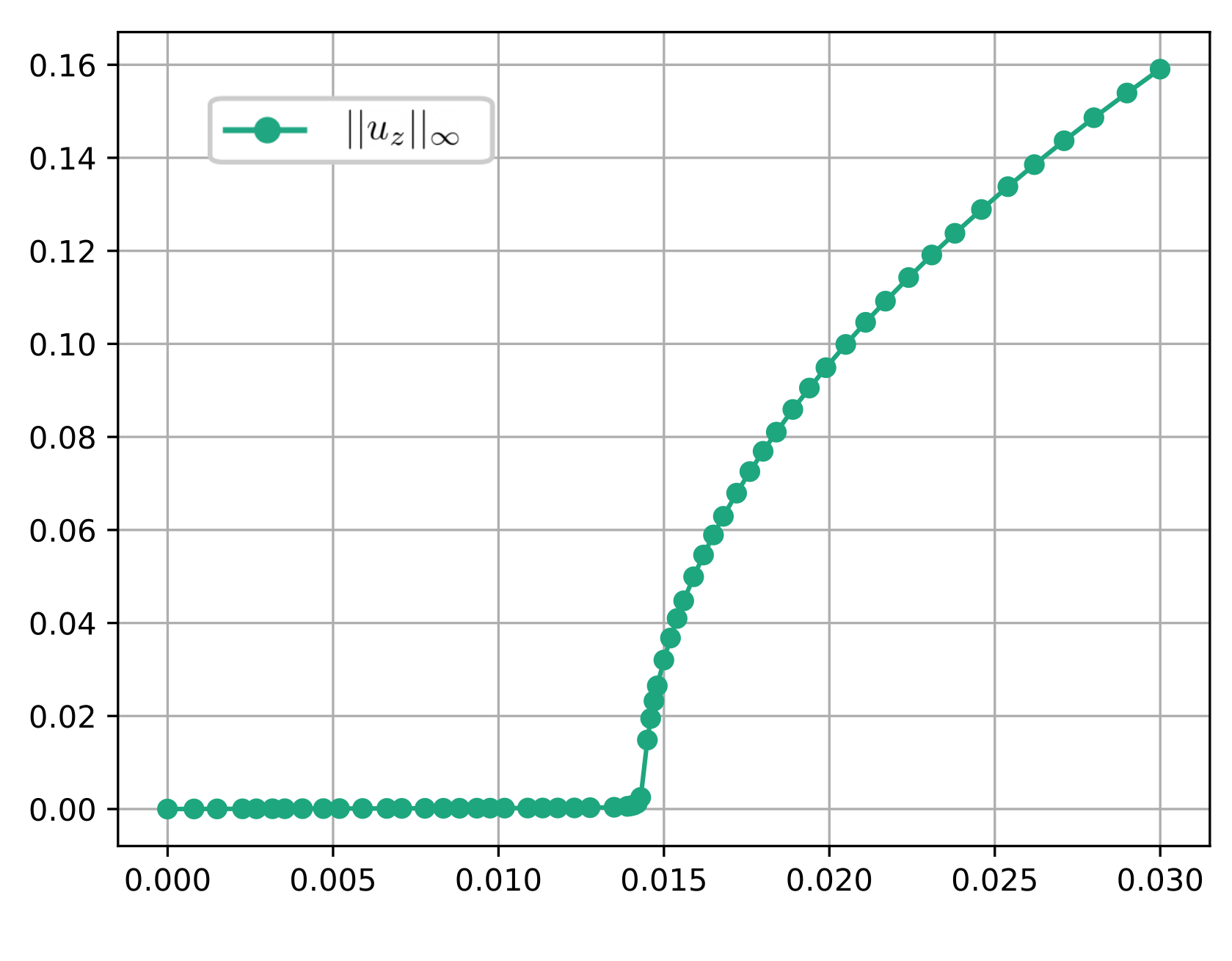}
\put(-100,4){\makebox(0,0){$\mu$}}
\caption{Reduced basis bifurcation diagram for the 3D SVK beam with $B = (0, 0, 0)$.}
\label{fig_20_BC}
\end{figure}

\begin{figure}[h]
\centering
\includegraphics[width=0.7\textwidth]{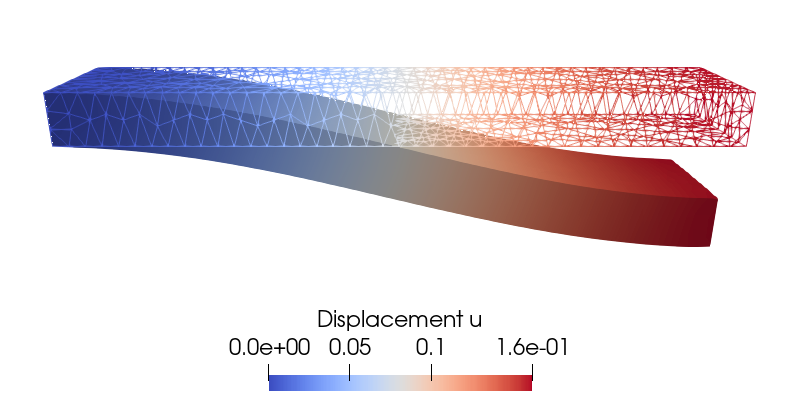}
\caption{High fidelity displacement $u$ for the 3D SVK beam with $B = (0, 0, 0)$ at $\mu = 0.03$.}
\label{fig_6_3d_bc}
\end{figure}

\begin{figure}[h]
\centering
\includegraphics[width=0.5\textwidth]{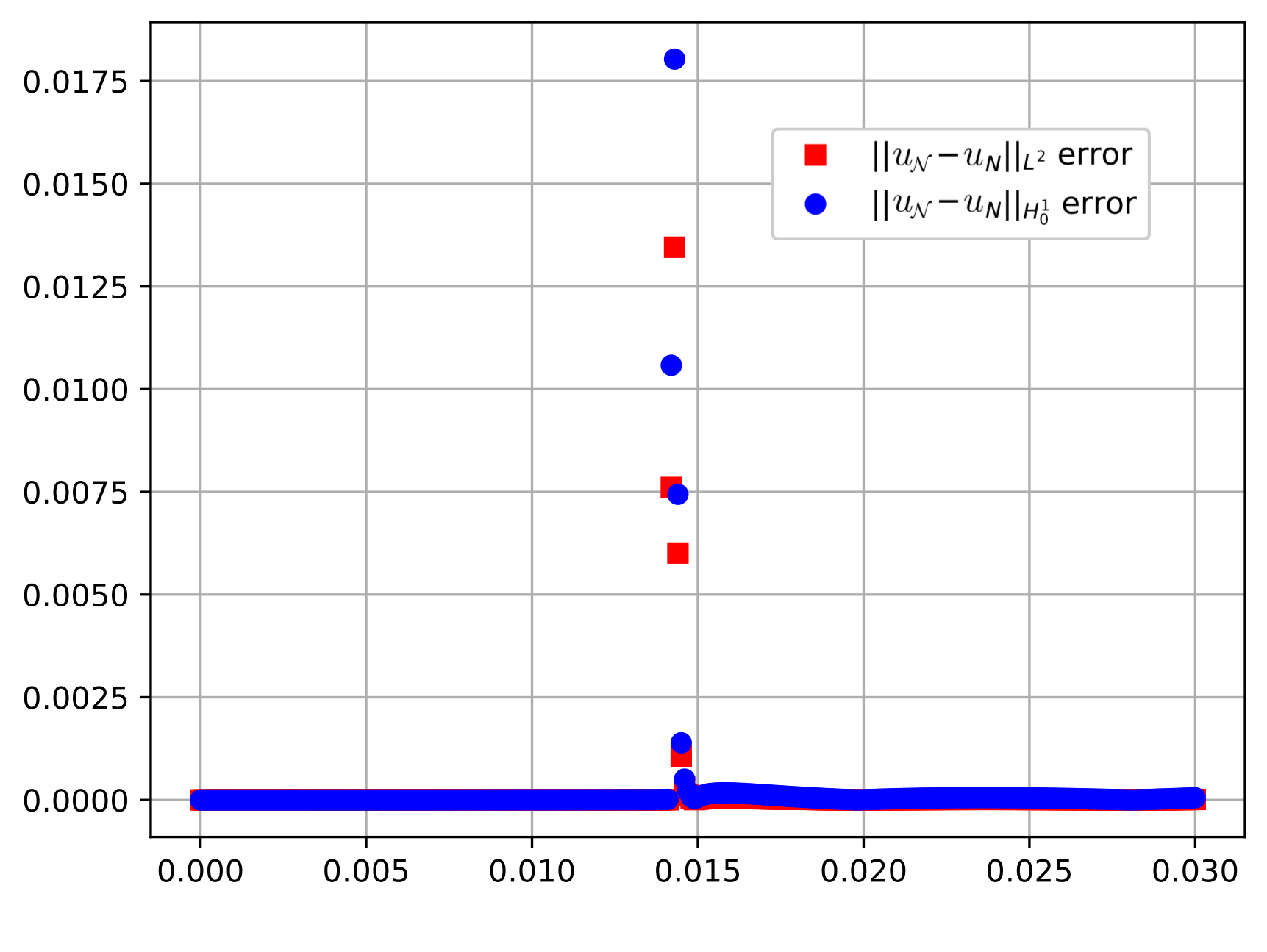}
\put(-100,4){\makebox(0,0){$\mu$}}
\caption{Reduced basis errors for the 3D SVK beam with $B = (0, 0, 0)$.}
\label{fig_6_3d_bc_err}
\end{figure}

Before ending this section we want to remark that as we have seen previously, also in this case the way in which we impose the compression and the boundary conditions chosen had a great influence on the buckling location.
In fact, if we allow the right end of the beam to move in the two perpendicular direction to the compression ($y$ and $z$ axis), the resulting buckling mode and its corresponding branching point change consistently.
In practice, we fixed the functional space as $$\X = \lbrace{ u \in (H^1(\Omega))^3: u = (0, 0, 0) \ \text{on}\  \Gamma^{l}_D,\ u_x = -\mu \ \text{on} \  \Gamma^{r}_D \rbrace} ,$$ and the parameter space as $\Pa = [0, 0.03]$.
Choosing the SVK constitutive relation with the trivial body force we obtained the reduced bifurcation diagram in Figure \ref{fig_20_BC}, where we can see the buckling occurring at $\mu^* = 0.014$.
The new buckled state solution is depicted in Figure \ref{fig_6_3d_bc} for $\mu = 0.2$. As we can see from Figure \ref{fig_6_3d_bc_err}, the reduced approximation follows the same behavior of the previous case, with $N = 4$ basis functions, maximum error of order $10^{-2}$ and a difference of 4 orders of magnitude with the average one.

\subsection{Industrial test case application to tubular element}\label{sec:norsok}

Finally, we can introduce a real test case scenario which comes from the Norwegian petroleum industry \cite{norsok}.
The results obtained in the previous sections allowed us to investigate the following more complex problem.
We want to investigate the deformation of a 3-D tubular geometry with annular section. 
It is clear that when dealing with real geometries the situation becomes more involved, since many factors have to be taken into consideration.
Here, for practical interest the domain is defined as a tubular member characterized by an annular cross-section with inner and outer radii $r = 0.28$(m) and $R=0.30$(m).
Thus, the thickness corresponds to $t = 0.02$(m) and the outer diameter is $D = 0.6$(m). This is important since the design condition on the ratio $D/t < 120$ required in \cite{norsok} has to be satisfied for the correct reconstruction of the physics at hand.
As before, different settings will be analyzed in order to understand their buckling properties. Despite the great difficulties and some unexpected behavior that such model exhibits, in the following we present a successful application of our reduced methodology.

\subsubsection{A comparison of constitutive relations}

In this section we fix the length of the tube as $L = 2$(m) and our focus will be on the behavior of the displacement
in the case of Neumann compression.
Hence, we fix the domain as $\Omega = A_r^R \times [0, 2]$ where we denoted the annular cavity as $$A_r^R = \lbrace{(x,y) \in \R^2 | r^2 \leq x^2 + y^2 \leq R^2\rbrace} ,$$ we choose a trivial body force $B = (0, 0, 0)$ and a compression given by the traction term $T = (0,0,-\mu)$. We remark now that the compression is acting along the $z$-axis. For what concerns the material properties we fixed $E = 2.1\cdot 10^5$(MPa) and $\nu = 0.3$. 

\begin{figure}[h]
\centering
\includegraphics[width=0.3\textwidth]{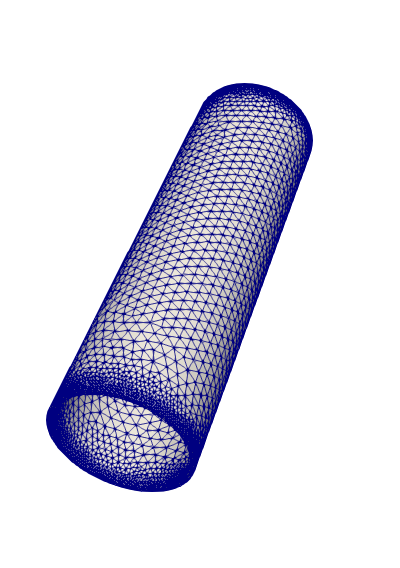}
\caption{Mesh for the Norsok test case.}
\label{fig:norsok_geom}
\end{figure}

We start with the one parameter test case, where homogeneous Dirichlet boundary condition is imposed at $\Gamma_D = A_r^R \times \lbrace{0\rbrace}$ and $\Gamma_N = A_r^R \times \lbrace{2\rbrace}$.
We built the FE space by means of $\mathbb{P}_1$ linear elements on the tetrahedral mesh in Figure \ref{fig:norsok_geom} with $22521$ cells, while the resulting number of degrees of freedom is $\N = 23442$.

Given the complexity of the geometry and since our main interest is to find the first failure mode, here we will focus the investigation only up to the buckling, disregarding the post-buckling behavior. Both constitutive relations, SVK and NH, have been analyzed, producing qualitatively different results also concerning the position of the critical value. 
For these reasons we considered different parameter spaces for the two hyperelastic models.

For the SVK model the parameter space is given by $\Pa = [0, 4500]$ and $N_{train} = 200$ snapshots were computed through the continuation method.
The reduced manifold was built choosing a tolerance $\epsilon_{POD} = 10^{-8}$ and this led to a reduced basis space of dimension $N = 5$.
The bifurcation diagram was recovered with an online continuation method based on $K = 451$ equispaced points in $\Pa$.

\begin{figure}[h!]
\centering
     \begin{subfigure}[b]{0.49\textwidth}
         \centering
         \includegraphics[width=\textwidth]{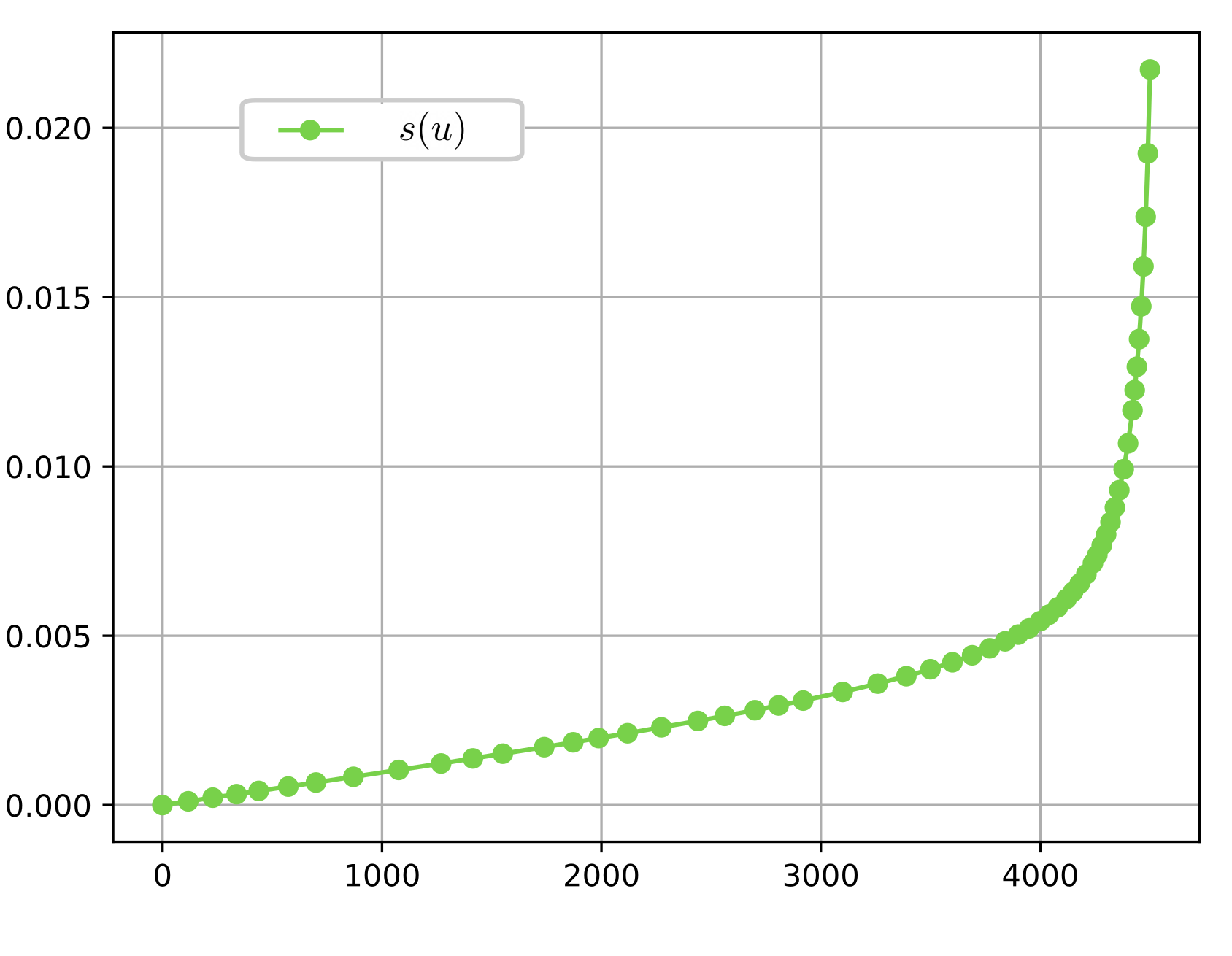}
                  \put(-95,4){\makebox(0,0){$\mu$}}
         \caption{SVK}
         \label{fig:15_a_3d_nor}
     \end{subfigure}
     \hfill
     \begin{subfigure}[b]{0.49\textwidth}
         \centering
         \includegraphics[width=\textwidth]{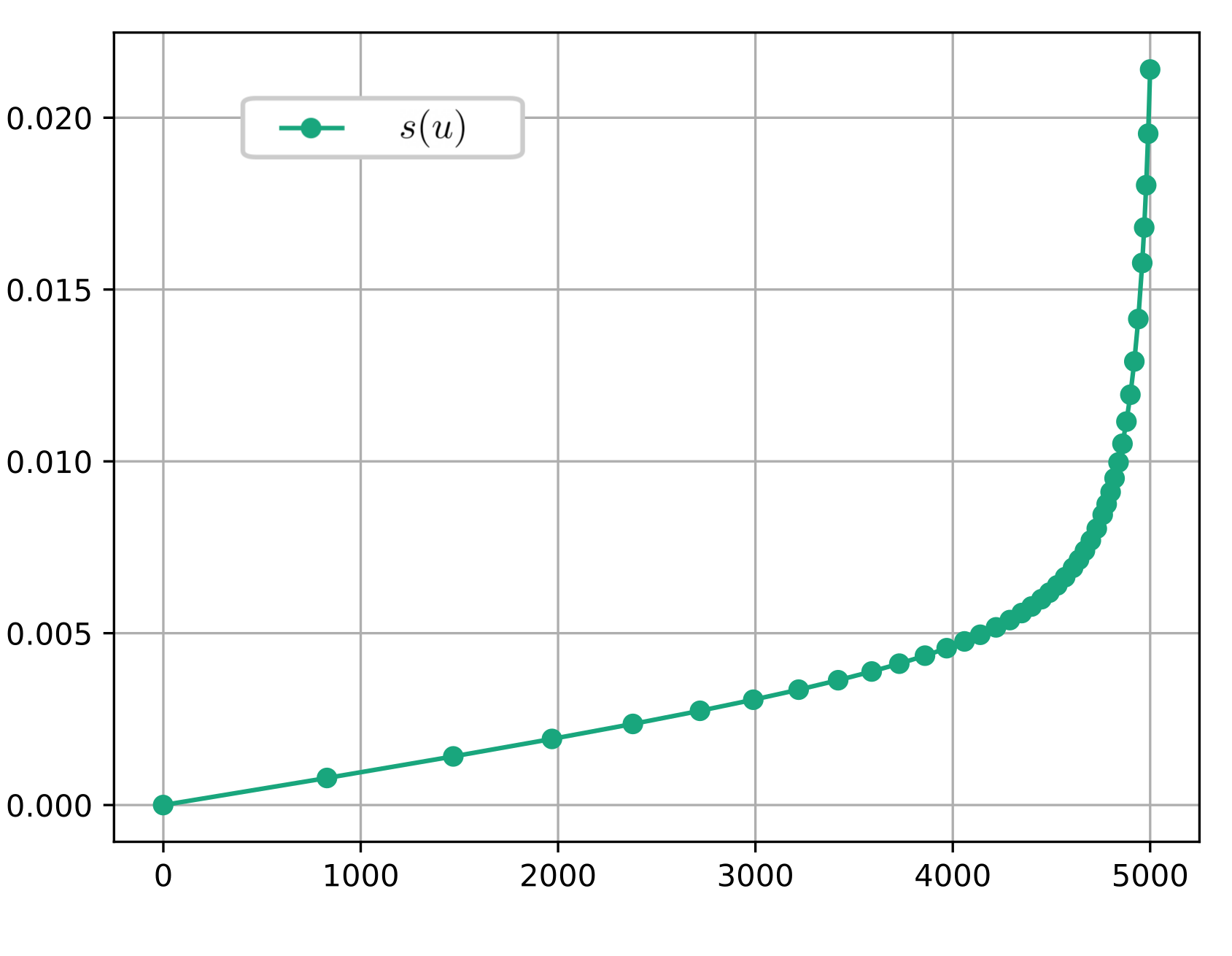}
                  \put(-95,4){\makebox(0,0){$\mu$}}
         \caption{NH}
         \label{fig:15_b_3d_nor}
     \end{subfigure}
	\caption{Reduced basis bifurcation diagrams for the 3D tubular geometry with $B = (0, 0, 0)$ and different constitutive relations.}
\label{fig_15_3d_nor}
\end{figure}

For what concerns the NH material, we computed again $N_{train} = 200$ snapshots, but the parameter space now is defined as $\Pa = [0, 5000]$. The same reduced setting was applied, providing a reduced basis space of the same dimension.
Similarly to the previous 3-D case, we remark that given the symmetry w.r.t.\ the tube axis, here the buckling can occur in any direction perpendicular to the axis itself. To detect the buckling behaviour, we consider as output functional the sum of the infinite norms of the $x$ and $y$ components of the displacement, namely $s(u) =  \norm{u_x}_{\infty} + \norm{u_y}_{\infty}$. 
In Figures \ref{fig:15_a_3d_nor} and \ref{fig:15_b_3d_nor} we can see the bifurcation plot for the SVK and NH constitutive relations, respectively. We can clearly observe that the buckling of the beam with different models occurs at different values for the compression $\mu$, in particular the buckling of the SVK beam is slightly anticipated.

\begin{figure}[h]
\centering
\includegraphics[width=0.3\textwidth]{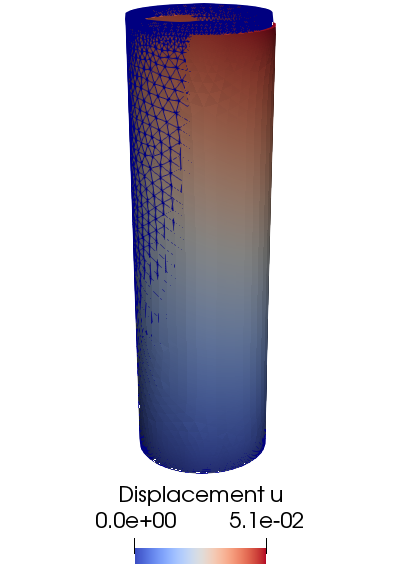}\hspace{1cm}
\includegraphics[width=0.3\textwidth]{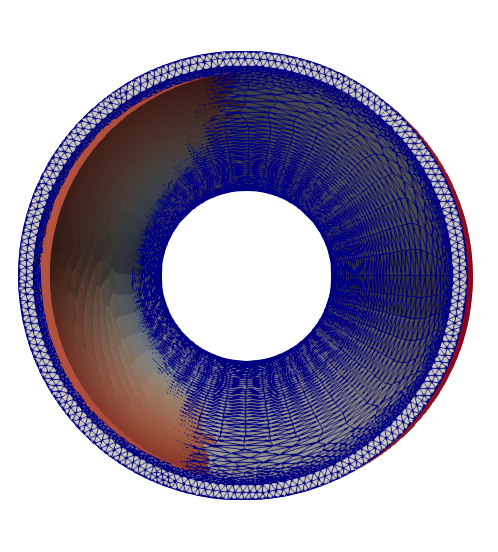}
\caption{High fidelity displacement $u$ for the 3D SVK tubular geometry with $B = (0, 0, 0)$ at $\mu = 4500$.}
\label{fig_6_3d_nor}
\end{figure}

As expected, the buckling is qualitatively similar to the one observed in Figure \ref{fig_err_sol_svk_force}, with Neumann compression in the 2-D geometry. A representative solution of the buckling mode for the SVK model is depicted in Figure \ref{fig_6_3d_nor} for $\mu = 4500$. 

\begin{figure}[h!]
\centering
     \begin{subfigure}[b]{0.49\textwidth}
         \centering
         \includegraphics[width=\textwidth]{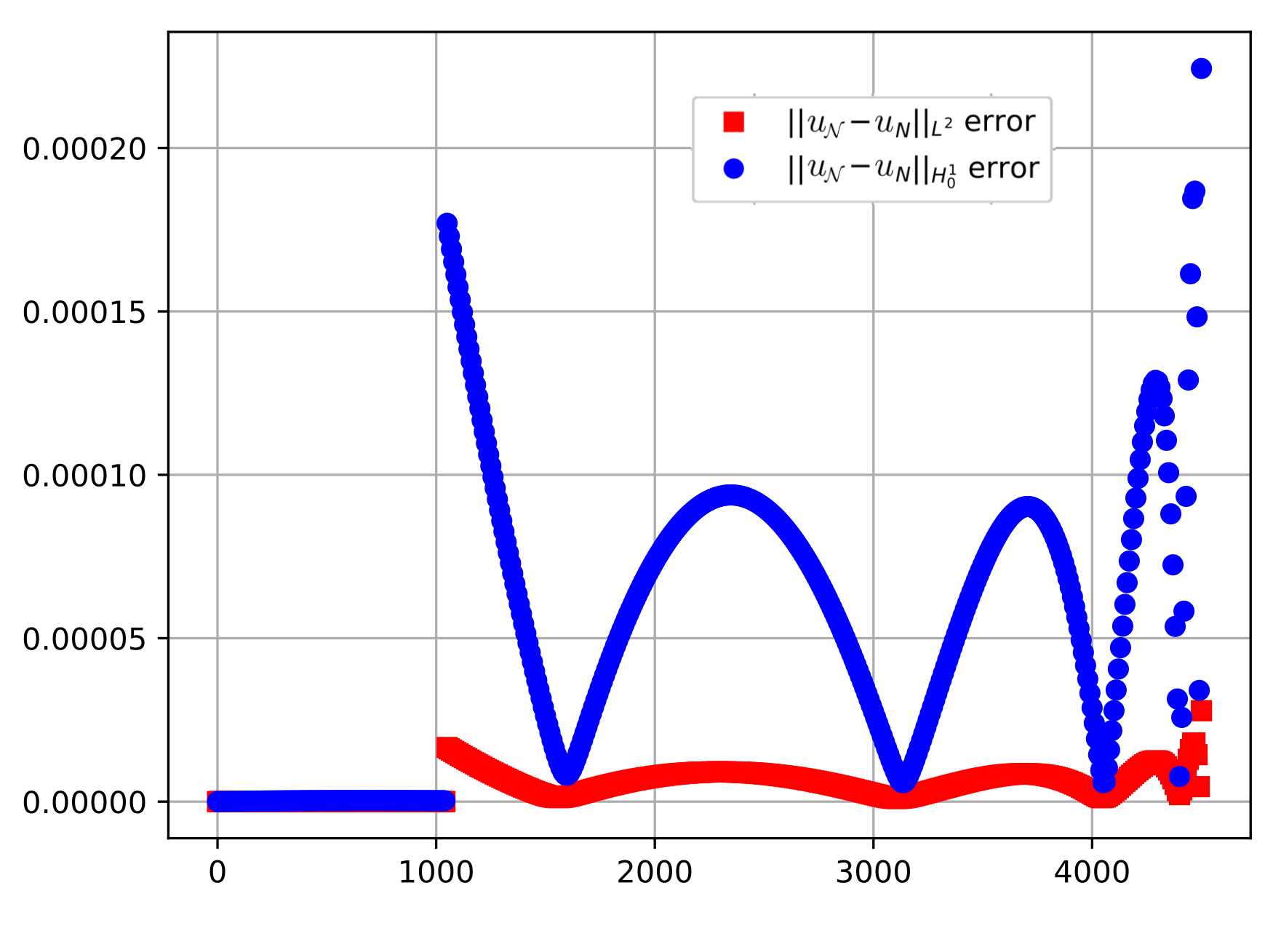}
                           \put(-95,4){\makebox(0,0){$\mu$}}
         \caption{SVK}
         \label{fig:10_a_3d_nor}
     \end{subfigure}
     \hfill
     \begin{subfigure}[b]{0.49\textwidth}
         \centering
         \includegraphics[width=\textwidth]{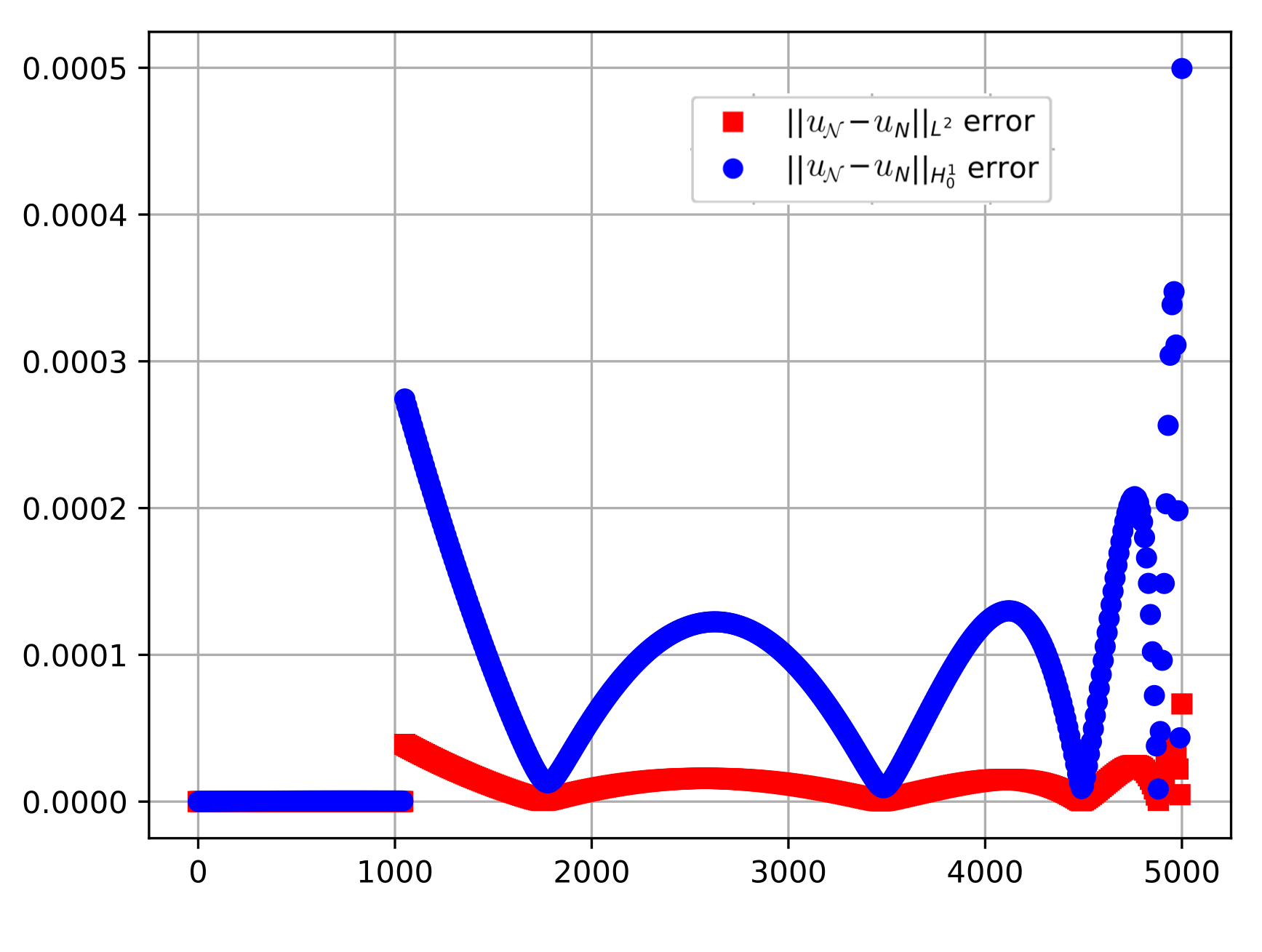}
                           \put(-95,4){\makebox(0,0){$\mu$}}
         \caption{NH}
         \label{fig:10_b_3d_nor}
     \end{subfigure}
\caption{Reduced basis errors for the 3D tubular geometry with $B = (0, 0, 0)$ and different constitutive relations.}
\label{fig_10_3d_nor}
\end{figure}

As we can see from Figure \ref{fig_10_3d_nor}, in both cases we were able to reach a good accuracy of the RB solution with respect to $\mu$. Once again the error increases near the buckling point, here, the right end of the parameter domains, where the solution changes more rapidly its behavior.
Regarding the computational time, we obtained a low speed-up of order $1.25$, in fact to plot the high fidelity versions of the bifurcation diagrams in Figure \ref{fig:15_b_3d_nor} we spent $t_{HF} = 1335$(s) while the reduced order one required $t_{RB} = 1076$(s). Of course this is not satisfactory, 
for this reason, we now present an application of the hyper-reduction strategies to recover the efficiency.

In fact, especially when dealing with real test cases, a real time evaluation of the solution is a key feature.  
Thus, we consider the geometry $\Omega$ with Neumann compression and we try to efficiently recover the bifurcation diagrams in Figure \ref{fig_15_3d_nor} by means of the DEIM. 

Within the same setting as before, we consider the DEIM with Greedy tolerance $\epsilon_{Gr} = 10^{-10}$ which splits the form into an affine decomposition made up by 2 terms, each one approximated by 15 interpolation basis functions. The reason for such low tolerance is that the complexity of the model, together with its buckling behavior, makes the approximation of the variational forms a difficult task. Indeed, we observed a non-convergence issue during the online phase when higher tolerances were chosen.
   
\begin{figure}[h!]
\centering
     \begin{subfigure}[b]{0.49\textwidth}
         \centering
         \includegraphics[width=\textwidth]{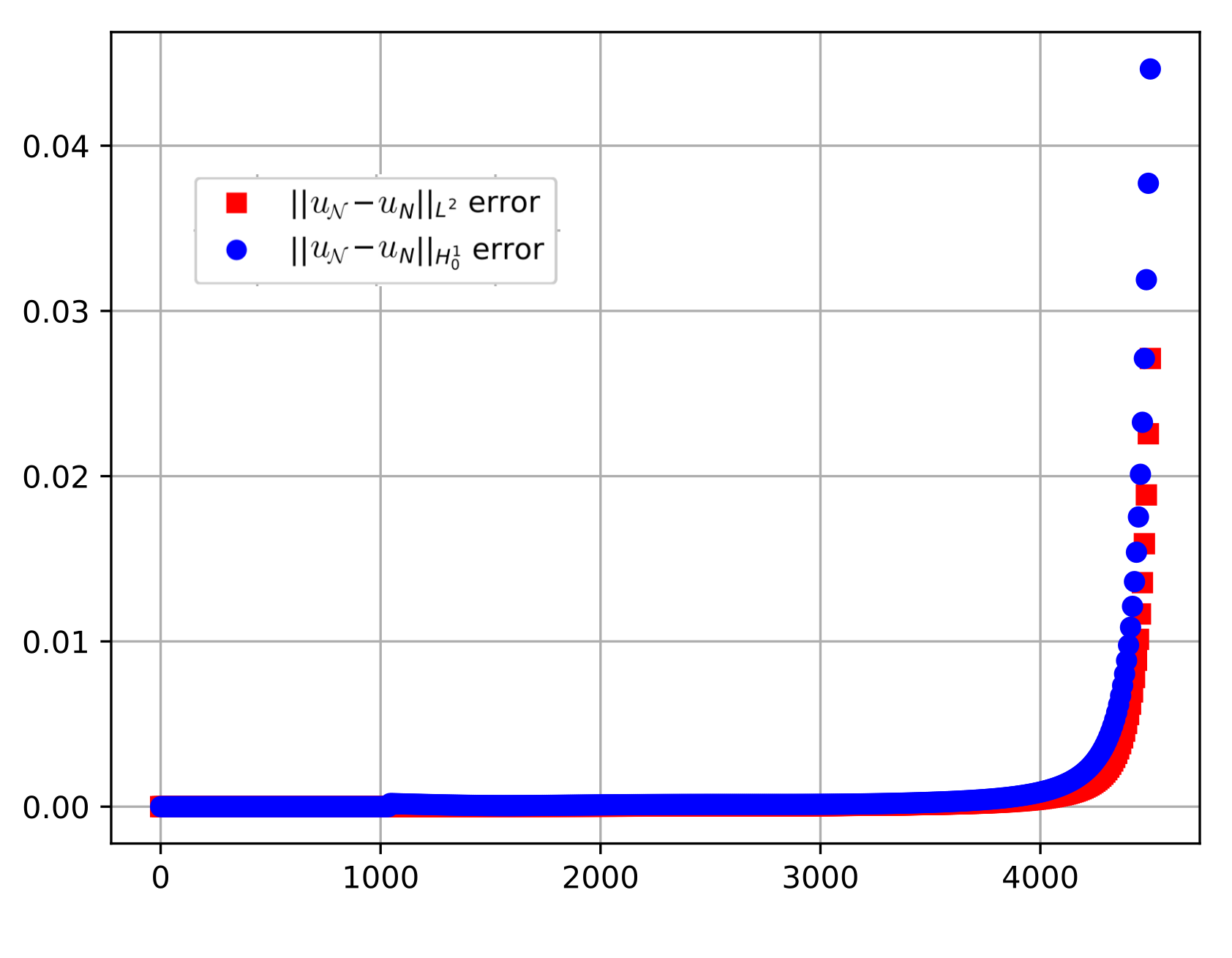}
                           \put(-95,4){\makebox(0,0){$\mu$}}
         \caption{SVK}
         \label{fig:10_a_3d_nor_deim}
     \end{subfigure}
     \hfill
     \begin{subfigure}[b]{0.49\textwidth}
         \centering
         \includegraphics[width=\textwidth]{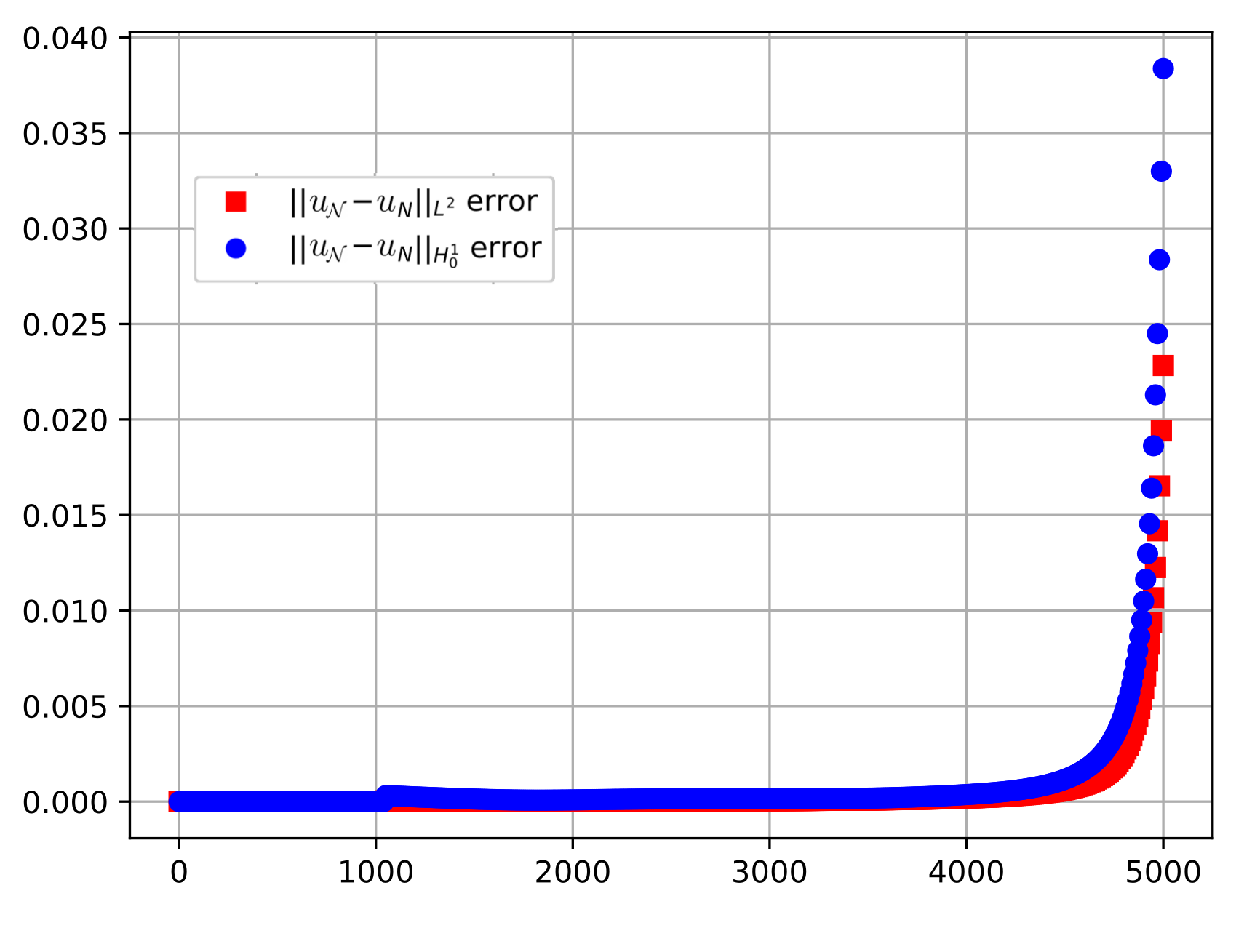}
                           \put(-95,4){\makebox(0,0){$\mu$}}
         \caption{NH}
         \label{fig:10_b_3d_nor_deim}
     \end{subfigure}
\caption{Reduced basis errors with DEIM for the 3D tubular geometry with $B = (0, 0, 0)$ and different constitutive relations.}
\label{fig_10_3d_nor_deim}
\end{figure}

Hence, let us show in Figure \ref{fig_10_3d_nor_deim} the reduced basis error, computed with DEIM for the SVK and NH constitutive relations. It is evident that a significant increment of both maximum errors, w.r.t.\ the case without interpolation strategies, occurs. Thus, the hyper-reduction approach seems to have some difficulties in approximating the forms near critical points. Despite this, the good news is that while the reduced basis error is still acceptable, the speed-up consistently increases. As an example, for NH material we pass from a speed-up of order 1.25 with $t_{RB} = 1076$(s) to a speed-up of order 38 corresponding to $t_{RB,DEIM} = 35$(s).

As we understood, usually one is not interested in the approximation of the bifurcating phenomenon for a fixed setting, rather the aim is the detection and reconstruction of the buckling modes varying physical or geometrical parameter. Since this is true especially for real test cases, the latter will be analyze in the next section.

\subsubsection{Multi-parametric study with varying geometries}\label{geom_nor_force}

Now we want to extend the previous analysis to the multi-parameter context, in which the additional parameter controls the length of the domain $\Omega$. To do so, we recall the investigation done in Section \ref{ssec:geom_param}.
For this reason, we now consider the geometrically parametrized tube, where again its semi-length $\mu_g \in \Pa_g = [1, 2]$ is taken into consideration.
Therefore, we can express the 3-D parametrized domain as $\widetilde{\Omega}(\mu_g) = \widetilde{\Omega}_1 \cup \widetilde{\Omega}_2(\mu_g)$, where $\widetilde{\Omega}_1 = A_r^R \times [0, 1]$ and $\widetilde{\Omega}_2(\mu_g) = A_r^R \times [1, \mu_g]$.
Therefore, the transformation map is simply given by the affine function $$\Phi(x; \mu_g) =  \begin{bmatrix} x \\  y \\ \mu_g(z - 1) + 1 \end{bmatrix} \quad \text{for } x \in \Omega_2 =  \widetilde{\Omega}_2(1) .$$

\begin{figure}[h]
\centering
\includegraphics[width=0.7\textwidth]{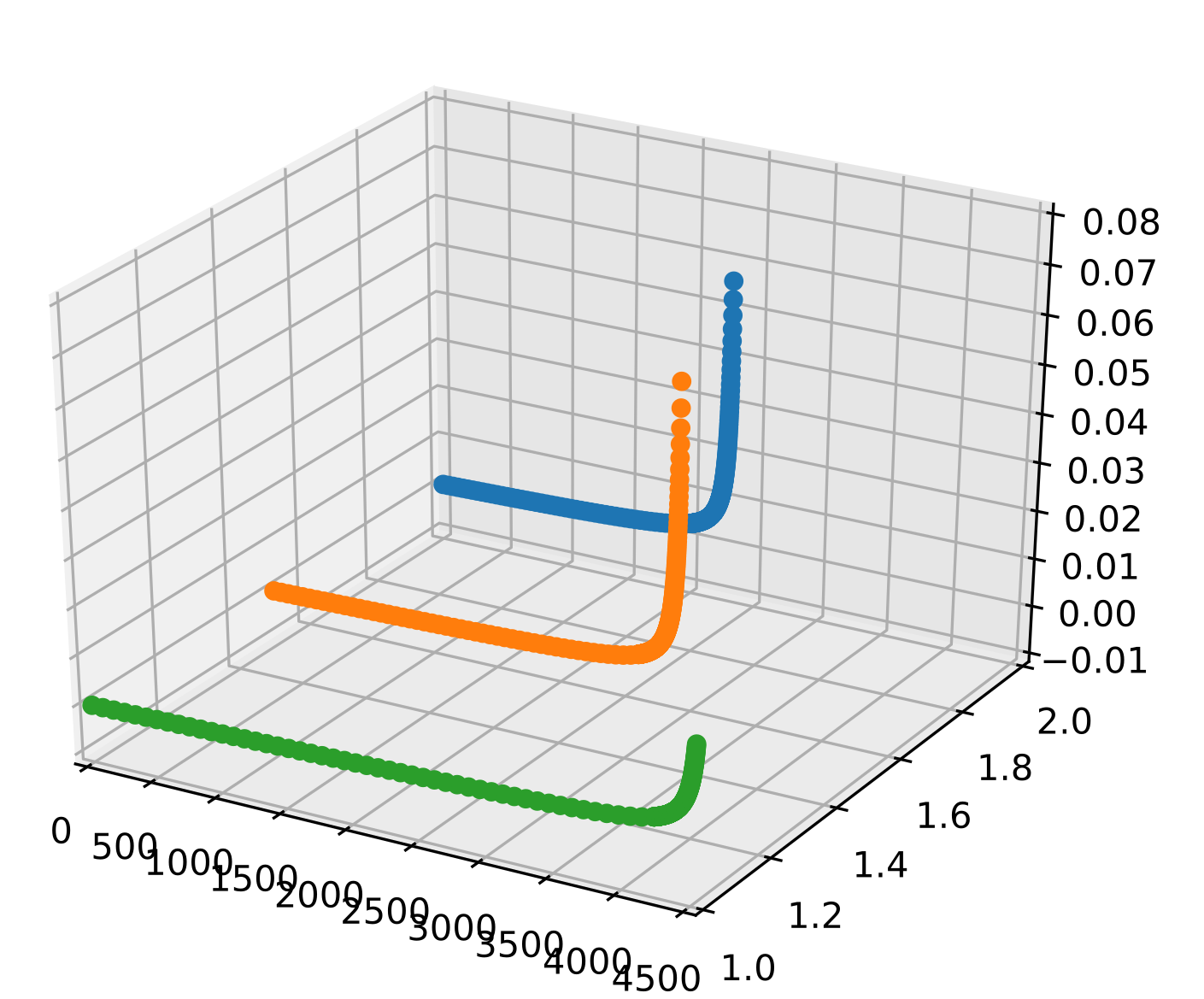}
    \put(-5,115){$s(u)$}
    \put(-50,28){$\mu_g$}
    \put(-220,12){$\mu$}
\caption{3D bifurcation plot for SVK beam with $B = (0, 0, 0)$ and $\mu_g \in \Pa_g$.}
\label{fig_17_force_param}
\end{figure}

\begin{figure}[h!]
\centering
    \begin{subfigure}[b]{0.32\textwidth}
         \centering
         \includegraphics[width=\textwidth]{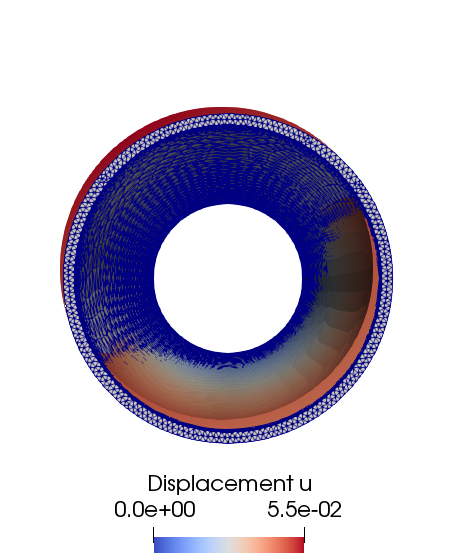}
         \caption{$\mu_g = 1.$}
         \label{fig:6_a_geom_1_force_nor}
     \end{subfigure}
     \hfill
    \begin{subfigure}[b]{0.32\textwidth}
         \centering
         \includegraphics[width=\textwidth]{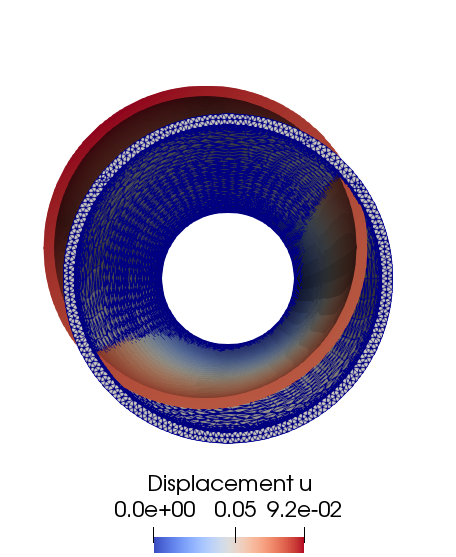}
         \caption{$\mu_g = 1.5$}
         \label{fig:6_b_geom_1_force_nor}
     \end{subfigure}
          \hfill
    \begin{subfigure}[b]{0.32\textwidth}
         \centering
         \includegraphics[width=\textwidth]{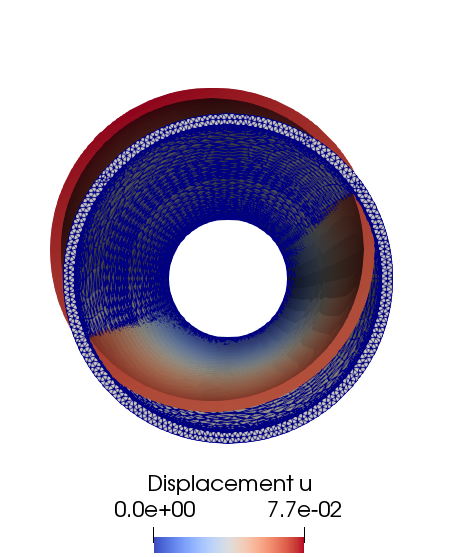}
         \caption{$\mu_g = 2.$}
         \label{fig:6_c_geom_1_force_nor}
     \end{subfigure}
\caption{Top-view of the high fidelity displacement $u$ for the SVK beam with $B = (0, 0, 0)$ for different geometries.}
\label{fig_6_geom_force_nor}
\end{figure}

The physical setting with Neumann BCs is the same as before, while for the offline phase we computed $N_{train} = 500$ snapshots, corresponding to the solutions computed for 3 random values of $\mu_g$, and values of $\mu$ up to the buckling parameter for each different geometry. Performing a global POD compression with tolerance $\epsilon_{POD} = 10^{-8}$ we obtained a reduced basis space of dimension $N = 9$. 
The online continuation method to reconstruct the 3-D bifurcation diagram in Figure \ref{fig_17_force_param} is based on $K = 350$ equispaced points in $\Pa_p = [0, 4500]$, for 3 equispaced values of the semi-length $\mu_g \in \Pa_g$. We remark again that due to the large strains and the complexity of the phenomena, we focus on the solution behaviors only up to their buckling point. For these reasons, the actual parameter space varies for different branches, in fact we truncated it after having detected the buckling.
The top-view of three representative solutions of the buckling modes for $\mu_g = \left\{1, 1.5, 2\right\}$ are depicted in Figure \ref{fig_6_geom_force_nor} for the last computed value of $\mu =  \left\{4530, 3150, 2300\right\}$, respectively.

\begin{figure}[h!]
\centering
     \begin{subfigure}[b]{0.49\textwidth}
         \centering
         \includegraphics[width=\textwidth]{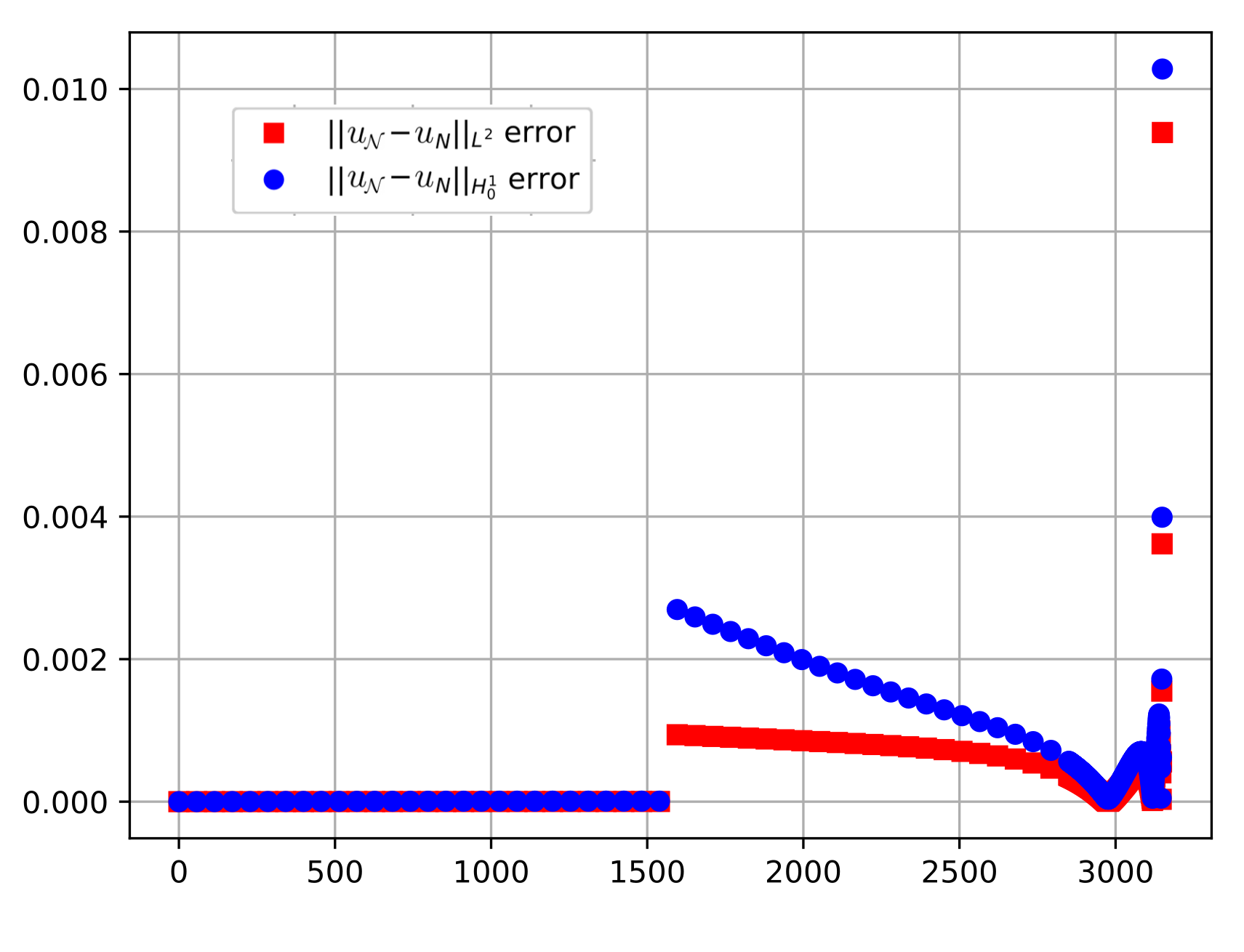}
                           \put(-95,4){\makebox(0,0){$\mu$}}
         \caption{$\mu_g = 1.5$}
         \label{fig:10_a_geom_force_nor}
     \end{subfigure}
     \hfill
     \begin{subfigure}[b]{0.49\textwidth}
         \centering
         \includegraphics[width=\textwidth]{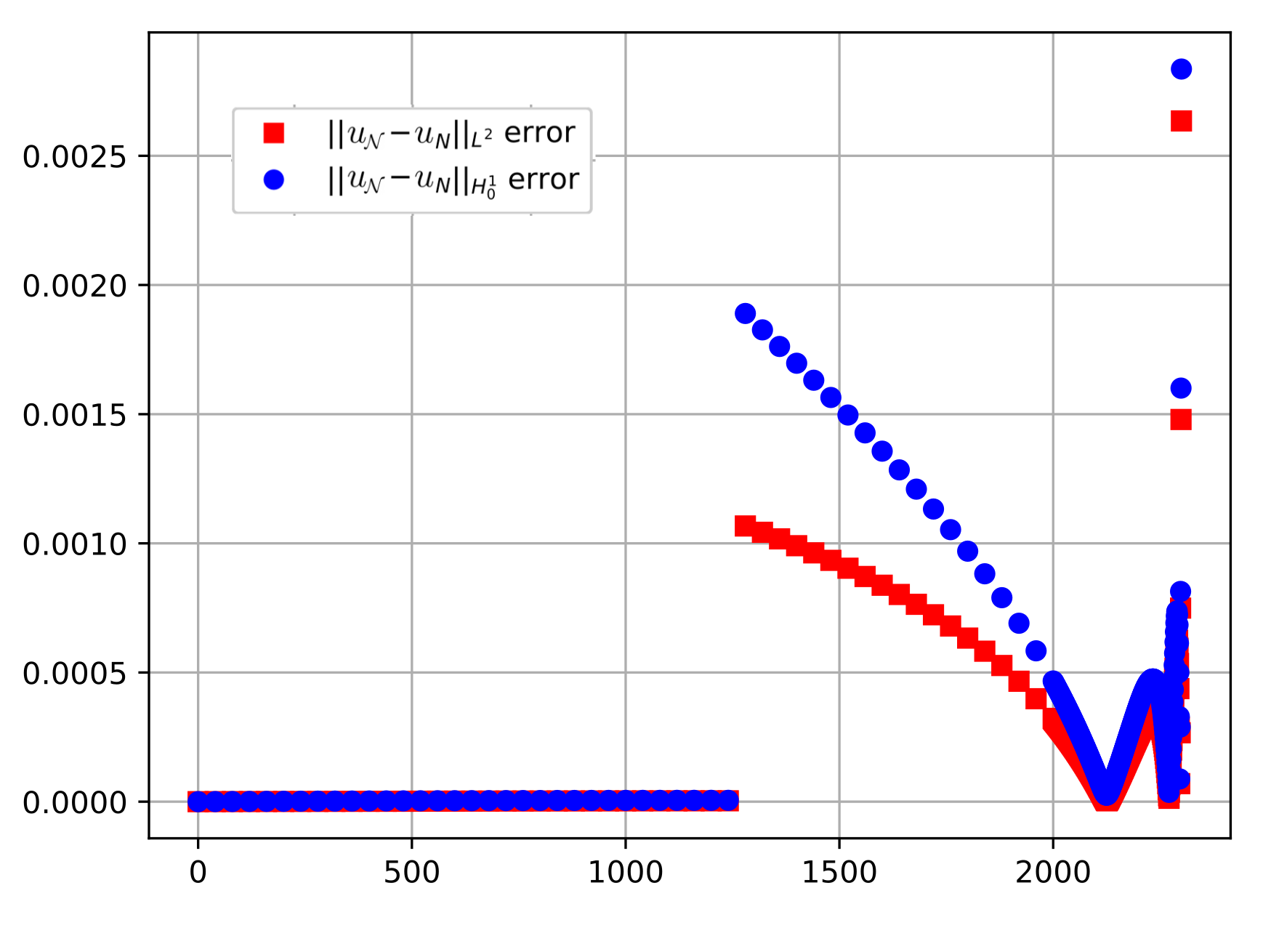}
                           \put(-95,4){\makebox(0,0){$\mu$}}
         \caption{$\mu_g = 2.0$}
         \label{fig:10_b_geom_force_nor}
     \end{subfigure}
	\caption{Reduced Basis errors with respect to $\mu \in \Pa$ for the SVK beam with $B = (0, 0, 0)$ for different geometries.}
\label{fig_10_geom_force_nor}
\end{figure}

We present in Figure \ref{fig_10_geom_force_nor} the reduced basis error.
We highlight that, as we can see from Figure \ref{fig:10_a_geom_force_nor}, also in this case the reduced manifold was able to approximate the buckling for unsampled geometries. 

The study of the evolution of the buckling, varying the length of the beam, confirm as expected that longer beams  need smaller forces for the buckling. We highlight that given the importance of the buckling location, already during the offline stage, one can not perform an extensive analysis on the parametrized geometries. Moreover, to prevent non-convergence issues due to bigger strains one should consider a much refined mesh, increasing the high fidelity dimension.

Once again, the speed-up equals only to 1.4 is not that much relevant, with $t_{HF} = 2363$(s) and  $t_{RB} = 1677$(s). Despite this, we remark that given the high dimensionality of the finite element space, especially for 3-D geometry with refined mesh and polynomial of order $P > 1$, the naive RB approach without empirical interpolation strategies can still provide consistent speed-up.

Up to now, we have highlighted the difficulties while modelling the buckling through the Neumann compression, in the next section the aim will be the investigation of the Dirichlet compression.

\subsubsection{Dirichlet compression}
Here, we want to extend the study of the Dirichlet compression for the tubular 3-D geometry in Figure \ref{fig:norsok_geom}.
So far, we have understood that the choice of compression which has to be modelled is a complex task. Another confirmation of this comes from the following scenario.
In fact, when we tried to apply a Dirichlet compression over $\Gamma_D = A_r^R \times \lbrace{2\rbrace}$, we encountered many difficulties.

\begin{figure}[h]
\centering
\includegraphics[width=0.3\textwidth]{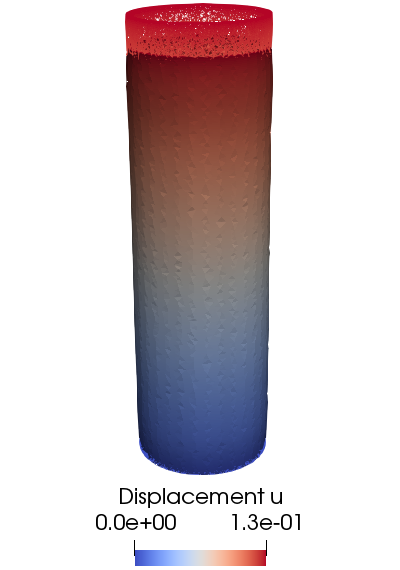} \hfill
\includegraphics[width=0.3\textwidth]{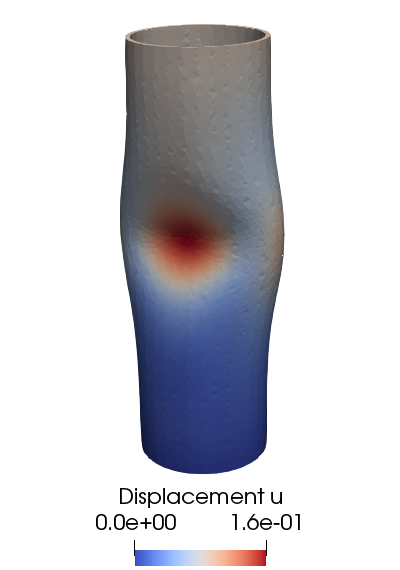} \hfill
\includegraphics[width=0.3\textwidth]{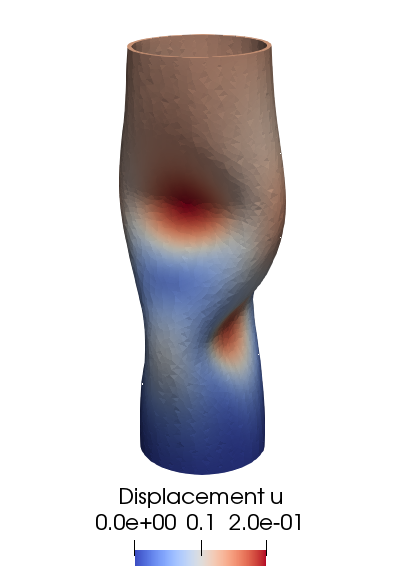}
\caption{Representative solutions of the 3D SVK model for the tubular geometry with Dirichlet compression.}
\label{fig_22}
\end{figure}

In particular, the cavity and the need for huge refinement of the mesh (connected to the remark done previously), made the reconstruction of the bifurcation diagram a too much difficult task. 
For these reasons, we were not able to fully recover the buckling within this compression context, but we only found the three different modes depicted in Figure \ref{fig_22}. 

From the analysis of the parametrized geometries we understood that the length of the tube plays a fundamental role in the buckling detection, thus we decided to investigate the tubular geometry with the same ratio, while considering a much increased length.

Therefore, let us consider the tube represented by the domain $\Omega = A_r^R \times [0, 20]$. We chose SVK material, trivial body and traction forces, and we fixed the material properties $E = 2.1 \cdot 10^5$ and $\nu = 0.3$. 
We first consider the one parameter test case, in which the parameter controls the compression through Dirichlet BC on  $\Gamma_D = A_r^R \times \lbrace{20\rbrace}$.
In this case the tetrahedral mesh consists of 147133 cells, resulting in a high fidelity dimension $\N = 147852$ when $\mathbb{P}_1$ linear elements are used.

Here, we fixed the parameter space as $\Pa = [0, 0.13]$, and its exploration is performed through the simple continuation method with step $\Delta\mu = 10^{-3}$, resulting in $N_{train} = 130$ snapshots.
The reduced manifold was built choosing a tolerance $\epsilon_{POD} = 10^{-8}$, which provides a reduced basis space of dimension $N = 5$.

\begin{figure}[h]
\centering
\includegraphics[width=0.5\textwidth]{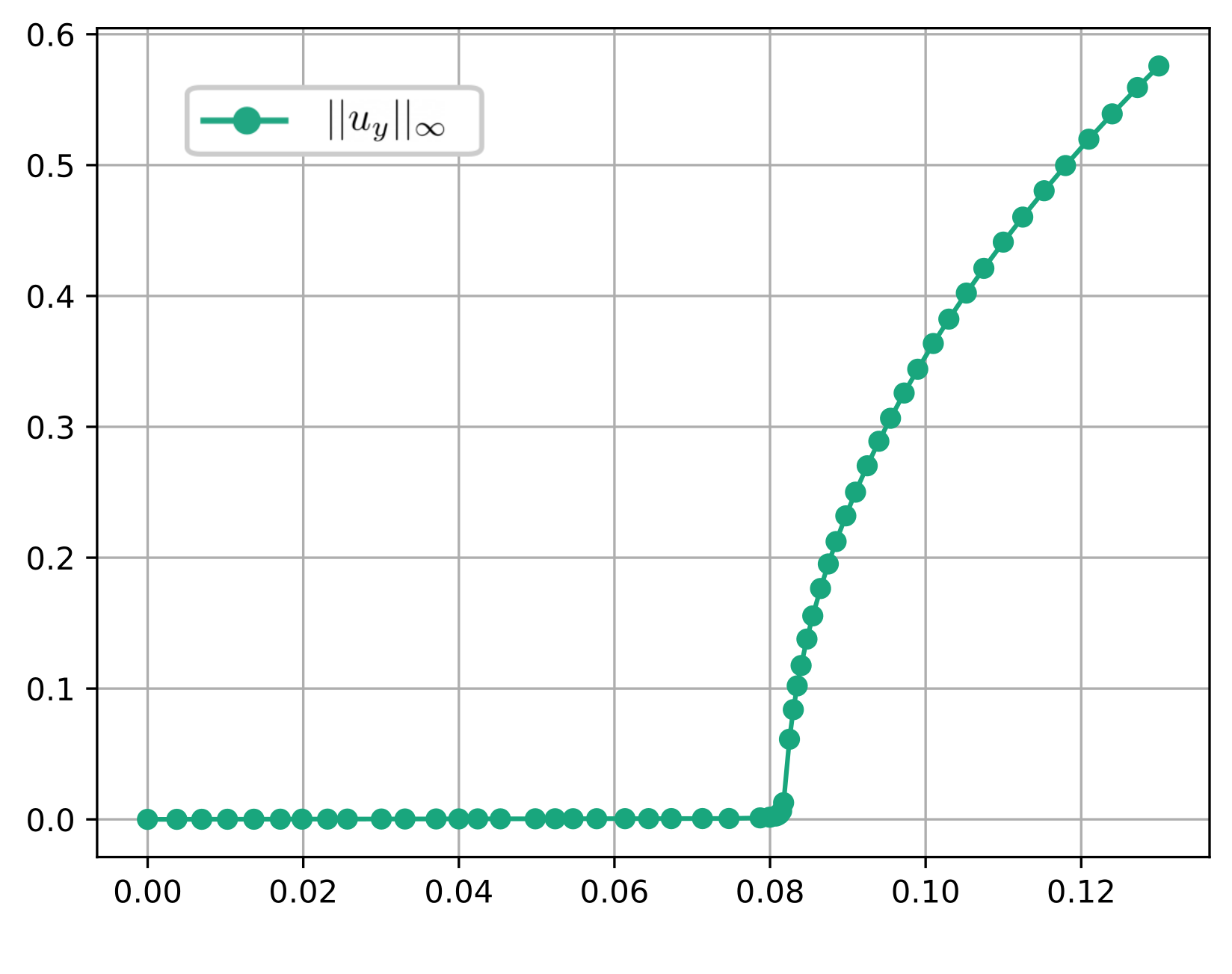}
\put(-100,4){\makebox(0,0){$\mu$}}
\caption{Reduced bifurcation diagram for the 3D SVK tubular geometry with $B = (0, 0, 0)$.}
\label{fig_20_3d_nor_dir}
\end{figure}

\begin{figure}[h]
\centering
\includegraphics[width=0.5\textwidth]{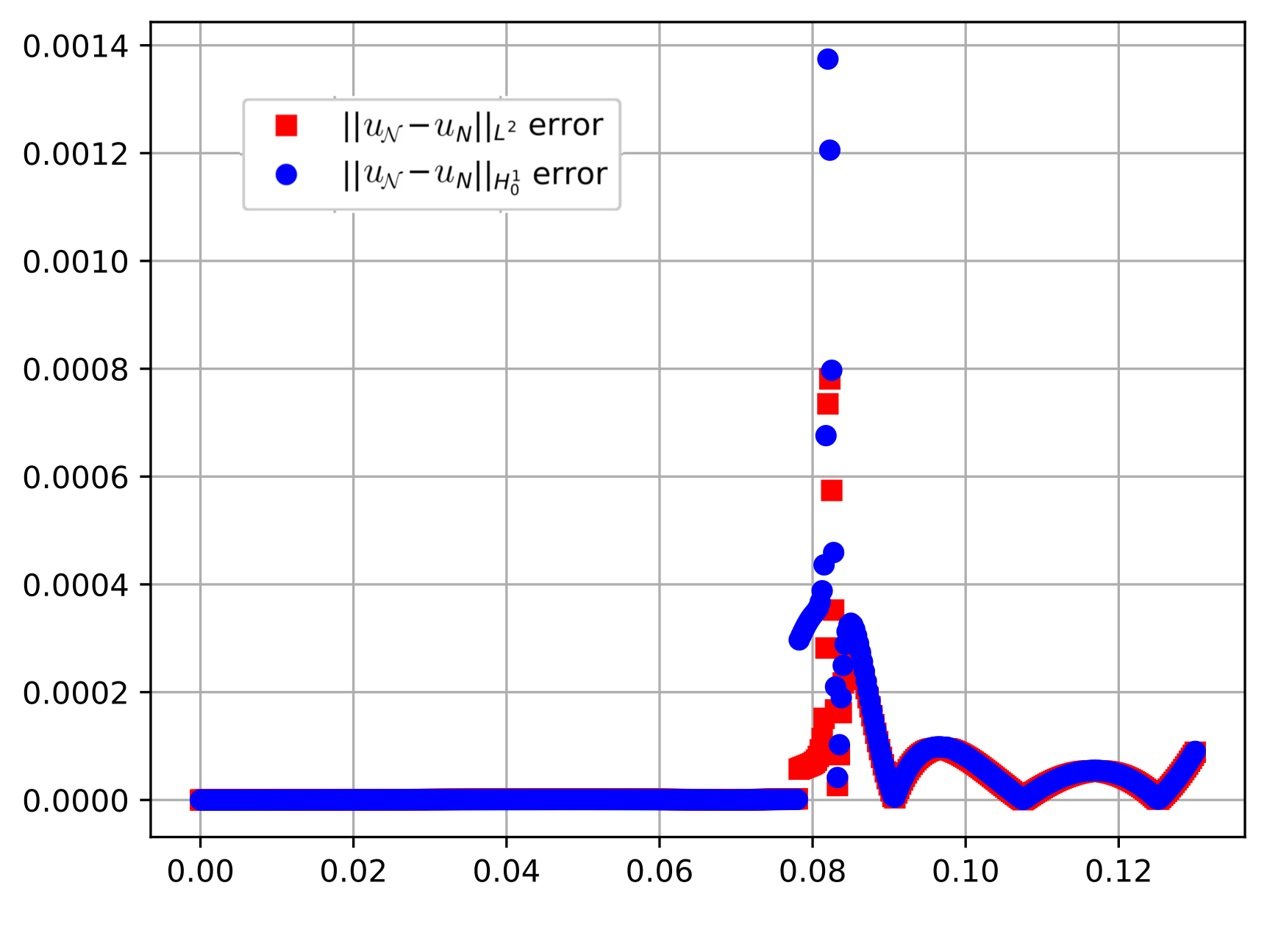}
\put(-95,4){\makebox(0,0){$\mu$}}
\caption{Reduced basis error for the 3D SVK tubular geometry with $B = (0, 0, 0)$.}
\label{fig_6_3d_nor_dir_err}
\end{figure}

In Figure \ref{fig_20_3d_nor_dir} we can see the bifurcation plot for the SVK tube with trivial body force $B = (0, 0, 0)$. This is recognizable since also here the bifurcation phenomenon is characterized by a sharp gradient in the sensitivity. This has the effect of compromising the RB accuracy in the buckling point, as we can observe from Figure \ref{fig_6_3d_nor_dir_err}. Moreover, we can clearly observe that the buckling of the tube occurs for the value $\mu^* = 0.082$.
Despite this, the reduced order model is able to reconstruct the post-buckling behavior with a maximum error of order $10^{-3}$ and an average one of order $10^{-5}$. As we remarked earlier the huge number of cells produces a slightly better speed-up of order 2, which is far from being a good result but highlights how much the empirical interpolation strategies could be useful.   

\begin{figure}[h]
\centering
\includegraphics[width=0.7\textwidth]{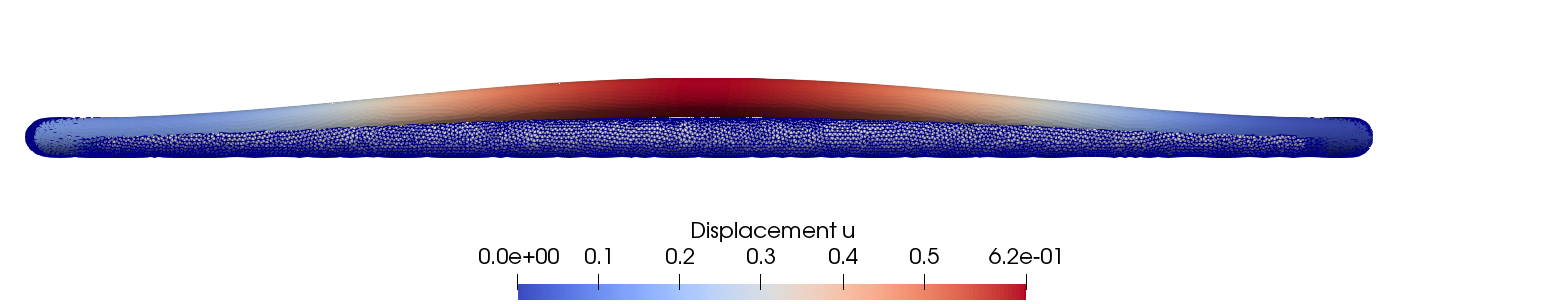}
\caption{High fidelity displacement $u$ for the 3D SVK tubular geometry with $B = (0, 0, 0)$ at $\mu = 0.13$.}
\label{fig_6_3d_nor_dir}
\end{figure}

\begin{figure}[h]
\centering
\includegraphics[width=0.45\textwidth]{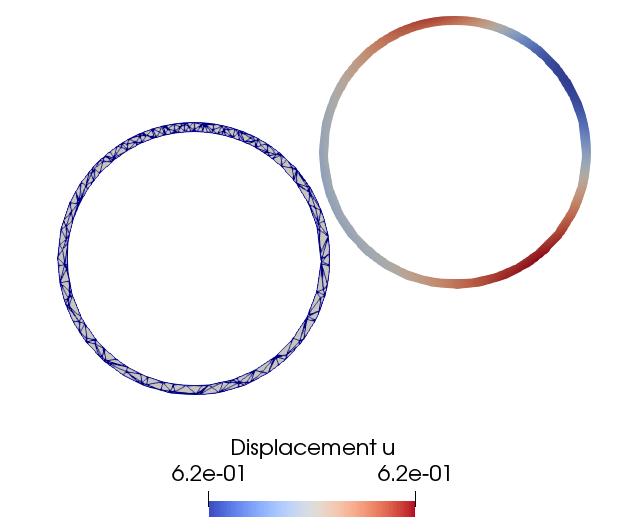}
\caption{Top view of the high fidelity displacement $u$ at $z =10$ for the 3D SVK tubular geometry with $B = (0, 0, 0)$ at $\mu = 0.13$.}
\label{fig_6_3d_nor_dir_ane}
\end{figure}

A representative solution of the post-buckling behavior is depicted in Figure \ref{fig_6_3d_nor_dir} for $\mu = 0.13$ with respect to the original un-deformed configuration (mesh wireframe). 
Moreover, we plot in Figure \ref{fig_6_3d_nor_dir_ane} the top-view sliced with respect to its axis at $z = 10$, where we can observe that the displacement of the cross-section for $\mu = 0.13$ completely exit from its (wireframe) original configuration. 

Now that we have a complete overview of the model, we can finally present the multi-parameter test case with geometrical parametrization of the Dirichlet compressed tubular beam. The aim here is to investigate the buckling behavior of longer geometries.

\begin{figure}[h]
\centering
\includegraphics[width=0.7\textwidth]{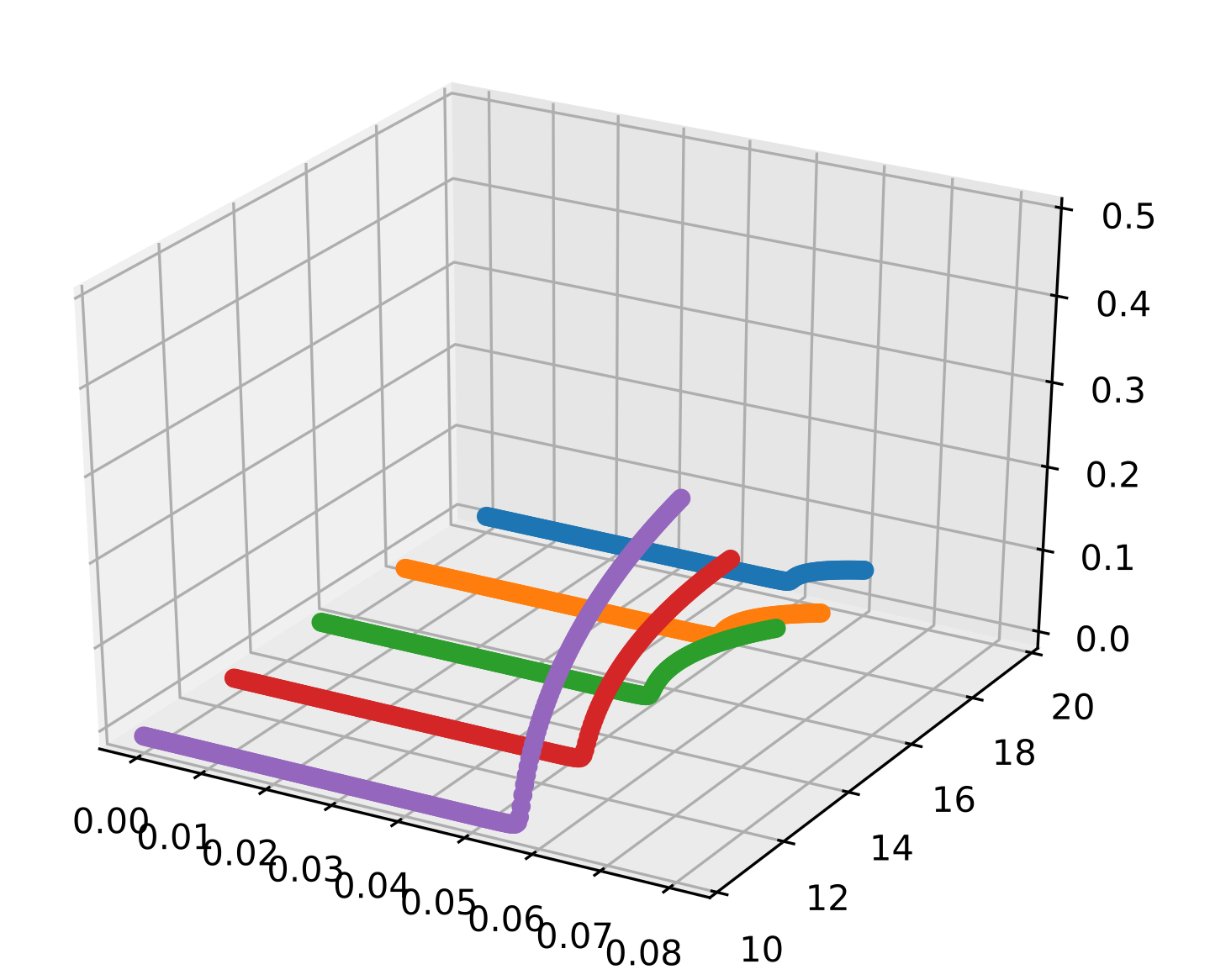}
    \put(-5,115){$\norm{u_y}_{\infty}$}
    \put(-50,28){$\mu_g$}
    \put(-220,12){$\mu$}
\caption{3D bifurcation plot for 3D SVK tubular geometries with $B = (0, 0, 0)$ and $\mu_g \in \Pa_g$.}
\label{fig_17_dir_param_nor}
\end{figure}

The domain is defined, consistently with Section \ref{geom_nor_force}, as $\widetilde{\Omega}(\mu_g) = \widetilde{\Omega}_1 \cup \widetilde{\Omega}_2(\mu_g)$, where $\widetilde{\Omega}_1 = A_r^R \times [0, 20]$ and $\widetilde{\Omega}_2(\mu_g) = A_r^R \times [20, 20 + \mu_g]$ where the geometrical parameter $\mu_g$ varies in $\Pa_g = [10, 20]$.
Keeping fixed the setting as before, we computed a global number of $N_{train} = 800$ snapshots, divided as equispaced points in $\Pa_p$ for each one of the three equispaced values in $\Pa_g$ and we obtained a reduced basis space of dimension $N = 12$ with POD tolerance $\epsilon_{POD} = 10^{-10}$. 
We present the 3-D bifurcation diagram in Figure \ref{fig_17_dir_param_nor}, in which we reconstructed the buckling behavior for 5 equispaced values of the semi-length $\mu_g \in \Pa_g$. From the figure it is evident the effect of the length on the buckling location, indeed we pass from $\mu^* = 0.048$ at $\mu_g = 10$ to $\mu^* = 0.056$ for $\mu_g = 20$.

\begin{figure}[h!]
\centering
    \begin{subfigure}[b]{0.32\textwidth}
         \centering
         \includegraphics[width=\textwidth]{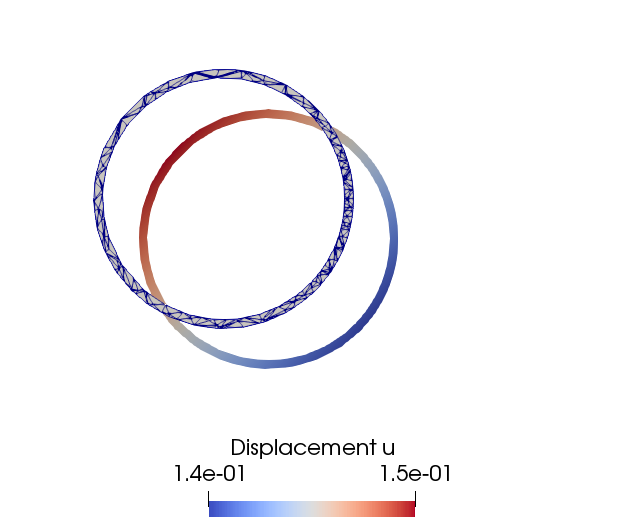}
         \caption{$\mu_g = 10$}
         \label{fig:6_a_geom_1_dir_param_nor}
     \end{subfigure}
     \hfill
    \begin{subfigure}[b]{0.32\textwidth}
         \centering
         \includegraphics[width=\textwidth]{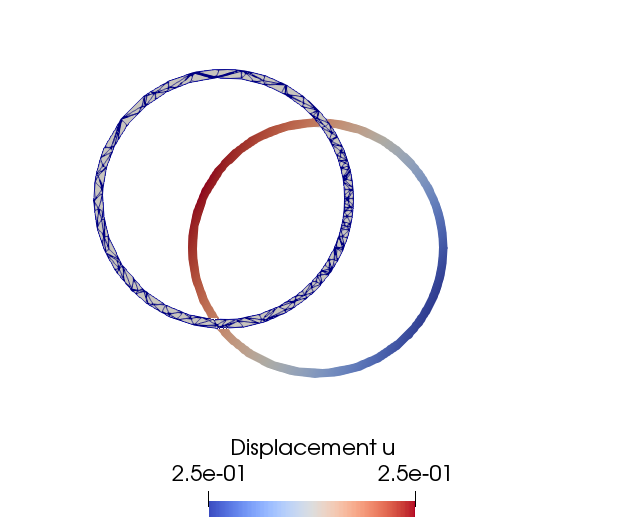}
         \caption{$\mu_g = 12.5$}
         \label{fig:6_b_geom_1_dir_param_nor}
     \end{subfigure}
          \hfill
    \begin{subfigure}[b]{0.32\textwidth}
         \centering
         \includegraphics[width=\textwidth]{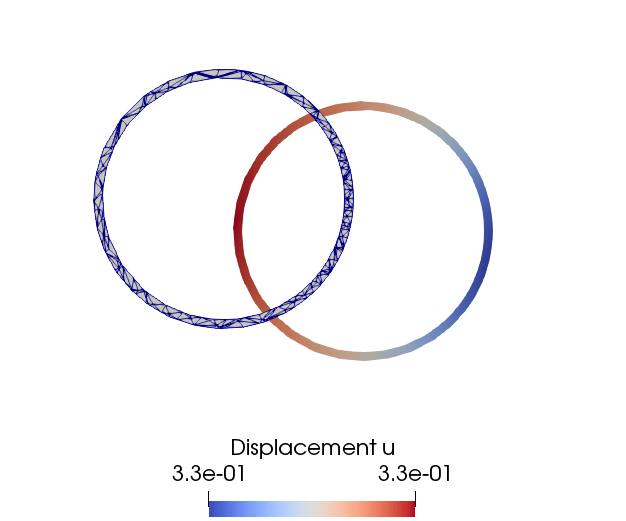}
         \caption{$\mu_g = 15$}
         \label{fig:6_c_geom_1_dir_param_nor}
     \end{subfigure}\\
     
    \begin{subfigure}[b]{0.32\textwidth}
         \centering
         \includegraphics[width=\textwidth]{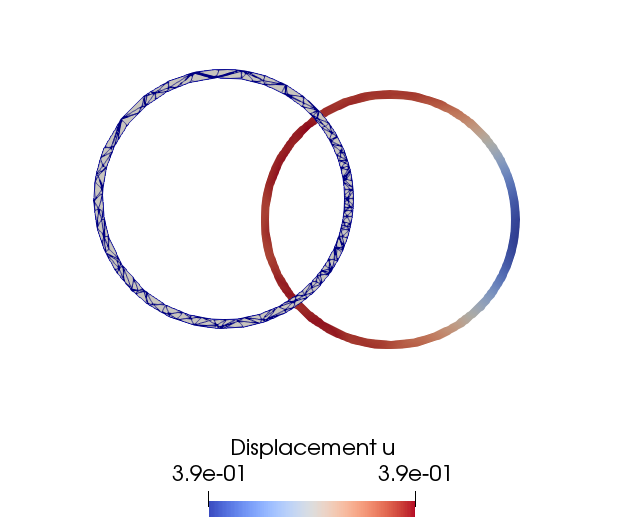}
         \caption{$\mu_g = 17.5$}
         \label{fig:6_d_geom_1_dir_param_nor}
     \end{subfigure}
          \hspace{2cm}
    \begin{subfigure}[b]{0.32\textwidth}
         \centering
         \includegraphics[width=\textwidth]{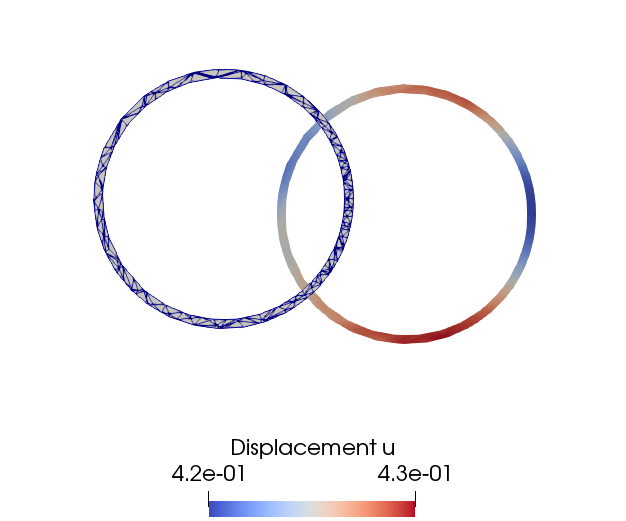}
         \caption{$\mu_g = 20$}
         \label{fig:6_e_geom_1_dir_param_nor}
     \end{subfigure}     
\caption{Top-view of the high fidelity displacement $u$ for the SVK beam with $B = (0, 0, 0)$ for different geometries at $\mu = 0.06$.}
\label{fig_6_geom_dir_param_nor}
\end{figure}

\begin{figure}[h!]
\centering
     \begin{subfigure}[b]{0.49\textwidth}
         \centering
         \includegraphics[width=\textwidth]{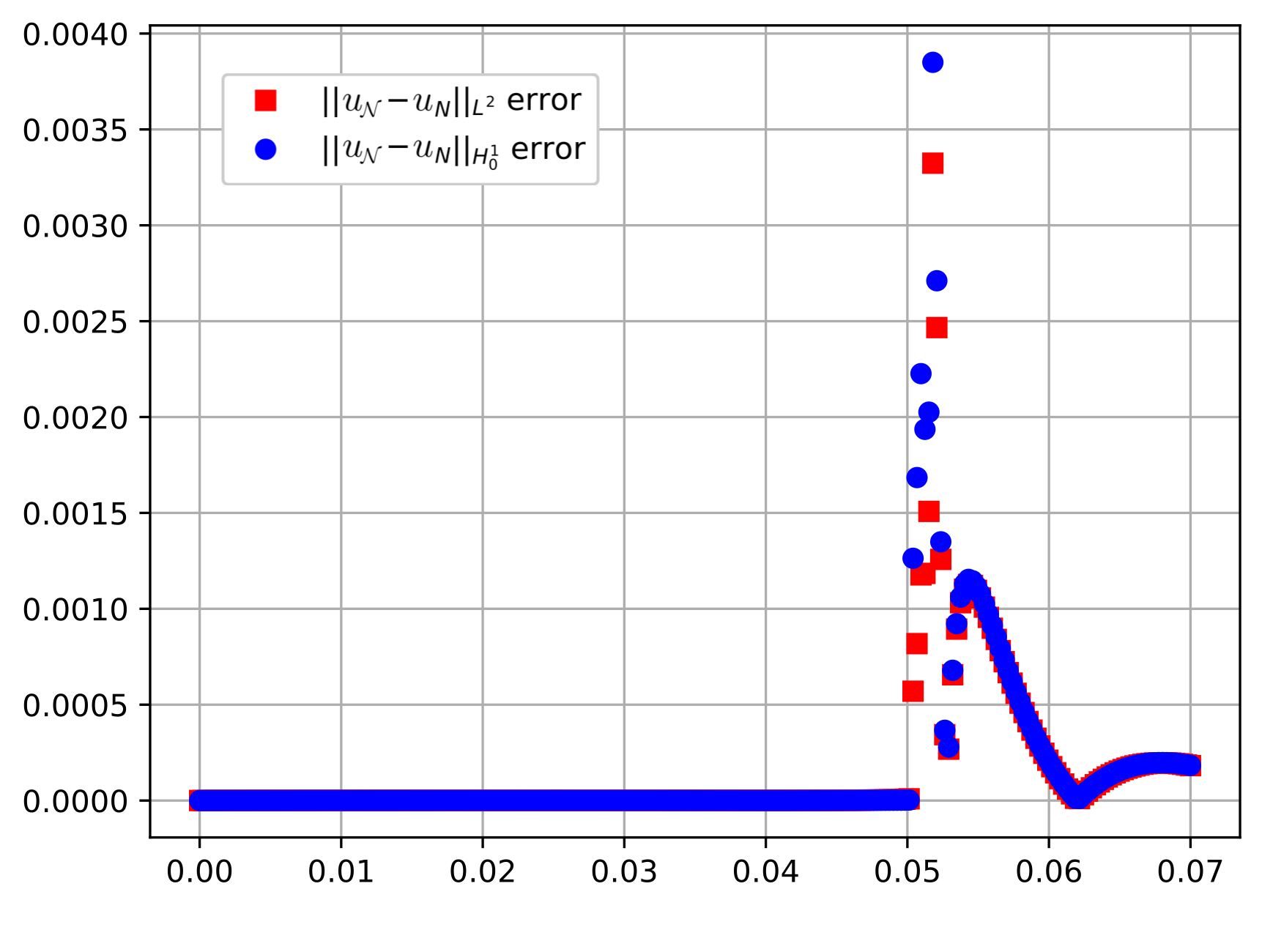}
                           \put(-95,4){\makebox(0,0){$\mu$}}
         \caption{$\mu_g = 15$}
         \label{fig:10_a_geom_dir_param_nor}
     \end{subfigure}
     \hfill
     \begin{subfigure}[b]{0.49\textwidth}
         \centering
         \includegraphics[width=\textwidth]{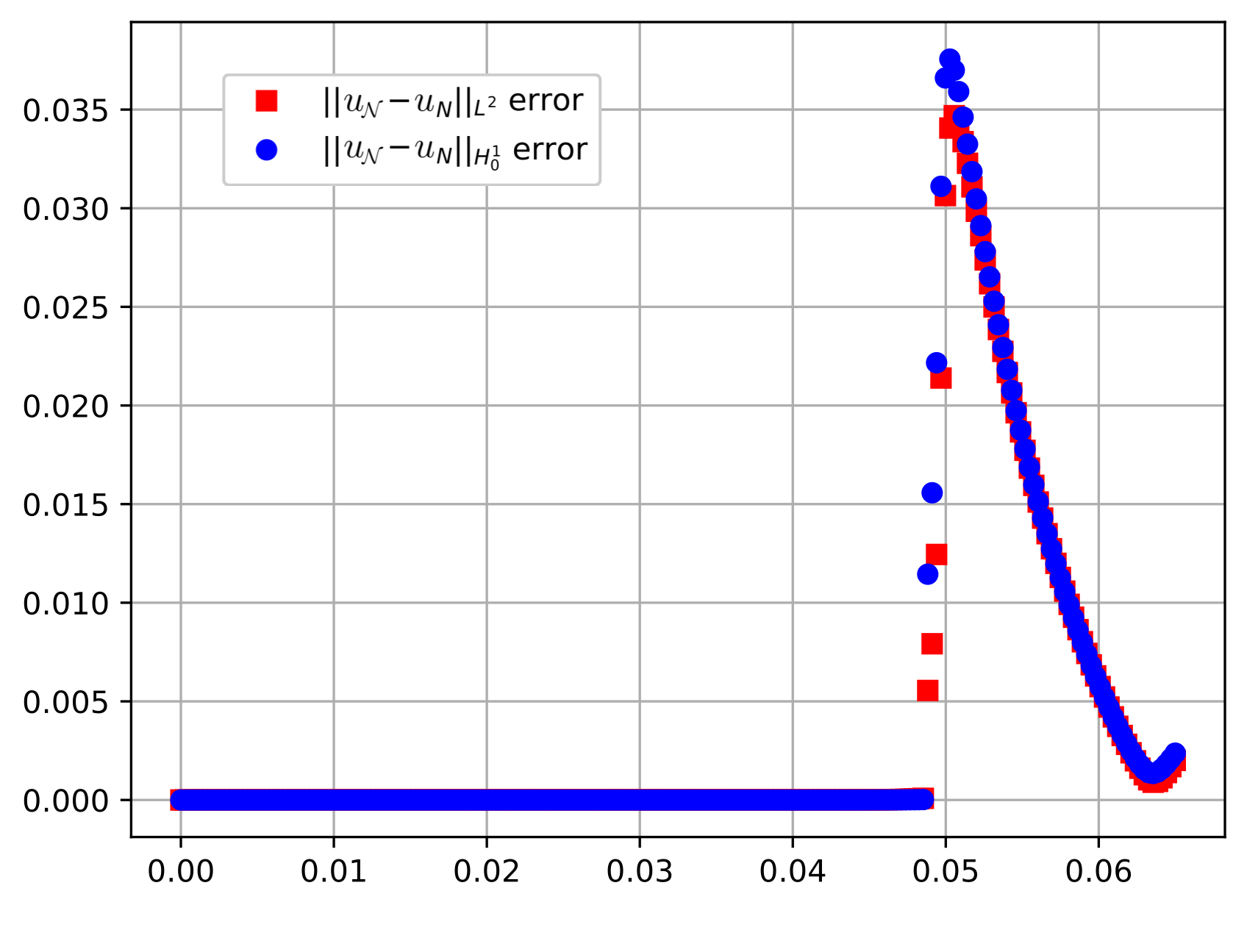}
                           \put(-95,4){\makebox(0,0){$\mu$}}
         \caption{$\mu_g = 17.5$}
         \label{fig:10_b_geom_dir_param_nor}
     \end{subfigure}
	\caption{Reduced Basis errors with respect to $\mu \in \Pa$ for the SVK beam with $B = (0, 0, 0)$ for different geometries.}
\label{fig_10_geom_dir_param_nor}
\end{figure}

The top-view of five representative solutions of the buckling modes for $\mu_g = \left\{10, 12.5, 15, 17.5, 20\right\}$ (colored by their magnitude) are depicted in Figure \ref{fig_6_geom_dir_param_nor} for the same value of the Dirichlet compression $\mu =  0.06$. 
We present in Figure \ref{fig_10_geom_dir_param_nor} the reduced basis error for two values of $\mu_g$, remarking that, also in this more complex context, the reduced manifold was able to approximate with good accuracy the buckling for unsampled values of the geometrical parameter space $\Pa_g$. 

The study of the evolution of the buckling, varying the length of the beam, confirms as expected that longer beams  need smaller forces for the buckling. We highlight that given the importance of the buckling location, already during the offline stage, one can not perform an extensive analysis on the parametrized geometries. Moreover, to prevent non convergence issues due to bigger strains one should consider a much refined mesh, increasing the high fidelity dimension.
Same conclusions on the speed-up hold here, where it increases up to 2.5, mainly due to the computational time $t_{HF} = 51102$(s) required to recover the high fidelity version of Figure \ref{fig_17_dir_param_nor}.

Finally, we remark that, due to the very high number of degrees of freedom within this test case we did not apply any empirical interpolation strategies. Indeed, it would cause an impracticable and too costly offline phase.
Within this setting we analyzed several test cases for the investigation of buckling beams. These allowed us to study how the critical points of the models vary in multi-parameter settings. Moreover, we understood that despite the simpler bifurcating behavior, if interested only in the first buckling mode, unexpected phenomena can occur, and a much deeper investigation is needed.

\section{Conclusions and perspectives}

In this manuscript we have developed and analyzed a reduced order methodology to investigate complex bifurcating phenomena in the context of hyperelastic equations.
In particular, a branch-wise algorithm has been investigated at the reduced level for an efficient reconstruction of the bifurcation diagrams. 
We performed an extensive study of beam-based buckling problems. In particular, the Saint Venant-Kirchhoff and neo-Hookean constitutive relations were analyzed and discussed in connection to different test cases. We studied the forced and the unforced, two-dimensional and three-dimensional rectangular beam subjected to a compression imposed by means of Neumann and Dirichlet boundary conditions. Two multi-parameter cases were studied, varying the physical properties of the beam, such as the Young modulus and the Poisson ratio, and its length through a geometrical parameter. We found that in a multi-parametric scenario, the  evolution of the bifurcation curve is much more complex, and as expected we noticed that longer beams buckle at lower compression.
These investigations allowed us to deal with a more complex benchmark coming from industrial application. A three-dimensional element with annular cross-section was studied, finding unexpected behaviors, investigating the buckling with respect to its length and to the type of compression imposed.

\section{Acknowledgements}
We gratefully thank Prof.\ Anthony T.\ Patera for the inspiring discussions during the visits at MIT, and  Dr.\ Jens Eftang for the fruitful exchanges regarding the test cases. 
We acknowledge the support of MIT-FVG project funded by Regione Friuli-Venezia Giulia ``ROM2S - Reduced Order Methods at MIT and SISSA".
We also acknowledge the support by European Union Funding for Research and Innovation — Horizon 2020 Program — in the framework of European Research Council Executive Agency: H2020 ERC CoG 2015 AROMA-CFD project 681447 "Advanced Reduced Order Methods with Applications in Computational Fluid Dynamics" P.I. Professor Gianluigi Rozza, PRIN 2017 ``Numerical Analysis for Full and Reduced Order Methods for the efficient and accurate solution of complex systems governed by Partial Differential Equations" (NA-FROM-PDEs), and ``GO for IT" program within a CRUI fund for the project ``Reduced order method for nonlinear PDEs enhanced by machine learning".

\bibliographystyle{abbrv}
\bibliography{bib}

\end{document}